\documentclass[twoside]{article}
\usepackage{amssymb, amsmath}
\usepackage[latin9]{inputenc}
\usepackage[francais, english]{babel}
\usepackage[affil-it]{authblk}
\usepackage{graphicx}
\usepackage{pdfpages}
\usepackage{hyperref}
\usepackage{color}
\usepackage[onehalfspacing]{setspace}
\usepackage{ulem}
\makeatletter
\def\imod#1{\allowbreak\mkern10mu({\operator@font mod}\,\,#1)}
\makeatother
\usepackage[babel=true]{csquotes}
\usepackage[T1]{fontenc}
\usepackage{fancyhdr}
\newcommand{\runningtitle}{The algebraic cast of Poincaré's ... }
\usepackage{geometry}
\begin{document}

\title{The algebraic cast of Poincaré's \\ \textit{Méthodes nouvelles de la mécanique céleste}}
\author{Fr{\'e}d{\'e}ric Brechenmacher
\thanks{Electronic address: \texttt{frederic.brechenmacher@euler.univ-artois.fr \\ \textit{Ce travail a  b\'{e}n\'{e}fici\'{e} d'une aide de l'Agence Nationale de la Recherche : projet CaaF\'{E} (ANR-10-JCJC 0101)}}}} 
\affil{Universit{\'e} d'Artois \\ Laboratoire de math{\'e}matiques de Lens (EA 2462) \\
rue Jean Souvraz S.P. 18, 62307 Lens Cedex France.
 \\ \& \\
 {\'E}cole polytechnique \\ D{\'e}partement humanités et sciences sociales \\
91128 Palaiseau Cedex, France. \\}
 \date{}
\maketitle

\begin{abstract}

This paper aims at shedding a new light on the novelty of Poincaré's \textit{Méthodes nouvelles de la mécanique céleste}. The latter's approach to the three-body-problem has often been celebrated as a starting point of chaos theory in relation to the investigation of dynamical systems. Yet, the novelty of Poincaré's strategy can also be analyzed as having been cast out  some specific algebraic practices for manipulating systems of linear equations. As the structure of a cast-iron building  may be less noticeable than its creative façade, the algebraic cast of Poincaré's strategy was broken out of the mold in generating the new methods of celestial mechanics. But as the various components that are mixed in some casting process can still be detected in the resulting alloy, this algebraic cast  points to some collective dimensions of the  \textit{Méthodes nouvelles}. It thus allow to analyze Poincaré's individual creativity in regard with the  collective dimensions of some algebraic cultures. 

At a global scale, Poincaré's strategy is a testimony of the pervading influence of what used to play the role of a shared algebraic culture in the 19th century, i.e., much before  the development of linear algebra as a specific discipline. This shared culture was usually identified  by references to the ``equation to the  secular inequalities in planetary theory.'' This form of identification  highlights the long shadow of the great treatises of mechanics published at the  end of the 18th century. 

At a more local scale, Poincaré's  approach can be analyzed in regard with  the specific evolution that Hermite's algebraic theory of forms impulsed to the culture of the secular equation. Moreover, this papers shows that some  specific  aspects of Poincaré's own creativity result from a process of acculturation of the latter to Jordan's practices of reductions of linear substitutions within the local algebraic culture anchored in Hermite's legacy .
\end{abstract}

\small
\tableofcontents

\newpage
\normalsize
\pagestyle{fancy}
\fancyhead{}
\chead[\runningtitle]{F. BRECHENMACHER}
\fancyfoot[C]{\thepage} 

\section*{Introduction}
\addcontentsline{toc}{section}{Introduction}

\subsection*{What's new in Poincaré's \textit{Méthodes nouvelles} ?}

%Considering the title Henri Poincaré attributed to his three volumes on celestial mechanics (1892-1899), it is quite natural to wonder about 
%What's new in Poincaré's \textit{Méthodes nouvelles de la mécanique céleste} ?
This issue  may seem quite straightforward to modern mathematicians.  Poincaré's \textit{Méthodes nouvelles} have indeed usually been celebrated since the 1950s for laying ground for the development of ``chaos theory''  in relation to the investigation of dynamical systems. Yet, chaos theory has taken various meanings in different times and social spaces over the course of the 20th century.\cite{AubinDahan}  Moreover, in relation to this multifaceted development, various readings of Poincaré's treaties have focused on different aspects of  the \textit{Méthodes nouvelles} :
\begin{itemize}
\item the qualitative investigation of differential equations, which Poincaré had already connected in 1881 to the description of the trajectories of celestial bodies,\footnote{See \cite{Poincar1881f}, \cite{Poincar1882f}, \cite{Poincar1885b},  \cite{Poincar1886h}.}
\item the consideration of the variation of differential systems in function of a parameter,
\item the issue of the global stability of sets of trajectories of celestial bodies, the notion of  ``bifurcation,''
\item the introduction of probabilities into celestial mechanics, 
\item  the recurrence theorem, which states that an isolated mechanical system returns to a state close to its initial state except for a set of trajectories of probability zero.
\end{itemize}
%\begin{quote}
%There is thus is an infinite number of particular unstable solutions  (...) as well as an infinite number  of stable solutions. Is shall add that the former solutions are exceptional (which allows to state that there is stability in general). Yet, the terminology ``exceptional'' does not have any proper signification. Here is what I mean by such a terminology : there is a probability zero for the initial conditions of the motion to correspond to an unstable solution.\cite{Poincar1891}\footnote{Il y a donc une infinité de solutions particulières qui sont instables, au sens que nous venons de donner à ce mot et une infinité d'autres qui sont stables. J'ajouterai que les premières sont exceptionnelles (ce qui permet de dire qu'il a stabilité en général). Voici ce que j'entends par là, car ce mot par lui même n'a aucun sens. Je veux dire qu'il y a une probabilité nulle pour que les conditions initiales du mouvement soient celles qui correspondent à une solution instable.}
%\end{quote}
Various readings  have thus put to the foreground different results, approaches, and  concepts developed in  the monumental three volumes of the \textit{Méthodes nouvelles de la mécanique céleste}.  As a consequence of these retrospective readings, some other aspects have been relegated to the background, including some issues that used to be considered as crucial ones as the time of Poincaré. Among these is the key role played by periodic trajectories in  the strategy the latter developed for tackling the three-body-problem :
\begin{quote}
These are [the trajectories] in which the distances of the bodies are periodic functions of the time ; at some periodic intervals, the bodies thus return to the same relative positions.\cite{Poincar1891} \footnote{Ce sont celles o{\`u} les distances des trois corps sont des fonctions périodiques du temps ; à des intervalles périodiques, les trois corps se retrouvent donc dans les mêmes positions relatives.}
\end{quote}
Periodic solutions were closely associated to the ``novelty'' of Poincaré's approach at the time of the publication of the \textit{Méthodes nouvelles}. At the turn of the 20th century, several astronomers  understood the ``new methods'' as pointing not only to Poincaré's  works but also to the ones of the astronomers  who had made a crucial use of periodic solutions, such as Georges William Hill and Hugo Gyldén.\footnote{This new understanding of the reception of Poincaré's works by astronomers has been communicated by Tatiana Roque at the conference organized for Poincaré's centenary at IMPA, Rio de Janeiro in November 2012.}
Moreover,  Jacques Hadamard, one of  the first mathematicians who  adopted Poincaré's approach to dynamical systems, \footnote{\cite{Hadamard1897}, \cite{Hadamard1901}. See \cite{Chabert1992}}
 attributed the main novelty of the \textit{Méthodes nouvelles}  to the classification of periodic solutions.\cite[p.643]{Hadamard1913}

Periodic trajectories have never been completely forgotten by later commentators of Poincaré. Yet, their significance has often been diminished to the one of an intermediary technical tool for the investigation of dynamical systems. To be sure, the classification of periodic solutions supports the investigation of families of more complex trajectories in their neighborhood, such as  asymptotic solutions, which are curves that asymptotically tend to a periodic solutions with increasing or decreasing time, or the famed doubly-asymptotic (or homoclinic) solutions, which are winding  around periodic solutions. But the role played by periodic solutions in Poincaré's strategy is nevertheless not limited to the one of an intermediary technical tool.  

As shall be seen in this paper, periodic solutions  allow  to introduce linear systems of differential equations with constant coefficients, and thereby to make use of some specific algebraic practices. Poincaré's specific use of periodic solutions is actually intrinsically interlaced with a specific algebraic culture. Moreover, we  claim that this algebraic culture  plays a key model role in the architecture of the strategy Poincaré developed in celestial mechanics.

\subsection*{Hardly new}

Let us  investigate further the issue of the novelty of Poincaré's methods. The following sentence, quoted from the introduction of the first volume of the \textit{Méthodes nouvelles}, exemplifies how Poincaré himself contrasted his approach with the ones of previous works :

\begin{quote}
The investigation of secular inequalities\footnote{As shall be seen in greater details later ``secular inequalities'' designate non-periodic oscillations of the planets on their keplerian orbits.} through a system of linear differential equations with constant coefficients has thus to be considered as rather related to the new methods than to the old ones.\cite[p.2]{Poincar1892} \footnote{L'étude des inégalités séculaires par le moyen d'un système d'équations différentielles linéaires à coefficients constants peut donc être regardée comme se rattachant plutôt aux méthodes nouvelles qu'aux méthodes anciennes.}
\end{quote}

Such a claim for the ``novelty'' of the use of linear systems with constant coefficients may seem quite paradoxical at first sight.  First, the use of such systems in  mechanics  dates back to the great treaties of the 18th century, e.g., the ones of Jean le Rond d'Alembert, Joseph-Louis Lagrange, and Pierre-Simon Laplace.  Second, it is well known that Poincaré's approach has often been celebrated as a starting point for ``chaos theory,'' which has been understood since the 1970s as the science of non linear phenomena. Yet, a similar insistence on the novel role played by linear procedures in the \textit{Méthodes nouvelles}  can also be found in Hadamard's 1913 eulogy of Poincaré. More precisely, Hadamard  presented the novelty of Poincaré's  methods as consisting in returning to some  ancient linear approaches, especially to the criterion of stability of mechanical systems that  Lagrange had stated in his 1788 \textit{Mécanique analytique} by appealing to the nature of the roots of the characteristic equation of a differential system with constant coefficients. 

In the present  paper, we shall thus analyze how some new methods have been cast out  traditional ones in Poincaré's \textit{Méthodes nouvelles}.

\subsection*{Hardly belonging to celestial mechanics}

Modern interpretations  have often  associated Poincaré's qualitative theory of differential equations with a topological approach. It may thus seem very puzzling that  Poincaré himself introduced  in 1881 his qualitative approach in analogy with the role played by Sturm's theorem in algebra :

\begin{quote}
In elementary situations, all the information we are looking for  is,  in general, easily provided through the expression of the unknowns by the usual symbols. [...] But if the question gets more complicated [...] there are two main steps in the reading - as I dare allow myself to say- that is made by the mathematician of the documents in his possession:  the qualitative one and the quantitative one. 

For instance, in order to investigate an algebraic equation, one  starts by looking at  the number of real roots with the help of Sturm's theorem ; which is the qualitative part. Then, one computes the numerical values of the roots, which consists in the quantitative study of the equation. [...] It is naturally by the qualitative part that one has to approach the theory of any function. For this reason, the first problem we shall deal with is the following : \textit{to construct the curves defined by differential equations}. \cite{Poincar1881f} \footnote{Dans les cas élémentaires, l'expression des inconnues par les symboles usuels fournit en général aisément à leur égard tous les renseignement que l'on se propose d'obtenir. [...] Pour peu que la question se complique [...] la lecture, si j'ose m'exprimer ainsi, faite par le mathématicien des documents qu'il possède, comporte deux grandes étapes, l'une que l'on peut appeler qualitative, l'autre quantitative. 
 Ainsi, par exemple, pour étudier une équation algébrique, on commence par rechercher, à l'aide du théorème de Sturm, quel est le nombre des racines réelles : c'est la partie qualitative ; puis on calcule la valeur numérique de ces racines, ce qui constitue l'étude quantitative de l'équation. [...] C'est naturellement par la partie qualitative qu'on doit aborder la théorie de toute fonction et c'est pourquoi le problème qui se présente en premier lieu est le suivant : \textit{Construire les courbes définies par des équations différentielles}.}
\end{quote}

We will see that the two  excerpts  quoted in the previous pages are directly connected one with another. For now, let us  simply remark in passing that the connection between the novelty of the use of linear systems and the qualitative model provided by Sturm's theorem is also highlighted in Hadamard's 1913 eulogy. The latter insisted that Lagrange's approach to differential systems with constant coefficients had been neglected by the predecessors of Poincaré, with the exceptions of Charles Sturm's works on both algebraic and differential equations,  Johann Peter Gustav Lejeune Dirichlet's proof of Lagrange's criterion of stability, as well as the approach later developed by Joseph Liouville.

\subsection*{The algebraic cast  of Poincaré's \textit{Méthodes nouvelles}}

 In a way, this paper  aims at shedding a new light on some aspects of Poincaré's celestial mechanics by highlighting some issues that were hardly new at the time, and which  hardly belonged to celestial mechanics, at least at first sight. More precisely, the approach we are developing in the present paper aims at looking up to  Poincaré's \textit{Méthodes nouvelles}. We shall especially focus on  some of the issues  which have often been considered as secondary issues by the mathematicians who have looked back at Poincaré's approach.

As we shall see, the reference to Sturm's theorem and the importance given to linear systems both implicitly point to an algebraic dimension of Poincaré's works. Even though it  has been overlooked by the historiography, this  \textit{algebraic cast} of the \textit{Méthodes nouvelles} nevertheless plays a key role in the architecture of  Poincaré's treaties.

 Analyzing this algebraic cast is  not only important for grasping the novelty of Poincaré's strategy but also for identifying  some of the temporalities and  collective frameworks in which the latter  took place. As it is used in the present paper, the term ``strategy'' aims at shedding light on   the individual creativity of Poincaré's works in analyzing the latter's  flexible uses of his resources in the constraint frameworks of some social and cultural contexts. 

Poincaré's reference to Lagrange's linear approach to  secular inequalities highlights the necessity to take into consideration some \textit{longue durée} issues, not only in regard with the longstanding  concerns for the stability of the solar system, but also because of the long shadow of the great treaties of mechanics that were published at the turn of the 19th centuries. But we will see also that Poincaré's allusion to the Sturm theorem referred implicitly to some much more recent local mathematical developments.   Along the line of ``scale games,''\cite{Revel1996}  we will therefore make use of different lenses for getting a better view at the various dimensions of a single phenomenon. Both the global and local scales of some specific mathematical cultures are indeed crucial for understanding the  individual originality of Poincaré's own algebraic practices. 

On the one hand, we will  highlight  the strong influence  of what used to be a  shared algebraic culture at the European level during the 19th century. This shared culture  used to be identified  by references to the ``equation to the secular inequalities in planetary theory.'' 

On the other hand, we will  also consider some more local cultures that developed in close connection to this global setting, in a back and forth motion between astronomy, algebra, geometry, analysis, and arithmetic. Among these, we will especially focus  on a specific  approach to Sturm's theorem that  circulated with the legacy of Hermite's ``algebraic theory of forms.'' In so doing, we shall  aim at  shedding a new light on the relationships between celestial mechanics and the other branches of the mathematical sciences in the 19th century.

\subsection*{Mathematical cultures}

As has already been alluded to before,  some new looks at  Poincaré's writings have played a key role  in the emergence of chaos theory in the mid-1970s. Various historical works have  aimed at accounting why such a great burst of activity only took place several decades after Poincaré's death.\cite{Kellert1993}  Yet, David Aubin and Amy Dahan Dalmedico have shown that  this discontinuity  is  mainly the consequence of the retrospective structuration the actors of the development of chaos theory have given  to their own history.
 Poincaré's works had never   been forgotten during the first half of the 20th century even though different aspects  of these works had been  developed  in various contexts.\footnote{See \cite{Dahan1996}, \cite{Mawhin1996}, \cite{Roque2008}, \cite{Roque2011}.}

The problem nevertheless remains of analyzing the collective dimensions of Poincaré's approach to  celestial mechanics. As a matter of fact, this approach  has often been celebrated retrospectively as a point of origin, and thereby for its individuality. In contrast, and as said before, the present paper  proposes a prospective perspective on the \textit{Méthodes nouvelles} in highlighting how Poincaré appealed not only to the century-old works of Lagrange and Laplace, but also to some of the collectives in which these works have been developed over the course of the 19th century. In so doing, we will tackle some issues that are complimentary to the ones considered by some previous historical investigations,\footnote{For an account of the historical studies about Poincaré in various domains, see \cite{Nabonnand2000} and \cite{Nabonnand2005}} such as  Jeremy Gray's wide account of how Poincaré uses the mathematics of his time,\cite{Gray:2000} especially topology, in his qualitative works,\cite{Gray1992}, the relationship between Poincaré's qualitative approach and  traditional research works on differential equations,\cite{Gilain1977} as compared to the approaches of other mathematicians,\cite{Gilain:1991} or in the light of the novelty of the geometric nature of Poincaré's approach,\cite{ChabertDahan1992} \cite{Roque2007} the context of Poincaré's researches on the three-body-problem,\cite{Barrow-Green1994} \cite{Barrow-Greene1996} the relationships between Poincaré and contemporary astronomers,\cite{Nabonnand1999} \cite{Poincar2013} or the influence of the works of Liouville and Ludwig Boltzmann on Poincaré's integral invariants.\cite{Lutzen1984}

Our approach is especially complimentary to the recent researches that have considered  observatories as ``practized places'' of science.\cite{AubinBigg2010} In this framework, Otto H. Sibum, Charlotte Bigg, and David Aubin have introduced the notion of ``observatory techniques''  for designating a  coherent set of physical, methodological, and social techniques rooted in the observatory.  Among these techniques, mathematical procedures figure prominently, whether concerned with astronomy, geodesy, meteorology, physics, or
sociology.  Moreover, an alliance between precise mathematical computations and precise observations lies at the root of the notion of ``observatory culture.'' This full alliance  was especially  demonstrated by the possibility to predict the presence of a missing planet by taking into account anomalies in the orbit of its neighbor. As is illustrated by the instant fame Urbain Le Verrier and John C. Adams acquired when they computed the orbit of Neptune to explain why Uranus was deviating from the orbit, the values of precision in observatory culture generated tremendous optimistic ideals about science, and especially Newton's gravitational theory. This optimism in a culture of precision was not only invested in the precision of the measurements made in the observatory but  also in the precision of the analytical method as represented by the series expansion of  functions. 
 
In addition to shedding new light on  key aspects of the evolutions of sciences in the 19th century,  the notion of ``observatory mathematics'' has also been used by  Aubin for analyzing some  of the collective dimensions of Poincaré's works and their reception.\cite{Aubin2008}  Even though the latter had not yet  been involved with the  observatory at the time  he published his \textit{Méthodes nouvelles}, Poincaré had been trained to the ``observatory techniques'' at the École polytechnique. Aubin has related this acculturation to  observatory culture  to the optimism shown by Poincaré's first approach to the three-body-problem, i.e., to his initial belief in the possibility to prove the stability of the solar system.  According to Aubin,  Poincaré's \textit{Méthodes nouvelles} illustrate that the extreme precision of observatory science provided incentives to re-examine the inner workings of its mathematical technologies.  Poincaré's discovery of homoclinic points has therefore been considered as a product of the ``mathematical culture of the observatory.'' The small reception of this discovery was then connected to the context of the increasing autonomization of both  mathematics and celestial mechanics in regard with one another : homoclinic points may have been considered as not  fundamental enough to modern mathematicians, yet too mathematically rigorous to the observatory community.

In the present paper we discuss some other forms of cultures, ones that  were neither directly  anchored in any specific social, institutional, or practized,  space, nor in any  theoretical or disciplinary framework. These can be identified by the investigation of the intertextual spaces of shared references in which some specific algebraic practices circulated. These ``algebraic practices''  were not limited to some operatory procedures. This terminology actually refers  to some specific intertwining of procedures, representations, meanings, values, and ideals.  Before the structuration of linear algebra as a mathematical discipline in the 1930s, such algebraic practices were usually not embedded within  any  theoretical (or even explicitly reflexive)  framework. In contrast, the circulation of algebraic practices constructed some specific algebraic cultures,  which, in turn, made  practices and cultures evolve in a dynamic interactional system between the collective and the individual.\footnote{The present paper thus appeals to a dynamic notion of culture, in contrast with  the acceptance of this notion in the framework of anthropological structuralism for the purpose of identifying some universal cultural invariants of any human society.\cite[p.XIX]{Strauss1950}} We shall see that the notion of culture is  a key analytical tool for investigating the time-period in which linear algebra did not exist as a mathematical discipline, with its sets of universal core objects  and methods. Indeed, until the 1930s  a  great variety of algebraic practices have circulated  in some specific, yet interlaced, networks of texts.\cite{Brechenmacher:2010a}

The true spaces of such cultures thus lay in the interactions between  texts. Yet, cultures in this sense should not be reified as pointing to actual elements of reality. Texts indeed only interact one with another through the  individuals and groups of individuals who read them. As shall be discussed more precisely in section 3, the interactionist notion of culture has  precisely been developed as a counterpoint to the substantialist conceptions of cultures.\cite{Sapir1949} 

Moreover, because any individual belong simultaneously to   various cultural systems,\cite[p.325]{Strauss1958} 
the interactionist notion of culture allows to take into account the  cultural diversity that is inherent to the unity of any individual's identity. It is thus well adapted to the analysis of individual creativity  in regard with some collective dimensions of mathematics that are not limited to institutions, nations, or local research schools. Yet, we shall see that the shared algebraic culture of the secular equation played no less an important role in Poincaré's approach than the mathematical culture of the observatory, anchored in a practized space ; or than the specific institutionalized technical world of diplomats, scientists and engineers, in which international conventions maps were used by modern states and businesses to control time and space, and to which Peter Galison has related Poincaré's conventionalism.\cite{Galison2003}

\section{Linear systems and periodic solutions}

\subsection{Poincaré's approach to the three-body problem}

According to Newton's law, the planets' mutual attractions disturb the keplerian ellipse that a single planet would run through if it was only subjected to the sun's attraction. These variations may be periodic (the planetary system then returns to  its initial situation). But there may also be some non-periodic  long-term variations in the planet's semi-major axes ; because these oscillations can only be noticeable on astronomical tables that range over a ``century,'' they were designated as ``secular'' variations (or inequalities). These secular inequalities not only make it difficult to compute ephemeris in the long run but they also raise the more theoretical issue of the universality of Newton's gravitational law  :

\begin{quote}
The issue is indeed not limited to the computation of the ephemeris of celestial bodies a few years in advance, for the needs of navigation or for the astronomers to retrieve some already known small planets. Celestial mechanics has a more elevated final goal ; that is to solve the following important question : can Newton's law explain by itself all astronomical phenomena ? [...] The mathematical expression [of this law] is a differential equation which has to be integrated in order to obtain the coordinates of celestial bodies. [...] What will be the motion of $n$ material points attracting each other in a direct ratio of their masses and in an inverse ratio of the square of their distances ? [...] The difficulty begin with a number $n$ of bodies equal to three : the three-body problem has challenged all the efforts of analysts until now.\footnote{Il ne s'agit pas seulement, en effet, de calculer les éphémérides des astres quelques années d'avance pour les besoins de la navigation ou pour que les astronomes puissent retrouver les petites planètes déjà connues. Le but final de la Mécanique céleste est plus élevé ; il s'agit de résoudre cette importante question : la loi de Newton peut-elle expliquer à elle seule tous les phénomènes astronomiques ? [...] elle a pour expression mathématique une équation différentielle, et pour obtenir les coordonnées des astres, il faut intégrer cette équation. [...] Quel sera le mouvement de $n$ points matériels, s'attirant mutuellement en raison directe de leurs masses et en raison inverse du carré des distances ? [...] La difficulté commence si le nombre $n$ des corps est égal à trois : le problème des trois corps a défié jusqu'ici tous les efforts des analystes.} \cite{Poincar1891}
\end{quote}

Moreover, despite their smallness, secular inequalities can accumulate with increasing or decreasing time,  thus producing great changes in the original aspect of the orbits, and thereby threatening the stability of the solar system :

 \begin{quote}
One of the main concerns of researchers is the issue of the stability of the solar system. In truth, such an issue is more a mathematical problem than a physical one. If we discovered a general and rigorous proof, one should nevertheless not conclude that the solar system is eternal. The solar system can indeed not only be subjected to some other forces than Newton's, but the celestial bodies are moreover not reduced to material points. [...] we are not absolutely certain of the absence of a resistant medium ; moreover, tides are absorbing some energy, which is shortly converted into heath by the ocean's viscosity, and which cannot but be borrowed from the celestial bodies' momentum. [...] Yet, all these causes of destructions would act much more slowly than perturbations, and if the latter were not able to alter stability, a much longer lifetime would be guaranted for the solar system.\footnote{Une des questions qui ont le plus préoccupé les chercheurs est celle de la stabilité du système solaire. C'est à vrai dire une question mathématique plutôt que physique. Si l'on découvrait une démonstration générale et rigoureuse, on n'en devrait pas conclure que le système solaire est éternel. Il peut en effet être soumis à d'autres forces que celle de Newton, et les astres ne se réduisent pas à des points matériels. [...] on n'est pas absolument certain qu'il n'existe pas de milieu résistant ; d'autre part les marées absorbent de l'énergie qui est incessamment convertie en chaleur par la viscosité des mers , et cette énergie ne peut être empruntée qu'à la force vive des corps célestes. [...] Mais toutes ces causes de destruction agiraient beaucoup plus lentement que les perturbations, et si ces dernières n'étaient pas capables d'en altérer la stabilité, le système solaire serait assuré d'une existence beaucoup plus longue.}  \cite{Poincar1891}
 \end{quote}

It is not the place here to go into any detail about the context in which Poincaré's works on celestial mechanics developed. Recall that the three-body problem was proposed for the price organized at the occasion of the sixtieth anniversary of  King Oscar II of Sweden and Norway.\cite{Barrow-Greene1996}.  It is well known that Poincaré's prizewinner memoir presented some erroneous conclusions in regard with stability.\cite{Poincar1889} It was in correcting this error that Poincaré introduced the notion of homoclinic trajectories,\cite{Poincar1890}  which has often been considered as the first description of a chaotic behavior.\cite{Anderson1994} A corrected version of the memoir was eventually published in 1890,\cite{Poincar1890} before Poincaré started working on the redaction of his \textit{Méthodes nouvelles}.\cite{Poincar1892} \cite{Poincar1893} \cite{Poincar1899b}

\subsection{Periodic solutions : approximations}

 Poincaré's treaties assimilates a celestial body to a point of coordinates $(x_1,..., x_n)$. The trajectory of such a point    is expressed in function of the time  $t$ by some analytic functions  $X_i$  of the coordinates, which are  the solutions of the following differential system :
\[
(*) \frac{dx_i}{dt}=X_i \ (i=1,..., n)
\]

Because the system (*)  cannot be exactly solved in general, one has  to approximate general solutions by some particular ones :

\begin{quote}
The motion of three bodies depends of their positions and of their initial velocities. If we set the initial conditions of a motion, we  have defined a particular solution. [...] The position and the initial velocity of our satellites could have been such that the Moon would be constantly full ; they could have been such that the Moon would be constantly new. [...] in one of the possible solutions, the new Moon starts to grow but, before it reaches its first quarter, it starts to decrease until  being new  again  and so on ; [the Moon] would then constantly have the shape of a crescent.\footnote{Le mouvement des trois astres dépend en effet de leurs positions et de leurs vitesses initiales. Si l'on se donne ces conditions initiales du mouvement, on aura défini une solution particulière du mouvement. [...] La position et la vitesse initiales de notre satellite auraient pu être telles que la Lune fût constamment pleine ; elles auraient pu être telles que la Lune fût constamment nouvelle [...]  dans une des solutions possibles, la Lune, d'abord nouvelle, commence par croître ; mais, avant d'atteindre le premier quartier, elle se met à décroître pour redevenir nouvelle et ainsi de suite ; elle aura donc constamment la forme d'un croissant.} \cite{Poincar1891}
\end{quote}

While some  particular solutions ``are only interesting because of their strangeness''\footnote{``ne sont intéressantes que par leur bizarrerie''}, others have some ``astronomical applications,'' such as the periodic solutions investigated in Hill's  theory of the moon.\footnote{See \cite{Hill1877} ,\cite{Hill1878}, \cite{Hill1886}.} Yet, even though Poincaré appealed to Hill's approach, the former made it clear that his aim was not to investigate periodic solutions for themselves :

\begin{quote}
Let us consider the example of the three-body-problem [...].\footnote{Here, Poincaré actually focuses on the restricted problem in which the third body,, assumed massless with respect to the other two bodies, cannot disturb the two others, which revolve around their center of mass in circular orbits under the influence of their mutual gravitational attraction. The restricted three-body- problem is then to describe the motion of the third body's trajectory in function of the ratio  $\mu$ of the weights of the two other bodies, which is supposed to be very small.} For  $\mu=0$, the problem is integrable, each of the two small bodies revolves around the third one in a keplerian orbit ; it is then plain to see that an infinity of periodic solutions exist. We will see later that we can conclude that the three-body-problem admits an infinity of periodic solutions, provided that $\mu$ remains small enough.

At first sight, this fact may seem  foreign to any practical interest. The probability is indeed zero that the initial conditions of the motion would be precisely the ones corresponding to a periodic solution. But it may happen that the differences between these initial conditions is very small, and this actually happens in the cases in which the ancient methods fail. One can thus consider with profit a periodic solution as a first approximation, or as an \textit{intermediary} orbit in M. Gyldén's parlance.

Moreover, here is a fact I could not prove rigorously but which seems very likely to me  [...] : we can always find a periodic solutions (which period may be very long), such as the difference between the two solutions is as small as we want. As a matter of fact, the reason why periodic solutions are so precious is that that these solutions are, so to say, the only breach by which we can attempt to storm a fortress which was until now believed to be unassailable.\cite[p.81-82]{Poincar1892}\footnote{Prenons pour exemple le Problème des trois corps [...]. Pour $\mu=0$, le problème est intégrable, chacun des deux petits corps décrivant autour du troisième une ellipse keplérienne ; il est aisé de voir alors qu'il existe une infinité de solutions périodiques. Nous verrons plus loin qu'il est permis d'en conclure que le Problème des trois corps comporte une infinité de solutions périodiques, pouvu que $\mu$ soit suffisamment petit.

Il semble d'abord que ce fait ne puisse être d'aucun intérêt pour la pratique.
En effet, il y a une probabilité nulle pour que les conditions initiales du mouvement soient précisément celles qui correspondent à une solution périodique. Mais il peut arriver qu'elles en diffèrent très peu, et cela a lieu justement dans les cas o{\`u} les méthodes anciennes ne sont plus applicables. On peut alors avec avantage prendre la solution périodique comme première approximation, comme \textit{orbite intérmédiaire}, pour employer le langage de M. Gyldén. 

 Il y a même plus : voici un fait que je n'ai pu démontrer rigoureusement mais qui me paraît pourtant très vraisemblable [...] on peut toujours trouver une solution périodique (dont la période peut, il est vrai, être très longue), telle que la différence entre les deux solutions soit aussi petite qu'on le veut, pendant un temps aussi long qu'on le veut. D'ailleurs, ce qui nous rend ces solutions périodiques si précieuses, c'est qu'elles sont, pour ainsi dire, la seule brêche par o{\`u} nous puissions essayer de pénétrer dans une place jusqu'ici réputée inabordable.} 
\end{quote}

In Poincaré's approach, the classification of periodic solutions is thus not an end in itself. As said  before, these particular solutions aim at approximating more complex trajectories, such as 
asymptotic solutions and homoclinic solutions :

\begin{quote}
I will show [...] how one can take a periodic solution as the starting point of a sequence of successive approximations, and thereby investigate the solutions which are only a little bit different [from the periodic solution].\footnote{ Je montrerai [...] comment on peut prendre une solution périodique comme point de départ d'une série d'approximations successives, et étudier ainsi les solutions qui en diffèrent fort peu.}
\end{quote}

This strategy of  ``approximations'' by periodic solutions  plays  a central role in the \textit{Méthodes nouvelles}.  Before getting into further details about the two types of approximations used by Poincaré (section 4), we shall first analyze the roots of his approach. We shall especially see that Poincaré's strategy of approximation was cast out a specific approach that had been developed  in the 18th century for dealing with secular inequalities.

\subsection{From the small oscillations of  swinging strings to the ones of  periodic trajectories}

In the 18th century, the secular inequalities in planetary theory have been investigated on the model of the mathematization that had been given previously to some problems of swinging strings. In this section we shall provide an overview of this development from the retrospective standpoint of Poincaré's \textit{Méthodes nouvelles}.\footnote{For this reason, we shall not consider in the present paper some  other developments of the 18th century that appealed to a linear approach to mechanical stability, such as the ones related to the stability of ships in navigation.}  We shall thus focus on the few texts Poincaré himself referred to, i.e., mainly  to Lagrange's \textit{Mécanique analytique} and Laplace's \textit{Mécanique céleste}.

Lagrange's approach was  rooted on the one d'Alembert developed in his 1743 \textit{Traité de dynamique}. The latter had investigated the small oscillations $\xi_i(t)$ of a string loaded with two bodies by neglecting the non linear terms in the power series developments of the equations of dynamics. The problem was thus mathematized by  a system of two linear differential equations with constant coefficients. Lagrange generalized this approach in 1766 to the small oscillations $\xi_i(t)$ of a system of $n$ bodies, and thereby to a system of $n$ linear equations :\cite[p.519]{Lagrange1766}
\[
\begin{matrix}
\frac{d^2\xi_1}{dt}= A_{1,1}\xi_1 + A_{1,2} \xi_2 + ... + A_{1,n} \xi_n \\ 
\frac{d^2\xi_2}{dt}= A_{2,1} \xi_1 + A_{2,2}\xi_2 + ... + A_{2,n} \xi_n \\ 
...\\
\frac{d^2\xi_n}{dt}= A_{n,1}\xi_1 + A_{n,2} \xi_2 + ... + A_{n,n} \xi_n \\ 
\end{matrix}
\]
The integration of the  above system was provided by the mathematization of a mechanical observation that had been made by Daniel Bernouilli, according to which the  oscillations of a swinging string loaded with  $n$ bodies can be decomposed into the independent oscillations of  $n$ strings loaded with a single body. The method of integration was thus based on the decomposition of the system into  $n$ independent equations $ \frac{d^2\xi_i}{dt}=\alpha_{i}\xi_j$.  Let  $S$ be the periodicity of such a proper oscillation,\footnote{In modern parlance, a proper oscillation corresponds to an eigenvalue of the matrix $A-SI$. The procedure of integration is thus  tantamount to reducing $A$ to a diagonal form. } $S$ is then the root of the following equation of degree $n$ :
\[
\begin{vmatrix}
A_{1,1}-S & A_{1,2} & ... & A_{1,n} \\ 
A_{2,1}  & A_{2,2}-S & ... & A_{2,n} \\ 
... &...& ... & ... \\ 
A_{n,1} & A_{n,2} & ... & A_{n,n} -S \\ 
\end{vmatrix}
=0
\]
In the 19th century, the above equation was usually designated as the ``equation in $S$''  (including in Poincaré's \textit{Méthodes nouvelles}).\footnote{In modern parlance, this equation corresponds to the characteristic equation of a pair of matrices, $det(A+SI)=0$. Yet, the latter perspective is based on linear algebra, which did not become an autonomous discipline until the 1930s. It thus  introduces some anachronistic conceptions in regard with some collective organizations of knowledge which did not correspond to the object-oriented mathematical disciplines we are used to nowadays. The notion of matrix especially introduces implicitly some anachronistic conceptions in regard with the articulation between objects, representations,  operatory procedures, and the various branches of mathematics.  As a matter of fact, in a given basis of $\mathbb{R}^n$, a matrix $A$ can be understood as representing various objects such as a differential system, a conic, or a quadratic form. But in contrast with  linear-algebra, which is based on structures such as vector spaces, it was usually the recognition of the special nature of the ``equation in $S$'' that supported analogies and permitted transfers of operatory procedures between mechanics, analytical geometry, arithmetic, algebra etc. As we shall see in greater details later, Lagrange's approach was actually based on some polynomial procedures which are very different from the ones of matrix decompositions.} To each distinct root $\alpha_i$ of this equation corresponds a proper oscillation $\xi_i(t)=C_i e^{\sqrt{\alpha_i} t} + C'_i e^{-\sqrt{\alpha_i}t}$. If the equation has $n$ distinct roots, one thus gets $n$ independent solutions through the linear combinations of which one can express all the solutions of the system :
\[
\xi_i(t)=C_1e^{\sqrt{\alpha_1t}}+C'_1e^{-\sqrt{\alpha_1t}}+ C_2e^{\sqrt{\alpha_2t}}+...+C'_ne^{-\sqrt{\alpha_nt}}
\]

In the 1770s, Lagrange and Laplace have transferred this mathematization to the investigation of the ``secular inequalities in planetary theory,'' i.e., to the small oscillations of the planets of the solar system on their orbits.\footnote{See \cite{Lagrange1774}, \cite{Lagrange1778}, \cite{Lagrange1783}, \cite{Lagrange1784}, \cite{Lagrange1788},  \cite{Laplace1775}, \cite{Laplace1776}, \cite{Laplace1787}, \cite{Laplace1789}, \cite{Laplace1799}.} From this point on, the ``equation in $S$'' has thus also been named  ``the equation to the secular inequalities in planetary theory'' (the secular equation for short).\footnote{Quite often, even though not systematically, the terminology ``secular equation'' was used for the case of symmetric systems, whereas the ``equation in $S$'' was used for general pairs of matrices.}

\subsection{From Poincaré to Lagrange and back : the equations of variations}

Let us now get back to the issue of the role played by periodic solutions in Poincaré's \textit{Méthodes nouvelles}. Even though the latter did not consider directly linear differential systems with constant coefficients, his  strategy of approximation by periodic trajectories was nevertheless molded on Lagrange's  method. This strategy indeed aimed at  introducing  linear systems :

\begin{quote}
It is unlikely that in any application, the initial conditions of the motion would be exactly the ones corresponding to a periodic solution ; but it may happen that the difference is very small. If we then consider the coordinates of the three bodies in their real motion, and, on the other hand, the coordinates that the three same bodies would have in the periodic solution, the difference remains very small at least for some  time, and one can thus as a first approximation neglect the square of this difference.\footnote{Il y a peu de chances pour que, dans aucune application, les conditions initiales du mouvement soient exactement celles qui correspondent à une solution périodique ; mais  il peut arriver qu'elles en diffèrent fort peu. Si alors on considère les coordonnées des trois corps dans leur mouvement véritable, et, d'autre part, les coordonnées qu'auraient ces trois mêmes corps dans la solution périodique la différence reste très petite au moins pendant un certain temps et l'on peut, dans une première approximation, négliger le carré de cette différence.} \cite[p.162]{Poincar1892}
\end{quote}

Let $\phi$ be a given periodic solution. Let $x_i(t)=\phi_i(t)+\xi_i(t)$ be a solution close to $\phi$. From
\[
\frac{dx_i}{dt}=X_i \ (i=1,..., n)
\]
 one gets the ``équations aux variations'' that express the difference $\xi_i (t)$ between the coordinates, $x_i(t)$ and $\phi_i (t)$, of the two trajectories,.

Suppose $\xi_i$ is very small and neglect all the terms of a degree higher than the first, one thus gets a system of linear equations with periodic functions of $t$ as coefficients (say of a period $2\pi$) :
\[
\frac{d\xi_i}{dt}=\sum_{j=1,n} \frac{\delta x_i}{\delta x_j} \xi_j \ (i,j=1,..., n)
\]
Now, linearity implies that any solution is a linear combination of $n$ independent solutions $\psi_i(t)$. By periodicity,  $\psi_j(t+2\pi)$ are  thus solutions of the above system also. Thus, $\psi_j(t+2\pi)$ can be expressed as a linear combination of $\psi_i(t)$. One thus eventually gets a linear system with constant coefficients.

The small  variation $\xi_i(t)$ of a periodic trajectory  $\phi_i(t)$ is  thus eventually mathematized by a system that can be integrated by  Lagrange's method. Consider the equation in $S$, and its  roots  $s_i=k_i e^{2\alpha_i \pi}$,  there exists a function $\theta_j$, which is a linear combination of the $\psi_i(t)$, and such that $\theta_j(t+2\pi)=s_j\theta_j(t)$.\footnote{In  modern parlance, $\theta_j$ is an eigenvector associated to the eigenvalue $S_j$.} Moreover, if the equation $S$ has no multiple roots, then
\[
\xi_i(t)=k_1e^{\alpha_1t}\lambda_{1,i}(t)+k_2e^{\alpha_2t}\lambda_{2,i}(t)+...+k_ne^{\alpha_nt}\lambda_{n,1}(t)
\]
with $\lambda_{i,j} (t)$ convergent trigonometric sums of the same periodicity as $\phi_i(t)$.

\subsection{Characteristic exponents}

The coefficients $\alpha_i$ were named ``exposants caractéristiques'' by Poincaré. Let us remark that these exponents interlaced some mechanical and algebraic meanings that were not identical in Lagrange's approach and in Poincaré's one. We have seen that the method of the former was based on a mechanical representation of the roots $\alpha_i$ as the proper periods of the small oscillations of the planets' elliptic orbits. In contrast, Poincaré did not take up Lagrange's \textit{a priori} of linearity. His proper oscillations operate on a given periodic trajectory, thereby generating a set of other trajectories in its neighborhood. Moreover, as we shall see in the next section, the behavior of such a  set of trajectories is controlled by the algebraic nature of the characteristic exponents, especially by the order of multiplicity of the roots of the equation in $S$.

\section{Mechanical stability and algebraic multiplicity }

In Poincaré's \textit{Méthodes nouvelles}, the issue of the ``stability of a periodic solution'' is a first step toward the analysis of  flows of trajectories in the neighborhood of  a given periodic solution. On the one hand, if the periodic solution remains stable, the approached trajectories remain close to it. On the other hand, unstable periodic solutions support the introduction of more complex trajectories, such as  asymptotic solutions which are expressed as power series of the periodic functions $k_i e^{\alpha_i t}$. 

But what, then, does the ``stability of  a periodic solution'' mean ?

\subsection{Stability in the sense of Lagrange}

\begin{quote}
The word stability has been understood with the most different meanings, and the difference between these various meanings will become clear if one recalls the history of Science.\cite[p.140]{Poincar1899b}\footnote{Le mot stabilité a été entendu sous les sens les plus différents, et la différence de ces divers sens deviendra manifeste si l'on se rappelle l'histoire de la Science.}
\end{quote}

Poincaré's works appealed to various notions of stability.\cite{Roque2011} The one that is relevant for our present investigation is what Poincaré designated as ``the stability in the sense of Lagrange'' of a periodic solution :

\begin{quote}
Lagrange proved that when neglecting the square of masses, the orbits' grand axes  remain invariable. He meant to say that with this degree of approximation, the grand axes can be developed in series of terms of the form $Asin(\alpha t+ \beta)$,  with $A$, $\alpha$ and $\beta$ constant.

Thus, if these series are uniformly convergent, the grand axes remain confined within certain limits [...]. Such is the complete stability.

[...] Pushing the approximation further, Poisson later stated that stability prevails if one takes into account the square of the masses but neglect their cubes. But this [stability] did not have the same meaning. [Poisson] meant to say that the grand axes can be developed in series not only with terms of the form $Asin(\alpha t+\beta)$ but also with terms of the form $Atsin(\alpha t+\beta)$.  The value of the grand axis is then subject to continuous oscillations, but nothing proves that the amplitude of these oscillation do not increase indefinitely with time. We can assert that the system will always return, an infinite number of times,  as close as we want to the initial situation  but we cannot assert that the system will not recede greatly. Thus, the word stability  does not have the same sense for Lagrange and Poisson. \footnote{Lagrange a démontré qu'en négligeant les carrés des masses, les grands axes des orbites demeurent invariables. Il voulait dire par là qu'avec ce degré d'approximation les grands axes peuvent se développer en séries dont les termes sont de la forme $Asin(\alpha t+ \beta)$, $A$, $\alpha$ et $\beta$ étant des constantes.

Il en résulte que, si ces séries sont uniformément convergentes, les grands axes demeurent compris entre certaines limites [...]. C'est la stabilité complète. 

Poussant plus loin l'approximation, Poisson a annoncé ensuite que la stabilité subsiste quand on tient compte des carrés des masses et qu'on en néglige les cubes. Mais cela n'avait pas le même sens. Il voulait dire que les grands axes peuvent se développer en série contenant non seulement des termes de la forme $Asin(\alpha t+\beta)$ mais des termes de la forme $Atsin(\alpha t+\beta)$. La valeur du grand axe éprouve alors de continuelles oscillations, mais rien ne prouve que l'amplitude de ces oscillations ne crois pas indéfiniment avec le temps. Nous pouvons affirmer que le système repassera toujours une infinité de fois aussi près qu'on voudra de sa situation initiale mais non qu'il ne s'en éloignera pas beaucoup. Le mot  de stabilité n'a donc pas le même sens pour Lagrange et pour Poisson.}\cite[p.140]{Poincar1899b}
\end{quote}

In the 18th century, Lagrange had generalized a criterion in which d'Alembert had related the stability of a mechanical system to the algebraic nature of the roots of the equation in  $S$ :\cite[p.532]{Lagrange1766}
\begin{enumerate}
\item  The  system is stable if and only if the  $\alpha_i$ are real, negatives and distinct. In this situation all particular solutions have the  form  $sin(\alpha_i t)$ ; their variations are thus  confined  within certain limits.
\item  If an imaginary root, or a real positive root, occurs, the system is unstable.  In this situation, some real exponential oscillations indeed appear in the solutions.
\item If the equation has a multiple root, then the oscillations are unbounded. In that case, it was  believed that the  proper oscillations would take the form $tsin(\alpha_it)$, implying that the amplitudes of variations of the semi-major axis can grow indefinitely with time, so that the system returns an infinite number of times to  its initial configuration, but  also goes far away from it.   This belief would be proven wrong by Karl Weierstrass  in 1858, and independently by Camille Jordan in 1871.\cite{Brechenmacher:2007a}
\end{enumerate}
Lagrange's discussion of  multiple root was modeled on the latter's investigations of a single linear differential equation of order $n$ with constant coefficients. In the latter case, one can indeed also associate to the differential equation an algebraic equation of degree $n$. If this equation has a multiple root of order $k$, a particular solution then takes the form $P(t)sin(\alpha_it)$, with $P$ a polynomial expression of degree $k$. But this situation does not occur in the case of the  systems of $n$ equations generated by mechanical situations. These systems are indeed symmetric. They  can therefore always be reduced to a diagonal form, whatever the multiplicity of the root.\footnote{In modern parlance, for a system to be stable in the sense of Lagrange, it is necessary that such a system can be turned to a diagonal form. Yet, this condition does not require the eigenvalues to be distinct but that  any eigenvalue of multiplicity $k$ is  associated to a vector space of dimension $k$. Such a vector space is generated by $k$  eigenvectors (i.e., proper oscillations) which are linearly independent even though they correspond to the same eigenvalue. To be sure, this situation was difficult to conceptualize in the absence of any formalized notion of vector space.}

 Expressions such as $tsin(\alpha_it)$,  in which ``the time gets out of the sinus'' were usually designated as \textit{secular terms} in the context of celestial mechanics. In the case of the small oscillations of a string, the stability of the system was a given hypothesis of  the mechanical situation investigated. In contrast, in the case of the secular inequalities in planetary theory, one of the main issue at stake was precisely the one of the stability of the solar system. The transfer of the  mathematization of swinging strings to the planets' oscillations thus triggered some discussions on the algebraic nature of the roots of the secular equation. 

Laplace's famous demonstration of the stability of the solar system (in the case of Lagrange's linear approximation) was especially based on a proof that the roots of the secular equations are real. This proof appealed  to the ``very remarkable property'' of symmetry of mechanical differential systems,\cite[p.14]{Hawkins1975} which  had been highlighted  by Lagrange's works on the secular equation.\cite{Brechenmacher:2007b}  Laplace thus concluded on the disappearance of all secular terms. Yet, his proof is not valid in the case of multiple roots which fails Lagrange's criterion of stability.

\subsection{Non linear approaches to stability in the 19th century}

The computations made by Lagrange and Laplace had eventually showed that, up to the first order in the planets' masses, all secular terms vanish. Over the course of the 19th century, the traditional treatment of stability would still consist in trying to eliminate secular terms in order to demonstrate that the variation in the elements of the planets' orbits would be confined within well-determined limits. Yet, in the framework of celestial mechanics, Lagrange's criterion of stability was quickly outdated. Indeed, after Laplace's proof, and starting with Siméon Denis Poisson's works in 1809,\cite{Poisson1809}  the issue of stability has usually been tackled  by taking into  consideration  some  of the  non linear terms in the series development of the coordinates of celestial bodies.

 Poisson's attention to non linear terms gave rise to a new conception of the notion of stability of a mechanical system, that Poincaré designated as ``stability as Poisson periodicity'' :  let $M$ be the point of a given trajectory corresponding to an instant $t$,  this trajectory is stable in the sense of Poisson if the  points in the trajectory of $M$ enter infinitely many times any circle of radius $r$ around $M$ even if $r$ is made arbitrarily small.\cite{Roque2011}  While for Lagrange, stable  solutions must be bounded in the neighborhood of the elliptical orbits, for Poisson, the solutions can go far away from the initial state, but at some time they return to its neighborhood.  

Poisson's approach has been very influential to later developments of celestial mechanics.  In 1856, Urbain Le Verrier  especially proved that non linear terms  in the series developments that depend on a parameter (such as mass, eccentricity, or inclination) can not only provide more precise approximations,\cite{LeVerrier1856}  but can also induce some important alterations of the orbits and can thus threaten the system's stability.\footnote{On Le Verrier's approach to  ``small divisors,'' see \cite{Laskar1992}.}

Over the course of the 19th century,  various attempts have been made to improve Le Verrier's approach in ruling out all secular terms, i.e., in  preventing time from appearing outside trigonometric terms. Various series have been introduced by astronomers such as  Simon Newcomb, Anders Lindstedt, Charles-Eugène Delaunay, Gyldén, Hill, etc.\footnote{See \cite{Tisserand1889}.} 
The notion of stability has thus been increasingly associated to the possibility of getting strictly periodic series development whose first terms decrease fast. It  gave rise to a conception which Poincaré designated as the ``astronomic convergence'' of series in order to distinguish it from the mathematical notion of convergence toward a finite limit. 

As is well known, one of Poincaré's first famous, even though controversial, result in celestial mechanics was to prove the mathematical divergence of the series used by  astronomers :

\begin{quote}
The methods of M. Gyldén, as well as the ones of M. Lindstedt, indeed  provide solely periodical terms, no matter how far the approximation, so that all the elements of the orbits can only oscillate around their mean value. The question would thus be solved if these developments were convergent. We unfortunately know that they are not.\footnote{Les méthodes de M. Gyldén et celles de M. Lindstedt ne donnent en effet, si loin que l'on pousse l'approximation, que des termes périodiques, de sorte que tous les éléments des orbites ne peuvent éprouver que des oscillation autour de leur valeur moyenne. La question serait donc résolue, si ces développement étaient convergents. Nous savons malheureusement qu'il n'en est rien.} \cite{Poincar1891}
\end{quote} 

\subsection{Poincaré's notion of stability of periodic trajectories}

As is documented from his correspondence with Lindstedt from 1883 to 1884,\cite{Poincar2013}
Poincaré tackled the issue of secular terms by getting  back to the linear case and to Lagrange's discussion  on the algebraic nature of the roots of the equation in $S$, i.e.,  the characteristic exponents, which he also designated as ``coefficients of stability'' :

\begin{quote}
There are three conditions to have complete stability in the three-body problem :
\begin{enumerate}
\item That none of the three bodies can recede indefinitely
\item That two of the bodies cannot shock and that the distance between the two bodies cannot decrease below a certain limit
\item That the system returns an infinite number of times as close as we want from its initial situation
\end{enumerate}
[...] There is one case [the one of the restricted three-body problem] in which we have for a long time proven that the first  conditions holds.  We will see that the third condition holds also.  As for the second one, I cannot say. \footnote{Pour qu'il y ait stabilité complète dans le problème des trois corps, il faut trois conditions :
\begin{enumerate}
\item Qu'aucun des trois corps ne puisse s'éloigner indéfiniment
\item Que deux des corps ne puissent se choquer et que la distance de ces deux corps ne puisse descendre au-dessous d'une certaine limite
\item Que le système vienne repasser une infinité de fois aussi près que l'on veut de sa situation initiale.
\end{enumerate}
[...] Il y a un cas [celui du problème restreint] o{\`u}, depuis longtemps, on a démontré que la première condition est remplie. Nous allons voir que la troisième l'est également. Quant à la deuxième, je ne puis rien dire.} \cite[p.343]{Poincar1892} 
\end{quote}

In contrast with  Poisson's stability (criterion 3 in the above enumeration), the stability  in the sense of Lagrange of a given periodic solutions is defined  through the ``equation in $S$'' generated by the equations of variations of the  periodic trajectory. A periodic solution is then said to be stable if all characteristic exponents are distinct purely imaginary numbers. This condition implies that the small variations $\xi_i$ of the periodic solutions will remain finite, since  in this case $\xi_i=(cos(b_i t) + i sin(b_i t))S_{i,k}$, where $S_{i,k}$ are periodic functions.

While in Lagrange's approach, stability used to be a property of an individual trajectory, Poincaré's  notion of stability  concerns the family of other solutions in the neighborhood of a given periodic solution. Yet, the discussions of the 18th century in regard with stability and the algebraic nature of the roots of the secular equation are nevertheless  reproduced almost word for word in the \textit{Méthodes nouvelles} :

\begin{quote}
In sum, $\xi_i$ can in all cases be represented by a convergent series. In this series, the time can enter under the sign sinus or cosinus,  through the exponential $e^{\alpha t}$, or out of the trigonometric or exponential signs.

If all the coefficients of stability are real, negatives, and distinct, the time will only appear under the signs sinus and cosinus and there will be temporary stability.  If one the coefficients is positive or imaginary, the time will appear under an exponential sign ; if two of the coefficients are equal, or if one of them is zero, the time will appear out of the trigonometric or exponential signs.  [...] 

We shall nevertheless not understand the word stability in an absolute sense. We have, indeed, neglected the squares of the $\xi$ [...] . We can express this fact in saying that a periodic solution  has, alternatively to the secular stability, at  least the temporary stability.\footnote{En résumé, $\xi_i$ peut dans tous les cas être représenté par une série toujours convergente. Dans cette série, le temps peut entrer sous le signe sinus ou cosinus, ou par l'exponentielle $e^{\alpha t}$, ou enfin en dehors des signes trigonométriques ou exponentiels. 

Si tous les coefficients de stabilité sont réels, négatifs et distincts, le temps n'apparaîtra que sous les signes sinus et cosinus et il y aura stabilité temporaire. Si l'un des coefficients est positif ou imaginaire, le temps apparaîtra sous un signe exponentiel ; si deux des coefficients sont égaux ou que l'un deux soit nul, le temps apparaît en dehors de tout signe trigonométrique ou exponentiel.   [...] 
Il ne faut pas toutefois entendre ce mot de stabilité au sens absolu. En effet, nous avons négligé les carrés des $\xi$ [...]  
Nous pouvons exprimer ce fait en disant que la solution périodique jouit, sinon de la stabilité séculaire, du moins de la stabilité temporaire.}  \cite[p.343]{Poincar1892}
\end{quote}

Lagrange's approach thus played a model role for the strategy Poincaré developed with the notion of periodic solution. This model role is also exemplified by the latter's discussion of the case of multiple roots. Recall that in the case of multiple roots, Lagrange's method of integration is not valid anymore, because this method is based on the decomposition of the linear system into $n$ independent equations, each associated to a distinct root. For discussing the case of double roots,  both d'Alembert and Lagrange had introduced  an ``infinitesimal variation'' $\xi$ to turn a double root $s_i$ into two distinct roots $s_i$ and $s_i+\xi$. They had  concluded that, if $\xi$ is made to tend toward zero, a particular solution has to take the form  $tsin(\alpha_it)$. Such a reasoning  could be made to fit nowadays criterions of rigor by the use of the Bolzano-Weierstrass theorem on the  set of symmetric matrices (which is bounded and closed, and therefore compact in $M_n(\mathbb{R}$)). Yet, its conclusion is erroneous : as would been shown by  Weierstrass and  Jordan, Lagrange's system can actually always be decomposed into $n$ independent equations because of its symmetric nature. The multiplicity of the roots has thus  no consequence on  stability in the case of symmetric systems. 

But in contrast with Lagrange's approach,  Poincaré's linear systems are not generated from the principle of dynamics but through a linearization of the equations of variations. They  do not have any property of symmetry in general, and therefore cannot be decomposed into $n$ independent equations. As shall be seen in greater details in section 4, Poincaré had developed some very efficient methods for dealing with such issues. These were based on the Jordan canonical form theorem. Yet, in contrast with his great mathematical memoirs of the 1880s, Poincaré did not display explicitly these methods in his works in celestial mechanics.  The \textit{Méthodes nouvelles} initially followed  Lagrange's approach by appealing to an infinitesimal variation for turning a multiple root into distinct roots.\cite[p.67-68]{Poincar1892} Poincaré concluded  that  a root of multiplicity $k$ generates a term $t^k$ out of the trigonometric or exponential functions.  For instance, for a double root $\alpha_1$=$\alpha_2$, two particular solutions are provided by  $\xi_k=e^{\alpha_1 t}\Psi_{1,i}$ and $\xi_i=te^{\alpha_1 t}\Psi_{1,i}+e^{\alpha_1 t}\Psi_{2,i}$. Yet, as shall be seen in the fourth section of this paper, Poincaré's approach to the issue of multiplicity was far to be reduced to this first discussion.

\subsection {Hardly new ...}

Let us end this section with some partial conclusions.  

As has been highlighted in the introduction of the present paper, Poincaré had presented the novelty of his approach in connection with the use of  linear differential systems with constant coefficients.  We  are now able to shed a new light on such a claim, which may have seemed quite paradoxical at first sight.  We have indeed seen the model-role played by Lagrange's approach to secular inequalities for the strategy Poincaré based on periodic solutions and linear systems. 

In a way, the ``new methods'' can  thus be understood as having been cast out  the ancient ones, or more precisely of the very ancient ones as opposed to the ``ancient ones.''  The introduction of the \textit{Méthodes nouvelles} indeed contrasts the novelty of the use of linear systems with the ``ancient methods'' consisting in looking for more and more precise series developments of the coordinates of the celestial bodies :
\begin{quote}
It would be wrong to believe that computing a great number of terms in the [series] developments resulting from  ancient methods  would be enough for computing   ephemeris with a great precision for  a great many years.

These methods, which consist in developing the coordinates of celestial bodies by power [series] of the masses, have indeed a mutual character, which conflict with their use for computing ephemeris in the long run. The series resulting from these methods contain some so called \textit{secular} terms, in which the time gets out the signs sinus and cosinus, and their convergence is thus  doubtful for large values of the    time $t$.

Yet, the presence of these secular terms does not result from the nature of the problem but only from the method at use. It is indeed easy to realize that if the true expression of a coordinate contains a term in
\[
sin \alpha mt
\]
with $\alpha$  constant, and $m$ one of the masses, then one would get the following secular terms when  developing in power series of $m$ :
\[
\alpha mt- \frac{\alpha^3 m^3 t^3}{6}-...
\]
and the presence of these terms would give a very false idea of the true form of the function investigated.

All astronomers have, for a long time,  had a feeling of the point made above ; especially in all the circumstances in which they have aimed at obtaining formulas relevant over a long time ; for instance, in the computation of secular inequalities, the founders of  Celestial mechanics themselves  had to operate differently in renouncing to simply develop along powers of masses.  The investigation of secular inequalities through a system of linear differential equations with constant coefficients must thus be considered as more related to the new methods than to the old ones.\cite[p.2]{Poincar1892}\footnote{Il ne faudrait pas croire que, pour obtenir les éphémérides avec une grande précision pendant un grand nombre d'années, il suffira de calculer un plus grand nombre de termes dans les développements auxquels conduisent les méthodes anciennes.
Ces méthodes, qui consistent à développer les coordonnées des astres suivant les puissances des masses, ont en effet un caractère commun qui s'oppose à leur emploi pour le calcul des éphémérides à longue échéance. Les séries obtenues contiennent des termes dits \textit{séculaires}, où le temps sort des signes sinus et cosinus, et il en résulte que leur convergence pourrait devenir douteuse si l'on donnait à ce temps $t$ une grande valeur.
La présence de ces termes séculaires ne tient pas à la nature du problème, mais seulement à la méthode employée. Il est facile de se rendre compte, en effet, que si la véritable expression d'une coordonnée contient un terme en
\[
sin \alpha mt
\]
$\alpha$ étant une constante et $m$ l'une des masses, on trouvera, quand on voudra développer suivant les puissances de $m$, des termes séculaires
\[
\alpha mt- \frac{\alpha^3 m^3 t^3}{6}-...
\]
et la présence de ces termes donnerait une idée très fausse de la véritable forme de la fonction étudiée.
C'est là un point dont tous les astronomes ont depuis longtemps le sentiment, et les fondateurs de la Mécanique céleste eux-mêmes, dans toutes les circonstances où ils ont voulu obtenir des formules applicables à longue échéance, comme par exemple dans le calcul des inégalités séculaires, ont dû opérer d'une autre manière et renoncer à développer simplement suivant les puissances des masses. L'étude des inégalités séculaires par le moyen d'un système d'équations différentielles linéaires à coefficients constants peut donc être regardée comme se rattachant plutôt aux méthodes nouvelles qu'aux méthodes anciennes.}
\end{quote}

At several occasions, Poincaré  presented his approach to  celestial mechanics in a direct  relationship with Lagrange's \& Laplace's  great treaties of the 18th century. This fact may come as no surprise. The historiography   has  indeed often presented the \textit{Méthodes nouvelles}  as the first treaties to reopen the issue of the stability of the solar system after Laplace's works.  Until now,  the present paper also has mainly investigated the direct relationship between Poincaré and Lagrange. Yet, this relation should not be considered as an exclusive one, and our analysis should not result in presenting Poincaré's new methods as jumping over most of the developments of the 19th century.  

Even though stability may have often been taken for granted after Laplace's proof, Lagrange's methods in celestial mechanics have nevertheless been developed over the course of the 19th century by a number of actors who were working in various domains.  Most of these works were  not directly dealing with  celestial mechanics, which may be the reason why they  have remained invisible to the historiography of the \textit{Méthodes nouvelles}. But despite the fact that the linear approximation underlaying Lagrange's criterion of stability had  been quickly outdated in celestial mechanics, Lagrange's criterion has had a long-standing influence in other branches of the mathematical sciences. 

Lagrange himself had shown that his criterion was tantamount to stating that the equilibrium of a mechanical system is stable if the potential function is in a minimum when the system is in equilibrium. He had then transferred his discussion of the nature of the roots of the equation in $S$  to the investigation of the stability of equilibrium figures, first in the case of the three mutually perpendicular principal axes of a rotating solid body,\cite{Lagrange1775} and later in the case of a rotating mass of fluid. The notion of stability of equilibria is nevertheless different from the one of the solar system : in the latter case, one considers the long-term stability of an individual approximate solution of a perturbed system, while in the former case  stability is a property of solution, i.e., an equilibrium point, that is verified by an analysis of the behavior of other solutions in its neighborhood.\cite{Roque2011} The issue of equilibrium figures of a rotating mass gave rise to a great many developments in the 19th century, among which we may mention the ones of Carl Gustav Jacobi, Liouville, and Bernhard Riemann.\cite{Lutzen1984}

In 1846, Lejeune-Dirichlet eventually reformulated Lagrange's proof, which was based upon linearization, by showing that higher order terms might also correspond to a minimum of the potential function.\cite{Dirichlet1846} Recall that in his 1913 eulogy of Poincaré, Hadamard insisted that one of the main specificity of the former's approach to celestial mechanics had been to take up some aspects of Lagrange's legacy that had almost been forgotten, except by mathematicians such as Dirichlet and Liouville.

As a matter of fact, Poincaré was not the only one to appeal to Lagrange's stability criterion in  the late 19th century. The former's  approach to celestial mechanics  has thus to be analyzed in a broader framework than the one of astronomy. One famous contemporary example is the one of the discussion of the stability of equilibrium of a rotating fluid in the second edition of  Thomson and Tait's treatise (1879-1883).\cite{ThomsonTait} \footnote{On this treatise, see \cite[p.528]{Wise2005} ; on the issue of stability, see \cite{Darigol2002}.}
This discussion was much analogous to Lagrange's approach to the small oscillations of mechanical systems, which the authors reformulated in the framework of Hamiltonian dynamics, i.e., in deriving the equations of motion from the energy principle written in variational form. The stability of equilibrium was thus developed in term of  the minimum of the potential energy function, as well as the maximum of the kinetic function. Yet, the behavior of the system was still determined by means of the roots ot the ``equation of $S$'' : it its roots are all real and negative, the equilibrium is stable.

But the legacy  of Lagrange's secular equation was far from being limited to the circulation of his criterion in connection with issues of equilibrium figures. As shall be seen in the next section, this circulation went with the one of some specific algebraic practices for dealing with linear systems. This situation is well exemplified by the fact that it was in direct connection with both Dirichlet's new proof of Lagrange's criterion (which became an appendix to the third edition of Lagrange's \textit{Mécanique analytique} in 1853), and with Dirichlet's, Hermite's, and Karl Wilhelm Borchardt's  discussions of Lagrange's procedures as applied to quadratic forms,\footnote{See \cite{Dirichlet1842} \cite{Borchardt1846}, \cite{Hermite1855}, \cite{Hermite1857}.} that Weierstrass investigated closely  Lagrange's claims  in 1858, and eventually showed that the multiplicity of the roots is irrelevant for the issue of stability.\cite{Weierstrass1858} 

As shall be seen in the next section, the algebraic practices generated by the mathematization of small oscillations in the 18th century have played a model role in a great variety of domains. Poincaré's \textit{Méthodes nouvelles} made a crucial use of some of these developments of Lagrange's approach in the long run. It was therefore through the prism of several works published over the course of 19th century that Poincaré  read the great treaties of the 18th century.

 We  thus now have to change the scale of our analysis, in order to take into account the collective dimensions of Poincaré's approach at various scales of spaces and times.

\section{The secular equation in the 19th century}

Over the course of the 19th century, a great many texts referred to the  secular equation. Such a reference even quite often appeared in the titles of such texts. Yet, most of these references showed  little interest for celestial mechanics. As has been showed in \cite{Brechenmacher:2007b}, the term ``secular equation'' was actually used to identify  a shared algebraic culture at the European scale.

\subsection{A network of texts}

In order to avoid the \textit{a priori} use of the retrospective categories of modern mathematical theories, our analysis is based on a careful study of  the ways texts are referring one to another, thereby constituting some networks of texts. Yet, such networks cannot be simply identified as webs of quotations.\cite{Goldstein:1999} Not only do practices of quotations vary in times and spaces but intertextual relations may also be implicit. Our approach to this problem consists in choosing a point of reference - here, Poincaré's \textit{Méthodes nouvelles} - from which a first corpus is built by following  systematically the explicit traces of intertextual relations.  A close readings of the texts involved then  gives access to some more implicit forms of intertextual references.\cite{Brechenmacher:2012c}   Because they provide  a heuristic for the construction of a corpus, and thus a discipline for reading texts,  intertextual investigations permit to identify the collective dimensions of mathematics which are shaped by some circulations of knowledge and practices.

 The graph below  provides a simplified representation of the intertextual relations of the texts referring to the secular equations.

\begin{center}
\includegraphics[scale=0.8]{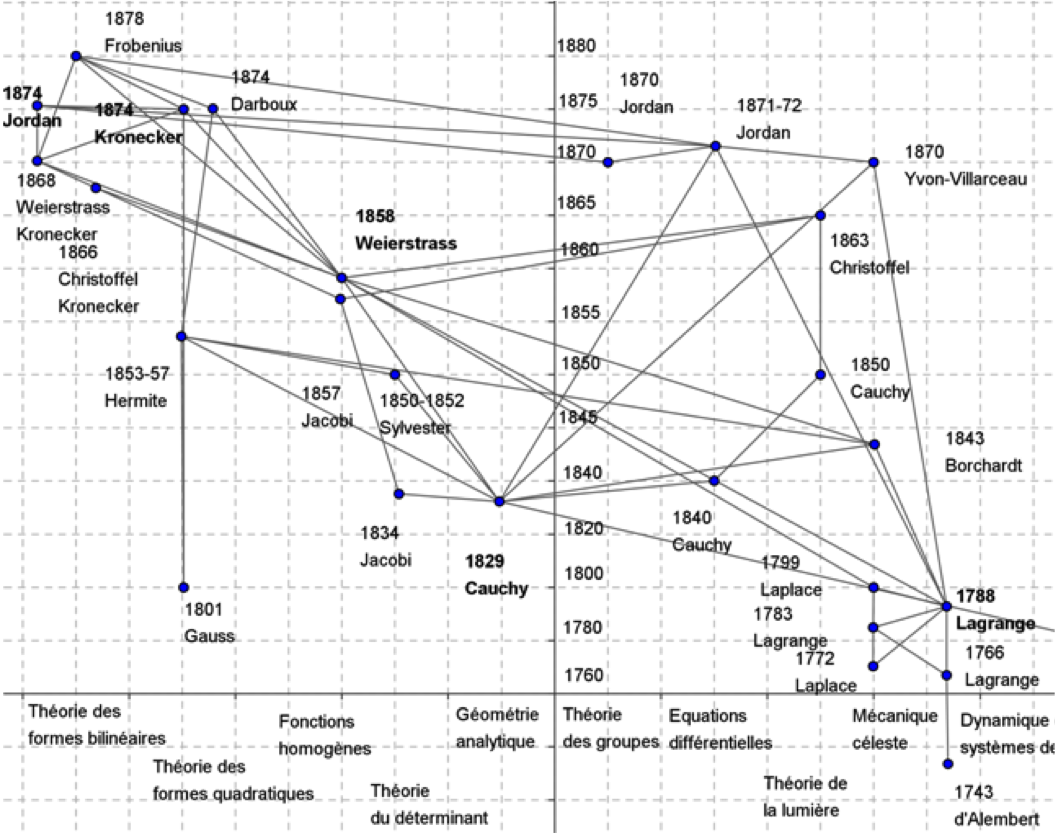}
\end{center}

It is not the place here to provide any detailed description of this graph.\footnote{See  \cite{Brechenmacher:2007b}.} Let us thus limit ourselves to the  following three short comments.

First, the references to the secular equation cannot be considered as aiming at the solution of a specific problem. Even though there was an initial problem, i.e.,  the one of  small oscillations, this problem was considered as having  been solved by Lagrange until Weierstrass and Jordan showed that the former's criterion of stability is erroneous in the case of multiple roots.\footnote{In Paris, the astronomer Antoine Yvon-Villarceau seems to have been the first to question Lagrange's conclusions for multiple roots, more than a century after the latter had stated his criterion of stability.\cite{Villarceau1870} Villarceau's intervention aroused Jordan's interest for this issue, which eventually lead to a contact between the latter's canonical form theorem and Weierstrass's elementary divisors theory.\cite{Weierstrass1868}. See \cite{Brechenmacher:2007b}.}

Second, as is plain to see from the above representation, the intertextual space identified by the secular equation can neither be identified to a theory nor to a discipline. In contrast, the latter equation has supported various analogies over the course of the 19th century between  different branches of the mathematical sciences, such as  dynamics, celestial mechanics,   analytical geometry, the theory of elasticity,  the theory of light, complex analysis, the algebraic theory of  forms,  group theory, etc. Yet, despite the diversity of the theoretical frameworks in which they were working, authors such as Cauchy, Jacobi, Sylvester, Dirichlet, Borchardt, Hermite, Weierstrass, Jordan, Darboux, etc., pointed to a specific algebraic identity laying at the roots of their works. 

Third, even though the  network  of the secular equation can neither be identified by a specific problem nor by a theory, this network shows a very strong  coherence : its texts not only refer frequently one to another  but also to a core of shared references, which corresponds to the main knots in the above graph (e.g., Lagrange 1788, Cauchy 1830, Hermite 1853, Weierstrass 1858).

We shall now see that the network of the secular equation  presents all the characteristics of  a shared culture, in the interactionist sense of the notion of culture.

\subsection{The interactionist notion of culture}

 Let us now  define more precisely the sense attributed to the notion of ``culture'' in the present paper.

 To be sure, the notion of culture has taken different meanings in various times and social spaces.\cite{Beneton1975} 
At the time of Poincaré, ``culture'' was usually considered in France as related to the intellectual development of individuals, in close connection with  the universal notion of ``civilization.'' The latter notion had developed during the 18th century Enlightment, when ``culture'' had been opposed to ``nature'' as a universal, distinctive, character of the human specie. In contrast with the French universalism, ``culture'' was usually understood in a national framework in Germany.  The notion of ``German culture'' had actually emerged at the time of the Napoleonic wars, that is prior to the political unification of Germany as a nation. This notion was initially developed by the intellectual bourgeoisie  in reaction to both the imperialism of French universalism and to the concept of ``civilization'' that was at the core of the court society of the aristocracy.\cite{Elias1939} 

The notion of culture  touches  the core of the symbolic order, i.e.,  to what makes sense. It is therefore no surprise that this notion  has been much debated over the course of  the 20th century. The opposition between the concepts of German culture and French civilization especially reached a climax during World War I, involving the scientific communities of both sides. More generally, this opposition  highlights the longue durée dichotomy between a particularist and a universalist approach to the notion of culture. This dichotomy  has especially been a structuring one for  the  concepts of culture that have been developed in the social sciences. 

In the present paper, the notion of algebraic culture is used in its  particularist sense. In a way, this notion points to the quite traditional meaning of learned cultures. Historically, when the action of ``cultivating'' one's  land or  cattle had been extended in the 17th century to the one of cultivating one's  mind, the term culture was at first always used in addition to a complement, such as the ``culture of the arts'', the ``culture of the sciences,'' etc, thereby identifying some particular forms of learned cultures.  Yet, the present paper appeals to the much more precise definitions and uses of the notion of culture that have been developed in researches in social sciences. 

In this context, focusing rather on particular cultures than on ``the Culture'' has usually aimed at avoiding the ethnocentric bias of applying the researcher's own categories to the field under investigation. 
This methodological principle suits well the purpose we have here to avoid the many anachronisms  that would implicitly come with the uses of the categories of modern linear algebra. Cultural discontinuities are indeed more to be found in time than in space.\cite{Bastide1970}

With the development  of both cultural anthropology and cultural history in North America, especially in the legacy of Franz Boas, specific cultures have  rather come to be identified as systems of interconnected elements than through a list of some distinctive cultural traits.\cite{Malinowski1944} 
The global organization of cultural configurations has thus been considered at least as relevant as their actual content.\cite{Benedict1934}
Yet, these cultural systems should not be reified. Cultures do not have  to  be considered as some existing elements of reality. They  can actually only be accessed through the concrete actions of the individual who create them, transform them,  transmit them, and appropriate them during their whole lives.
 It is  what is shared collectively in terms of behaviors and actions that defines a specific culture. Any culture both presents the relative independence of a collective system in regard with individuals, and the dynamic nature of a system that is always embodied in their lives, and thereby changed by them.

 Cultures, thus, can  be considered as systems of communications, or interactions, between individuals and groups of individuals.\cite{Sapir1949}
Following Edward Sapir's approach, communication has actually been rather conceived as an orchestra  than as a transmitter/receiver type of situation.  The orchestral model points to a situation in which  a collective of individuals play together within a sustainable, yet ever evolving, form of interaction, i.e., their culture.

\subsection{A shared algebraic culture}

The culture of the secular equation is rooted on  a space of intertextual relations, in the sense of the interactions between  various readings of a corpus of texts.  This culture also presents the nature a global system, which is characterized by  specific representations, procedures, ideals, values and organizations of knowledge.

We shall now detail more precisely the specificities of the culture of the secular equation. The first is  a widely shared form of representation : the analytic representation. The second is the more specific use of Lagrange's procedures for manipulating linear systems by the decomposition of the polynomial form of the secular equation. The third is a specific ideal of generality, which is instrumental to the special nature of the secular equation. The fourth is a type of interconnections between branches of the mathematical sciences  through the formal analogies allowed by the secular equation. The fifth is the value attributed to  issues related to multiple roots.

\subsubsection*{The analytic representation}

One of the main issues at stake in the  model of linearization associated to  the secular equation is  the  explicit analytic representations this model provides to  solutions of differential equations. This issue still plays a key role in the strategy Poincaré developed with his use of periodic trajectories in celestial mechanics.

 This situation can  be analyzed as a part of a large scale phenomenon, i.e.,  the crucial role played  by polynomial representations of functions in the long run of the 18th and 19th centuries (with extensions to infinite sums or products). It is well known that such a conception of functions has been challenged in the 19th century, especially in connection to the issues raised by representations by  Fourier's series, from which Georg Cantor's set theory would emerge  in the 1870s. Yet, analytic representations continued to play an important role even after the introduction of the more general notion of functions as applications, as is exemplified by  Poincaré's efforts in the 1880s to provide a  representation of fuchsian functions by infinite sums or products.  It is remarkable that the latter addressed the issue of the analytic representation right from the start of the first volume of the \textit{Méthodes nouvelles}. In the introduction of this volume, Poincaré first indicated that he ``had, as much as possible, complied himself''  to the ``habit'' of expressing  the coordinates of celestial bodies as explicit functions of the time.  This habit, the latter admitted, is  most suitable  for the issue of the computation of ephemeris. But such claims actually aimed a legitimating the frequent use Poincaré also made of ``implicit relations'' between coordinates and time, by resorting to ``integral forms.''  Poincaré indeed argued that the use of such implicit relations is legitimate because these relations allow to deal with the question of the universality of Newton's law.\cite[p.5]{Poincar1892}

 Weierstrass's factorization theorem  is another example of the lasting influence of  analytic representations.   Recall that Weierstrass's theorem states that any analytic function - i.e. the sum of a power series - can be expressed as an infinite product which factors contains the zeros of the function considered. This theorem   illustrates that analytic representations are not limited to a form of notation. They actually cannot  be dissociated from some specific algebraic procedures modeled on the factorization of polynomial expressions.\footnote{The factorization theorem also highlights the limitations of analytic representations. It was indeed in attempting to generalize Weierstrass's theorem to infinite products of rational expressions that G\"osta Mittag-Leffler was drawn to Cantor's set theory. In the case of functions with singular points, one can provide some global analytical representations only in some specific cases while,  in general, one has to consider a function as an application between two sets of points. As Cantor wrote to Mittag-Leffler in 1882 : ``In your approach, as well as in the path that Weierstrass is following in his lecture, you cannot access to any general concept because you are dependent of analytical representations.'' (cited in\cite{Turner2012})}

The analytic  representation plays a key role in all the various lines of developments that emerged from the shared culture of  the secular equation, including in  substitutions group theory, a context in  which it had remained unnoticed until recently.\footnote{See \cite{Brechenmacher:2011}, \cite{Brechenmacher2013b}.} As shall be seen in greater details in section 4,  Jordan's approach to the  reduction of the analytic representation a linear substitution to its canonical form especially played a key role in the development of Poincaré's own algebraic practices in the early 1880s.  

\subsubsection*{Lagrange's procedures for manipulating linear systems}

We have   seen in section 2, that, at the turn of the 20th century, Lagrange's procedures still irrigated in depth Poincaré's \textit{Méthodes nouvelles}.  Yet, Lagrange's procedures are very different from the ones of manipulations of matrices which are familiar to all modern mathematicians. Actually, one of the reasons why Poincaré's  algebraic practices have been overlooked by the historiography  may be that  the procedures of manipulation of linear systems with constant coefficients may have seemed an elementary issue to  commentators in the 20th century.  But in contrast with modern linear algebraic methods, Lagrange's procedures cannot be dissociated from  analytic representations. They appeal to the following polynomial expressions of the coordinates $(x_i^{\alpha_j})$ of the solutions of  symmetric linear systems of $n$ equations with constant coefficients, 
\[
x_i^{\alpha_j}=\frac{\Delta_{1i}}{\frac{\Delta}{S-\alpha_j}}(\alpha_j)
\]
which involve :
\begin{itemize}
\item   $\Delta (S)$, the (polynomial) characteristic determinant of the system $A$, i.e. the one that generates the equation in $S$, $\Delta=det(A-SI)$, 
\item the (polynomial) successive minors  $\Delta_{1i} (S)$, obtained by developing the first line and i\up{th} column of $\Delta (S)$
\end{itemize}

In modern parlance, $x_i^{s_j}$ are the coordinates of the eigenvector associated to the eigenvalue  $\alpha_j$. Lagrange's procedure is thus tantamount to giving a polynomial expression to the eigenvector of a symmetric matrix $A$, as provided by the non-zero column of the cofactor matrix of  $A-SI$. For instance, given
\[ A=
\begin{vmatrix}
1 & - 1 & 0 \\
-1 &  2 & 1 \\
0 &  1 & 1 \\
\end{vmatrix}
\]
of characteristic equation 
\[
det(A-SI)=\Delta (S)=S(3-S)(1-S)
\]
Then, $\Delta_{11} (S)=(1-S)(2-S)-1$,  $\Delta_{12} (S)=(1-S)$,  $\Delta_{13} (S)=1$, e.g., for the eigenvalue $s_1=1$, the coordinates of an eigenvector are :
\[
x_1^{s_1}=\frac{\Delta_{11}}{\frac{\Delta}{s-1}}(1) =\frac{1}{2}
\ , \ x_2^{s_1}=\frac{\Delta_{12}}{\frac{\Delta}{s-1}}(1) =0
\ , \ x_3^{s_1}=\frac{\Delta_{13}}{\frac{\Delta}{s-1}}(1) =-\frac{1}{2}
\]

\subsubsection*{A specific ideal of generality}

The specificity of the algebraic culture of the secular equation was not limited to the technicity of some  polynomial procedures. These procedures also supported a specific ideal of generality. Indeed, the secular equation made it both compulsory and legitimate to deal with  $n$ variables, even in geometric issues. This ideal of generality, in turn, participated to the special nature of the secular equation : even though this  equation is of degree $n$ and thus cannot be solved by radicals in general,  the real nature of its roots can be demonstrated by appealing to the symmetry of the  linear system which generates the equation. 

Let us  exemplify this situation  with a paper published by Augustin Louis Cauchy in 1829. This paper was entitled ``Sur l'équation à l'aide de laquelle on détermine les inégalités séculaires des planètes.'' Yet,  Cauchy did not develop any concern for celestial mechanics. The latter actually appealed to the secular equation for legitimating the generalization to $n$ variables of some methods he had developed for the determination of the principal axis of conic curves and quadric surfaces in the framework of his teaching at the École polytechnique, as well as in his works on the ellipsoids of  the theory of elasticity
in the legacy of Augustin Fresnel's approach to the double refraction of light.\cite{Dahan1980b}

In the general case of $n$ variables, the problem of the determination of the principal axis of a surface of the second degree is tantamount to transforming  the following homogeneous function  (with real coefficients) into a sum of squares :
\[
f(x_1, x_2,...,x_n)=A_{11}x_{12}+A_{22}x_{22}+...+A_{nn}x_{n2}+2A_{12}x_1x_2+2A_{13}x_1x_3 + ... 
\]
Cauchy pointed out that this problem involves considering an equation of degree $n$, which he recognized as ``analogous'' to the secular equation. Cauchy thus mixed up Lagrange's procedures with his own methods, especially the ones he had developed in connection with determinants since 1815. The polynomial expressions involved in Lagrange's formulas were then considered as the changes of coordinates which allow to turn the initial equation of the surface into the following sum of squares :
\[
f(x_1, x_2,...,x_n)=\Delta_{n-1}X_{1}^2+\frac{\Delta_{n-2}}{\Delta_{n-1}}X_{2}^2+...+\frac{\Delta}{\Delta_{1}}X_{n}^2
\]
In modern parlance,  returning  to the example developed above, one can associate to the matrix
\[ A=
\begin{vmatrix}
1 & - 1 & 0 \\
-1 &  2 & 1 \\
0 &  1 & 1 \\
\end{vmatrix}
\]
the quadratic form
\[
A(x_1, x_2, x_3)=x_1^2-2x_1x_2+x_2^2+2x_2x_3+x_3^2
\]
Lagrange's expressions 
\[
x_i^{\alpha_j}=\frac{\Delta_{1i}}{\frac{\Delta}{S-\alpha_j}}(\alpha_j)
\]
then provide the change of basis to get :
\[
A(x_1, x_2, x_3)==1.X_1^2+0.X_2^2+3.X_3^2
\]

\subsubsection*{Specific interconnections between various branches of the mathematical sciences }

The special nature of the secular equation has supported   analogies between various branches of the mathematical sciences.  Among these were the analytic geometry of conics, quadrics, and more general ellipsoids, in connection with  the theories of light and of  elasticity, Charles Fourier's approach to heat theory, Sturm's theorem,  Cauchy's complex analysis, Hermite's algebraic theory of quadratic forms,  Cauchy's 1850 molecular theory of light,\cite{Cauchy1850} (which, as  taken up  by Elwin Bruno Christoffel,\cite{Christoffel1864} eventually gave rise to the theory of bilinear forms), and Jordan's group theoretical  approach to  algebraic  forms. 

Let us return to the analogy that was at the root of Cauchy's 1829 memoir, i.e., between the secular equation and the algebraic equations arising in the determination of the principal axes of the rotation of a solid body, or the ones of a quadric surface. In the late 1820s, this analogy had actually been pointed out to Cauchy by Sturm, who had been especially interested in the secular equation in connection with the statement of his theorem on the number of real roots of an algebraic equation.\cite[p.22]{Hawkins1975}

Recall that the Sturm theorem had been stated in the late 1820s in the framework of researches on linear differential equations.\footnote{Sturm's theorem was stated in a memoir submitted to the Academy. This memoir has been reviewed by  François Arago on May 25, 1829. Even though, this memoir  remained unpublished until 1835 (a situation that was not unusual at the time)\cite{Sturm1835}, Sturm had published a summary of his memoir in the \textit{Bulletin de Férussac} in 1829,\cite{Sturm1829a}. On the connection between Sturm's theorem and linear differential equations, see especially \cite{Sturm1829} and \cite{Sturm1836}.}\cite{Sinaceur1991} Sturm himself stressed that his theorem was ``discovered'' through an investigation of Descartes's rule of sign, i.e., the rule that provides an approximation of the number of roots of an algebraic equation by counting the variations of signs in the coefficients of such an equation. This rule had been generalized by Fourier to the resolution of any determined equation, that is both to algebraic and transcendental equations ; it was then considered as giving rise to a ``general notion'' of analysis, that one could apply to the transcendental functions encountered in problems of celestial mechanics, swinging strings, heat theory, waves propagations, etc. In the case of polynomial equations, Fourier combined Descartes's rule with Rolle's method of cascades for providing an upper bound to the number of real roots of such an equation.

After his arrival in Paris in the mid 1820s, Sturm followed closely Fourier's analytic approach to physical problems, and  especially to the systems of linear differential equations arising in connection with celestial mechanics and problems of heat conduction. In this framework, Sturm extended Fourier's upper bound theorem by applying Euclid's algorithm to the sequence constituted by the polynomial function and its successive derivatives, i.e., to what is nowadays designated a the Sturm sequence of a polynomial equation.\cite{Sinaceur1992} 

But Sturm also took up Fourier's epistemological standpoint, according to which this approach gave rise to a general tool of \textit{a priori} analysis, that is a ``qualitative'' analysis of equations analyzed ``by themselves'' (``en elles-mêmes''), whether these  were algebraic, transcendental or differential equations. For the resolution of differential equations specifically, one had to know the ``march and the characteristic properties'' of the integral functions before actually computing them.\footnote{One may recognize here an expression almost identical to the ones used by Galois at the  time.}  Investigations of  variations of signs  formed the deep unity of Sturm's approach to both algebraic and differential equations. It was in this framework that the latter applied his theorem to the problem of demonstrating the reality of the roots of the secular equation in a memoir he submitted to the Académie des sciences de Paris in 1829. 

Sturm's approach was limited to the case of 5 planets.\cite[p.126]{Hawkins1977} He thus considered a symmetric system of 5 differential equations with constant coefficients.  Appealing to Lagrange's polynomial expression, he  computed the expressions, $\Delta_1$, $\Delta_2$,  $\Delta_3$,..., which, in modern parlance, are tantamount to the minors extracted from $\Delta$  by deleting respectively the first row and column, the first two rows and columns, ..., etc. Sturm then concluded that all these functions ``will have all their roots real and unequal'' and that ``the roots of each of these will comprise in their intervals the roots of the preceding function.'' Yet,  his proof fails in  the case in which ``two consecutive functions  have one or more common roots.''\cite[p.317]{Sturm1829}
As Weierstrass  proved it in 1858, this situation actually never occur in the case of the secular equation because of the symmetry of the system.

Cauchy's approach to the reality of the roots of the equation in $S$ is very similar to the one of Sturm : it shows that all the roots of $\Delta =0$, $\Delta_1 =0$, $\Delta_2 =0$, ..., are real, and that if the roots of  $\Delta_1 =0$ are $r_1 \leq r_2 \leq ... \leq  r_{n-1}$, then the roots of $\Delta =0$ are comprised, respectively, between the limits $- \infty , r_1 , r_2 , ... , r_{n-1} , + \infty$.\cite[p.187]{Cauchy:1829} 

In contrast with  Sturm, whose paper was not published by the Academy, Cauchy not only handled the general case of $n$ variables but had his own memoir published immediately. His method was very influential for later developments.\cite{Hawkins1972} Actually, in the 1830s-1840s, Cauchy's 1829 memoir played a role quite similar to the one of a textbook for the acculturation of many European mathematicians to the algebraic culture  of the secular equation. His approach was especially very quickly developed by Jacobi.\cite{Jacobi1834} \footnote{See also \cite{Jacobi1857}.} As shall be discussed later with the case of Cambridge, the links the secular equation provided between mechanics, geometry, and analysis has exerted a strong  fascination, which was instrumental to the circulation of Cauchy's approach in various contexts of teaching of mathematics.

\subsubsection*{Shared values : a focus on the multiplicity of the roots of the secular equation}

A focus on issues related to multiple roots was shared by the great variety of  the works related to the secular equation over the course of the 19th century.  We have already seen that Lagrange's method for integrating linear differential systems resorted to a decomposition into $n$ distinct equations which is not valid in case of multiple roots. Let us now  have a closer look into this problem. 

Multiple roots may be common roots between the determinant $\Delta$ and its minors $\Delta_{1i}$ and therefore turn into $\frac{0}{0}$ both   the  expression
\[
x_i^{\alpha_j}=\frac{\Delta_{1i}}{\frac{\Delta}{S-\alpha_j}}(\alpha_j)
\]
and 
\[
f(x_1, x_2,...,x_n)=\Delta_{n-1}X_{1}^2+\frac{\Delta_{n-2}}{\Delta_{n-1}}X_{2}^2+...+\frac{\Delta}{\Delta_{1}}X_{n}^2
\]
Cauchy is often presented as one of the first mathematicians who blamed the abusive generality of  algebraic formulations which, such as the ones above,  lose any meaning in some singular cases. In contrast with Sturm, Cauchy indeed dealt carefully with the occurrence of multiple roots that limited the validity of the analytic method he had developed in his 1829 memoir. His first approach to this problem was to introduce a specific reasoning in  the case of multiple roots, which he based on some limit considerations.\cite[p.127]{Hawkins1977}

Multiple roots  gave rise to issues as serious as the ones of imaginary roots. Both actually participated to  the development of complex analysis.  It was indeed for overcoming the problems posed by  multiple roots  that Cauchy developed his Residue theory. In the mid 1820s, Cauchy especially investigated the case of the characteristic equation of  a  linear differential equation with constant coefficient of order $n$.\cite{Cauchy1826} \footnote{On the connection between Sturm's theorem and Cauchy's Residue theory, see also the approach developed in \cite{Cauchy1831} on the localization of the roots of an equation. Sturm and Liouville provided in 1836  a new proof  to Cauchy's localization theorem,\cite{Sturm1836b} one that would turn into a classic with Serret's textbooks of algebra.\cite[p.117-131]{Serret1849}} Later on, in the late 1830s, Residue theory allowed him to provide a fully homogeneous and general solution to  systems of $n$ linear differential equations with constant coefficients, whatever the multiplicity of roots.\cite{Cauchy1839a} \cite{Cauchy1839b} Considering the expression $\frac{\Delta_{1i}(S)}{\Delta(S)}$, with $S$ running on the complex plane, if $\alpha_1$, $\alpha_2$, ..., $\alpha_n$ denote the roots of the secular equation, then the solution $y_i(t)$ of the system of linear differential equation, satisfying $y_i(0)=a_i$, is given by :
\[
y_i(t)=\Sigma_{j;k=1}^n a_j Res_{s=\alpha_k}[\frac{\Delta_{1i}(s)}{\Delta(s)}] e^{st}, \ i=1,2,...,n
\]
Given the multiple combinations of particular cases of multiplicity of roots, it may have seemed hopeless to achieve through   algebraic  methods a solution as homogeneous as the one provided by complex analysis. From Cauchy to Leopold Kronecker, several mathematicians actually appealed to the example of the secular equation to blame the generic tendency of algebraic reasonings which pay little attention to singularity.\footnote{See especially \cite[p.404]{Kronecker1874}. On this issue, see \cite[p.122]{Hawkins1977} \cite{Brechenmacher:2008}.}  

Yet, from the 1850s on,  different  homogeneous algebraic approaches to the secular equation have been developed, e.g.,  Hermite's algebraic theory of quadratic forms, James Joseph Sylvester's matrices and minors, Weierstrass's elementary divisors theorem for quadratic and bilinear forms, and Jordan's canonical form in finite groups theory. In 1858, Weierstrass especially showed that multiple roots actually do not interfere with Lagrange's procedures, which is tantamount to proving that the poles of the rational function $\frac{\Delta_{1i}(S)}{\Delta(S)}$ are simple.  Indeed, in the case of the secular equation the  symmetry of the linear system implies that a root of order $k$ of the determinant is  a root of order $k-1$ of its minors. Lagrange's expressions are thus always valid ; powers of $t$ never arise to destroy stability in case of multiple roots.

\subsection{The slow weathering of the  culture of the secular equation}

%That Poincaré cast out some new methods from the traditional culture of the secular equation is  typical of the time period in which the latter was working. Indeed,
 In the 1850s-1860s, the shared character of the culture of the secular equation was  slowly torn apart by some distinct local lines of developments,  which  gave rise to some local algebraic cultures and eventually resulted into a strong structuration of the algebraic practices at use at the turn of the 20th century.

 Yet, these local algebraic cultures did not develop in isolation one from another. In contrast, they all  developed from within the shared framework of the secular equation. This situation can be understood as   resulting from the intertwining of two phenomena. The first is due to the inner tendencies of evolutions which exist within any given  culture. The second results from the  interplay between interactionist cultures, such as the one of the secular equation, and  other forms of mathematical cultures anchored into  social and institutional contexts. Recall that  a network of texts does not have any existence by its own :  the texts are not only read  by individuals but the culture underlained by the interactions between the texts is also interpreted within the spatialized cultures into which individuals have been  educated and into which they are actually living. The permanent interplay between these two forms of cultures played a key role in the slow divergence of  several lines of developments. Historically, contacts between cultures are indeed anterior to the distinction that produces cultural differences.\cite{Balandier1955}
 The historical, social, and institutional contexts in which interactions took place is therefore crucial.
We shall now look more closely into such evolutions, which we analyze as resulting from processes of acculturations.

As has already been alluded to before, Cauchy's approach  played a crucial role for the extension of the culture of the secular equation to the European continent in the 1830s. This extension is marked by the context of the increasing development of the teaching of mathematics, which  played an important role in the increasing autonomization of mathematics in regard with other branches of the mathematical sciences, such as celestial mechanics. This context is especially documented by the foundings of journals revolving around teaching issues.  The evolution of the populations of contributors to these journals over the course of the 19th century shows an increasing proportion of students and  professors in regard with other professionals trained in mathematics.

Let us consider in more details the founding of the \textit{Cambridge mathematical journal} in the late 1840s.  An important proportion of  the papers published in the first issues of this journal were  devoted to an acculturation to Cauchy's approach into the specific local context of the teaching in Cambridge. Even though most of these publications did not present any new mathematic results, these papers were insistently presented as ``original'' contributions by their authors. This insistence highlights a phenomenon to which we shall return in greater details in section 5, i.e., that processes of acculturations play a key role in the dynamic nature of any culture. It was indeed the whole system of interactions underlained by the secular equation, i.e., between mechanics, algebra, differential equations, analytic geometry, etc., that the authors of the \textit{Cambridge mathematical journals} transposed into their own local culture. This acculturation thus had consequences on both cultural systems in contact, i.e. the one of Cambridge and  the one of the secular equation. As a result, the  secular equation was eventually embedded into a new  local algebraic cultural system that would eventually give rise to Sylvester's matrices and minors, \cite{Sylvester1850} \cite{Sylvester1851} as well as Arthur Cayley's famed theory of matrices. \cite{Cayley:1858}

Let us discuss briefly Sylvester's introduction of the notions of matrices and minors in 1850-1851.\cite{Brechenmacher:2006d} On the one hand, because of their   connection to the secular equation, the procedures of extractions of minors  immediately circulated on the continent.  For instance, shortly after he had coined his terminology of minors and matrices, Sylvester summed up his work in a memoir entitled ``Sur l'équation à l'aide de laquelle on détermine les inégalités séculaires des planètes.''\cite{Sylvester1852} This memoir was published  in French in the \textit{Nouvelles annales de mathématiques}, a journal with a broad audience, including especially teachers and students. As a tool for dealing with multiple roots, Sylvester's minors  were shortly incorporated into other lines of developments, such as Hermite's theory of forms, or Riemann's approach to linear differential equations. 

On the other hand, the issues of symbolical algebra that underlained Sylvester's matrices did not circulate on the continent until the 1880s. These pointed to a symbolical algebraic culture anchored in the specific institutional academic framework of Cambridge, i.e., to the  legacy of  the ``British algebraic school'' which had developed in the first third of the 19th century 
 for legitimizing the symbolic operations of differential calculus.\cite{Durand-Richard:1996}
This  algebraic culture shows a spatialized nature. After having remained for decades specific to authors in Cambridge, such as  Cayley, it first circulated in  Oxford in the 1860s,\cite{Smith1861}
as well as in other academic settings in the U.K. and  in the U.S.A. in the 1870s, and eventually on the continent in  the 1880s.\cite{Brechenmacher:2010a} This situation highlights that the non symbolic, i.e. the technical and material, elements of a culture, such as Sylvester's extractions of minors, circulate more easily than symbolic elements,\cite{Barnett1940} which are charged with ideals and values, such as Cayley's operations on matrices.

Moreover,  each context of interaction imposes its rules, conventions, and expectations to individuals. In his continental publications of the 1850s,  Sylvester himself did not present any  of the symbolic algebraic considerations he was  simultaneously developing  in British journals. The plurality of  contexts of interactions  actually accounts for the inner heterogeneity and instability of any culture, and therefore of any individual. 

Before moving on to the conclusions of this section, let us make a short side comment. Since the beginning of the 19th century, some  procedures of iterations of operators had grown from issues related to differential calculus into taking a prominent part into the specific  algebraic culture developed in Cambridge. Issues of iterations were actually at the core of Cayley's  theory of matrices.\cite{Brechenmacher:2006d} The role played by such issues in Hill's approach to the periodic trajectory of celestial bodies is a testimony of the acculturation of a great many American scientists to this specific culture. This approach may in turn have influenced the iterative procedures of Poincaré's surface-of-section method ; yet, this issue would require further historical investigations.

\subsection{Some partial conclusions}

The crucial role played by linear systems in the \textit{Méthodes nouvelles} is a testimony of the longstanding legacy of the  shared algebraic culture of the secular equation. Let us now draw some partial conclusions from this situation.

\subsubsection*{Local and global cultures}

Let  first pause briefly to consider the various scales of mathematical culture we have discussed in this section : an interactionist, shared, algebraic culture, some spatialized mathematical cultures, anchored in institutions or social spaces, as well as some  local lines of developments generated by the acculturation of an interactionist   culture into a spatialized one.

 It may be tempting to describe this situation  in analogy with the subdivision of biological species into subspecies, i.e., in the framework of a hierarchy between  a global culture  and its sub-cultures. Yet these so called sub-cultures are actually the ones that work as cultures per se, i.e., as the systems of values, ideals, representations, and behaviors by which any group identify itself and act in the surrounding social space.  The prime concept here is thus the culture that stems from immediate interactions, instead of  the more global culture of a large community. 

But such immediate interactions can either take the form of interactions between texts or between actors, thereby giving rise to two different types of local cultures, which are nevertheless always interlaced one with another. Because of this situation, the secular equation can be understood as both a local and global  culture. On the one hand, the interactions underlaying a network of texts are always embodied  into a specific social and institutional setting, as has been exemplified with the context of Cambridge.  But on the other hand,  these intertextual spaces are also shared at a much larger scale, therefore connecting various local cultures whose interrelations, in turn,  construct a global culture. In the  sense of a global culture, the secular equation provided a flexible shared  model to groups and individuals, which allowed the coexistence of some different ways of thinking and of actions.

It was actually the global cultural system of the secular equation that slowly weathered in the second half of the 19th century, while  its various local interpretations evolved more and more independently one from another. Yet, this slow weathering did not occur homogeneously in time and space. References to the secular equation especially continued to play an important role in textbooks, until a new shared culture slowly took over from the 1930s to the 1950s, i.e., the one of linear algebra. 

Retrospectively, the time period extending from the 1860s to the 1930s can be considered as a period of cultural mutation, i.e., as a discontinuous evolution of forms and structures. Destructuration followed by restructuration is indeed the normal evolution caused to any cultural system by cultural contacts, and  processes of acculturations.\cite{Bastide1956} \cite{Bastide1970b}
 During this time period, one can identify various lines of developments specific to some networks of texts.\cite{Brechenmacher:2010a} To be sure, these various lines interacted one with another, but not to the point of constituting any global algebraic culture. For this reason, one can  find  from the 1860s to the 1950s  the subsistence of some ancient elements of the global culture of the secular equation. References to the secular equation  especially continued to play a key role in establishing connections between  various distinct local  cultures.  Yet, communication was partial : authors participating to different lines of developments were usually able to grasp the issues one another tackled but nevertheless always remained faithful to their own algebraic practices.  For instance, when dealing with linear groups, Félix Klein always appealed to computations of invariants based on Weierstrass's elementary divisors theorem while Poincaré  made a systematic use of the Jordan canonical form of a linear transformation.

\subsubsection*{On Poincaré as a mathematician}

 In the preliminary sections of his treatise, when Poincaré  first  dealt in details with linear systems of differential equations with periodic coefficients (i.e. much before this issue was applied to  the equations of variations),  the ``equation in $S$''  was presented as an elementary method. Poincaré even  insisted on the ``extreme simplicity'' of his results, for which he  referred to  the ``well known'' works of Floquet, Callandreau, Bruns, and Stieltjes, with no more precision.\cite[p.68]{Poincar1892}

This situation sheds a new light on the much debated issue of Poincaré's identity as a mathematician in regard with his contributions to celestial mechanics. The notion of identity is  closely related to the one of culture. In the perspective of the relational approach the present paper is building on, the identity of an individual, or group or individuals, is not something that is consubstantial to a culture but one that is  constructed by  relations.\cite{Barth1969}
Identity is thus an ever changing modality of categorization that is used by individuals and groups to organize their communications and exchanges. In this sense, Poincaré's identity as a mathematician is not absolutely determined by the latter's culture, but rather refers to the significations he developed in the various relational situations in which he evolved. It is thus customary to consider closely  the cultural aspects that Poincaré put to the fore to affirm and maintain his own specific identity as a mathematician. 

We have seen in section 2 that it was by proving the  divergence of the series used by the astronomers, and by distinguishing  ``mathematical convergence'' from ``astronomical convergence,'' that  Poincaré  legitimated that he, as a ``mathematician,'' could develop his own approach to celestial mechanics by returning to methods based on linear systems. Moreover, as is illustrated by a quotation given in section 1, Poincaré  insisted repetitively that the issue of the stability of the solar system has to be regarded as a ``mathematical problem'' because of the many physical phenomena that one is neglecting when investigating stability.

But such claims did not imply that Poincaré resorted to the modern methods of mathematicians of his time. The mathematical nature of the problem was mainly associated to the consideration of long run issues, such as the universality of Newton's law, or the stability of the solar system :
\begin{quote}
The final aim of Celestial mechanics is to solve the grand question which is to know if Newton's law explains by itself all the astronomical phenomena ; the only way to achieve this aim is to make observations as precise as possible, and to compare these afterward to the results of computations. Such computations are limited to   approximations [...]. It is thus useless to expect  a greater precision from  computations than from observations ; but we should not expect less from the former than from the latter.

Therefore, we have, for now, to content ourselves with an approximation that will come to be insufficient in a few centuries [...]. To be sure, the epoch in which we will have to renounce to ancient methods is still very distant from us ; but the theorician has to forestall it, because his works have to precede, usually by a great number of years, the one of the numerical computer.\cite[p.1-2]{Poincar1892} \footnote{
Le but final de la Mécanique céleste est de résoudre cette grande question de savoir si la loi de Newton explique à elle seule tous les phénomènes astronomiques ; le seul moyen d'y parvenir est de faire des observations aussi précises que possible et de les comparer ensuite aux résultats du calcul. Ce calcul ne peut être qu'approximatif [...]. Il est donc inutile de demander au calcul plus de précision qu'aux observations ; mais on ne doit pas non plus lui en demander moins.

Aussi l'approximation dont nous pouvons nous contenter aujourd'hui sera-t-elle insuffisante dans quelques siècles. [...] Cette époque, où l'on sera obligé de renoncer aux méthodes anciennes, est sans doute encore très éloignée ; mais le théoricien est obligé de la devancer, puisque son oeuvre doit précéder, et souvent d'un grand nombre d'années, celle du calculateur numérique.}
\end{quote}

Poincaré thus constructed his identity as a theoretical mathematician in contrast with both the figures of the observer and the computer. In questioning mathematically the optimistic faith of observatory culture into Newton's law, the core value of ``precision'' of this culture was turned into the one of ``rigorous approximations.'' The aim was thus to  achieve :
\begin{quote}
[...] the rare results relative to the Three-body-problem that can be established with the absolute rigors that Mathematics demand. It is this rigor that gives some value to my theorems on periodic, asymptotic, and doubly asymptotic solutions.\cite[p.1-2]{Poincar1892} \footnote{[...] les rares résultats relatifs au Problème des trois Corps, qui peuvent être établis avec la rigueur absolue qu'exigent les Mathématiques. C'est cette rigueur qui seule donne quelque prix à mes théorèmes sur les solutions périodiques, asymptotiques et doublement asymptotiques.}
\end{quote}

Actually, Poincaré's works in celestial mechanics  avoided to make an explicit use of some recent mathematical methods for dealing with the key issue of the case of multiple roots in the equation in $S$, such as the Jordan canonical form theorem. As a matter of fact, it was mainly in reference to Lagrange that Poincaré constructed his specific identity as a mathematician getting involved in celestial mechanics, i.e., through the reference to some works that  did not naturally belong to a mathematical culture as opposed to an astronomical one. The issue of the cultures to which Poincaré belonged as an individual has thus to be considered at a finer scale than through the quite rough opposition between ``mathematicians'' and ``astronomers.''
 
To be sure, Poincaré's own identity as  a ``mathematician'' working on  celestial mechanics is not limited to an identity constructed for legitimizing the latter's specific approach, i.e., as a process of communication with fellow ``astronomers.'' This identity can also  be understood in the framework of Pierre Bourdieu's notion of ``habitus,'' which is at the core of the latter's approach to the anthropological notion of culture. In contrast with the interactional nature of a culture, such as the one of the secular equation, Bourdieu's habitus characterizes a social group in regard with other groups that do not share the same social conditions.  Habitus especially works as the ``materialization of the collective memory that reproduces in  successors the elements acquired by  predecessors.''\cite[p.91]{Bourdieu1980b}\footnote{L'habitus fonctionne comme la matérialisation de la mémoire collective reproduisant dans les successeurs l'acquis des devanciers.}
These acquired elements are often so deeply interiorized by individuals that they do not require consciousness to be effective. They are, especially, ``able, in presence of new situations,  to invent some new means to fulfill ancient functions.''\footnote{Ils sont capables d'inventer en présence de situations nouvelles des moyens nouveaux de remplir des fonctions nouvelles.}  Habitus thus makes it possible for individuals to adopt some anticipatory strategies in order to explore new  spaces in accordance with their own social belonging, i.e., as guided with the schemes resulting from their ``primary experiences'' of education and socialization, which weight is enormous in regard with ulterior experiences.

That Poincaré incorporated the collective memory of the shared algebraic culture of the secular equation during his training as a mathematician in the 1870s is shown by his loose reference to the ``classic'' works of Floquet and  Callandreau. More importantly, Bourdieu's approach allows to understand that the traditional dimension of this algebraic culture is not in contradiction with the fact that Poincaré eventually developed a new  approach to linear systems, one that was much related to the legacy of the secular equation, but which differed from it in its details.

\section{The algebraic cast of Poincaré's new methods}

In order to analyze  Poincaré's own individual creativity, it is  customary to locate precisely the position of the latter  
 in the complex algebraic landscape of the late 19th century.  We shall see in this section that Poincaré's algebraic practices resulted from the contact of two local cultures revolving around the works of   Jordan and Hermite respectively. More precisely,  the works of Poincaré's in the early 1880s show the latter's acculturation to Jordan's approach to linear substitutions  into a  cultural system marked by the legacy of  Hermite's approach to algebraic forms in 1850s-1860s.\footnote{On Hermite's theory of forms, see .\cite[p.391-396]{Goldstein:2007}.}

\subsection{Hermite's algebraic forms and Sturm's theorem}

The legacy of Hermite's specific approach to the secular equation  illuminates the key role Poincaré  attributed  to Sturm's theorem when introducing his qualitative approach to differential equations. 

In the early 1850s, both Sylvester and Hermite  were looking for a purely algebraic proof of  Sturm's theorem.\cite[p.124-132]{Sinaceur1991} The secular equation provided a special model case for their investigations. As all the roots of this equation are real,  the counting of the number of real roots is thus limited to the one of distinct roots. More importantly, in the framework of Cauchy's 1829 paper, one can associate to the secular equation a real quadratic form,
\[
f(x_1, x_2,...,x_n)=A_{11}x_{12}+A_{22}x_{22}+...+A_{nn}x{n2}+2A_{12}x_1x_2+2A_{13}x_1x_3 + ... 
\]
which can be turned into a sum of squares:
\[
f(x_1, x_2,...,x_n)=\Delta_{n-1}x_{1}+\frac{\Delta_{n-2}}{\Delta_{n-1}}X_{2}^2+...+\frac{\Delta}{\Delta_{1}}X_{n}^2
\]
The coefficients of such a sum of square are not uniquely determined. Yet, as was shown by Sylvester, the number of positive and negative signs in the sequence of the coefficients is an invariant of the quadratic form  (i.e., Sylvester's inertia law in modern parlance). Moreover, this invariant  actually provides the number of real distinct roots of the secular  equation. In generalizing this approach to any algebraic equation,  Hermite and Sylvester eventually provided an algebraic proof 
of the Sturm theorem.

The role played by  Sturm's theorem in Hermite's early work was at the root of what the latter designated as the ``algebraic theory of forms.''\footnote{See especially  \cite{Hermite1853}, \cite{Hermite1854}, \cite{Hermite1855}, \cite{Hermite1857}.}  In contrast with the ``arithmetic theory of forms,'' which, in the tradition of Carl Friedrich Gauss's \textit{Disquitiones arithmeticae},  concerns classes of equivalences up to substitutions with integers as coefficients, the ``algebraic theory of forms'' investigates the classes of equivalences up to real substitutions that are relevant for the secular equation.

In the context of the development of the algebraic theory of forms in the 1850s, the traditional issues related to the secular equation eventually resulted into a new definition for the notion of multiple root of any algebraic equation :   a root is of order $p$ if all the minors of order $p-1$ of the invariant $\Delta$ vanish. This definition was later used by Hermite's followers, as is exemplified by the  ``Mémoire sur la theorie algébrique des formes quadratiques''  published by Gaston Darboux in 1874:

\begin{quote}
A multiple roots will thus be considered as a simple root if it does not cancel all the minors [of $\Delta$] of the first order ; as a double root if, cancelling all the minors of the first order, it does not cancel all the minors of  the second order, and so on.\cite{Darboux:1874} \footnote{Ainsi une racine multiple pourra être considérée comme simple si elle n'annule pas tous les mineurs du premier ordre ; comme double si, annulant tous les mineurs du premier ordre, elle n'annule pas tous ceux du second et ainsi de suite.}
\end{quote}

In his \textit{Méthodes nouvelles}, Poincaré transferred this notion of multiplicity from the roots of algebraic equations to  the periodic solutions of differential equations :

\begin{quote}
If the determinant of a linear substitution is zero, as well as all its minors of the first [order], the second [order], etc., of the (p-1)\up{th}  order, the equation in $S$ will have $p$ roots equals to zero.\cite[p.174]{Poincar1892} \footnote{Si le déterminant d'une substitution linéaire est nul, ainsi que tous ses mineurs du premier, du second, etc., du (p-1)\up{e} ordre, l'équation en $S$ aura $p$ racines nulles.}
\end{quote}

Before analyzing further this aspect of Poincaré's approach, it is customary to recall that it was very common  in the 19th century to resort to analogies between algebraic and differential equations. Yet, different forms of such analogies have been developed in various contexts. For instance, in 1913, Hadamard  compared  Poincaré's concerns for sets of trajectories to the investigation of the relationships between algebraic roots in Galois theory.\footnote{This comparison has to be understood in the framework of the great variety of analogies that have been made between Galois theory and differential equations  in the late 19th century.\cite{Archibald:2011}} But we have seen in the introduction of the present paper that both Poincaré and Hadamard also pointed to a very different analogy when they referred to Sturm's theorem as a model of ``qualitative approach.''

At the beginning of the 19th century, the ``theory of equations,'' and more generally ``algebra,'' were usually considered as a ``specie'' of an ``analytic gender'' altogether with ``differential analysis,'' ``infinitesimal analysis,'' ``geometric analysis,'' and the ``analysis of curves.'' These various species of analysis were often crossbred one with another, with little concern for their actual specificity. Such crossbreedings were even theorized by some authors, as is exemplified by the introduction of Fourier's posthumous treatise on equations. 
Sturm's theorem is a typical product of such crossbreeding between algebra and analysis. This theorem was presented by Alfred Serret as ``one of the most brilliant discovery by which the mathematical Analysis  has enriched itself''.\cite{Sinaceur1991}\footnote{L'une des plus brillantes découvertes dont se soit enrichie l'Analyse mathématique.} Sturm himself presented his theorem as exemplifying a principle stated by Fourier : ``the complete resolution of numerical equations [is] [...] one the most important application of differential calculus.''\footnote{La résolution complète des équations numériques [est] [...] une des plus importantes applications du calcul différentiel.} Moreover, Sturm's interest for the determination of the number of real roots of an algebraic equation took place in his ``general analysis'' of differential equations, which, as seen before, consisted in considering first the ``appearance'' and the ``march'' (``l'allure et la marche'') of the integral before trying to compute its numerical values, that is by the ``sole consideration of the differential equations themselves, with no need for their integration.''\footnote{La seule considération des équations différentielles en elles-mêmes, sans qu'on ait besoin de leur intégration.}

We shall now investigate another specific form of  analogy between algebraic and differential equations in Poincaré \textit{Méthodes nouvelles}, which focuses on the distinction between simple roots and multiple roots.

\subsection{Poincaré's  approximations by periodic solutions}

Recall that Poincaré developed a strategy of successive approximations of the trajectories of celestial bodies by periodic trajectories. Moreover, we have seen in section 2 that the investigation of sets of trajectories in the neighborhood of a given periodic solution was based on the notion of stability of periodic solutions. But   stability was precisely determined by the multiplicity of the characteristic exponents, i.e., the roots of the equation in $S$. It is therefore no wonder that the transfer of the notion of multiplicity from roots of algebraic equations to periodic solutions of differential equations plays a key role in Poincaré's work.  Let us now return to Poincaré's strategy of approximations by periodic trajectories, which we shall  analyze more closely in the light of the legacy of Hermite's approach to the secular equation.

For the sake of clarity, we shall distinguish between two distinct meanings in Poincaré's  ``method of approximations'' by periodic trajectories. The first  consists in investigating sets of solutions of the same differential system. This method, which we shall designate as the method of variations,   revolves around the following elementary problem :  given  two periodic solutions with close initial conditions, do these solutions have similar behaviors over time ? We shall designate  the second method as the method of perturbations. It consists in investigating the variations of a  differential system in function of a small parameter $\mu$.  An important part of Poincaré's work is devoted to the proof of the existence of periodic solutions for some given initial conditions and to the analysis of their behavior by perturbation, which implies considering simultaneously the  solutions of distinct differential systems.

The method of perturbation  is legitimated by the delimitation of what Poincaré designated as the   ``restricted three-body problem.'' In this case,   the third body cannot disturb the two others, which revolve around their center of mass in circular orbits under the influence of their mutual gravitational attraction, while the third body, assumed massless with respect to the other two bodies, moves in the plane defined by the two primaries and, while being gravitationally influenced by them, exerts no influence of its own. The restricted problem is then to describe the motion of the third body's trajectory in function of the ratio  $\mu$ of the weights of the two other bodies, which is supposed to be very small :

\begin{quote}
The [restricted three-body] case of the problem is the one in which it is supposed that the motions of the three bodies take place in the same plane, that the weight of the third one is zero, and that the first two bodies revolve in concentric circles around a shared center of gravity. If  $\mu=0$, the situation is very simple. As a matter of fact, $m_1$ is motionless while the motion of $m_3$ is a keplerian ellipse of which $m_1$ is a focus.  What happens if $\mu$ is not zero but  very small ? [...] Do we have the right to conclude [that a system with some periodic solutions for $\mu = 0$]  will still have such solutions for the small values of $\mu$ ? [...] The first periodic solution to have been   pointed out is the one discovered by Lagrange, in which the three bodies run around three similar keplerian ellipses, while their mutual distances remain equal to a constant ratio. [...] M. Hill investigated another [periodic solutions] in his remarkable researches on the theory of the Moon [...]. \cite[p. 106 \& 153]{Poincar1892} \footnote{Ce cas est celui du problème o{\`u} l'on suppose que les trois corps se meuvent dans un même plan, que la masse du troisième est nulle, que les deux premiers décrivent des circonférences concentriques autour de leur centre de gravité commun.
Lorsque $\mu=0$, la situation est très simple. En effet, $m_1$ est alors immobile tandis que $m_3$  décrit une ellipse képlérienne dont la position de $m_1$  est un foyer. Que se passe-t-il lorsque le paramètre $\mu$ n'est pas nul, mais simplement très petit ? [...] Avons-nous le droit [de] conclure [qu'un système admettant des solutions périodiques pour  $\mu = 0$] en admettra encore pour les petites valeurs de $\mu$ ? [...] La première solution périodique qui ait été signalée pour le cas o{\`u} $\mu>0$ est celle qu'a découverte Lagrange et o{\`u} les trois corps décrivent des ellipses képleriennes semblables, pendant que leurs distances mutuelles restent dans un rapport constant. [...] M. Hill, dans ses très remarquables recherches sur la théorie de la Lune en a étudié une autre [...]. }
\end{quote}

The position of the third body (in phase space, in modern parlance, see \cite{Chenciner2007}) is described by two linear and two angular variables, $x_i$ and $y_i$ respectively, $y_i$ being periodic with period $2\pi$, connected by the integral $F(x_1, x_2, y_1, y_2)=C$. The differential equations can then be put down into the following Hamiltonian form :
\[
\frac{dx_i}{xt}=\frac{\delta F}{\delta y_i}, \  \frac{dy_i}{xt}=\frac{\delta F}{\delta x_i}, \ (i=1,2)
\]
which can be considered as defining flows on a three-dimensional surface in the framework of the qualitative approach Poincaré developed from 1882 to 1886.

Poincaré himself did not make a clear distinction between what we have designated above as the methods of variations and perturbations. The first four chapters of the \textit{Méthodes nouvelles} are actually rather structured by the strategy based on periodic trajectories than by a distinction between these two methods. The treaties opens with two first  chapters devoted to classic results regarding the existence of solutions of differential equations, the canonical forms into which various types of equations can be reduced, and complex analysis. The introduction of the equation in $S$ of a linear system with periodic functions as coefficients concludes these preliminary chapters.  The third chapter is devoted to the introduction of the notion of periodic solutions,\cite[p.79]{Poincar1892}  with a focus on the distinction between ``simple'' and ``multiple'' periodic solutions.\cite[p.83]{Poincar1892} What follows is then structured by the various types of situations that may occur in regard with periodic solutions. Poincaré first discussed the existence of periodic solutions in distinguishing between the cases in which the time $t$ occurs explicitly or not in the functions $X_i$ of the equation (*). Indeed, if $t$ is explicitly contained in the $X_i$, periodic solutions must have the same period as the $X_i$. Otherwise, periodic solutions may have any period,  and the issue is then to investigate the variation of this period in function of $\mu$. Poincaré's main statement in this respect is that if a periodic solution of period $T$ exists for $\mu=0$, then periodic solutions will still exist for small values of $\mu$, with a period close to $kT$ (with $k \in \mathbb{N}$).\cite[p.95]{Poincar1892} The issue of ``perturbation'' therefore already occurs in chapter III, i.e., before the introduction of the issues of approximations by periodic solutions. These general considerations are then applied to various particular ``applications,'' in which Poincaré discusses the existence of periodic solutions,  especially  in relation to the three-body-problem\cite[p.95-108]{Poincar1892}  and to the general problem of dynamic.\cite[p.109-152]{Poincar1892} The fourth chapter is devoted to approximations by periodic solutions. It opens with the ``equations of variations,''\cite[p.156-159, 162-264]{Poincar1892} which  aim at  introducing the ``characteristic exponents.''\cite[p.176]{Poincar1892}  The structure of the third chapter is then repeated, i.e., the distinction between the equations in which $t$ occurs explicitly or not, as well as the list of particular ``applications.''

 We have  already discussed  the connection between characteristic exponents and stability, and therefore between the roots of the equation in $S$ and the behavior of sets of trajectories. Yet, characteristic exponents play an even more crucial role in issues of perturbation in function of a parameter  $\mu$. Given a periodic solution for $\mu  = 0$, the multiplicity of the roots of its equation in $S$ plays not only a key role in the existence of a periodic solution for small values of $\mu >0$, but, as was proved by Poincaré, the characteristic  exponents  can actually  be developed in power series in $\sqrt{\mu}$.

\subsection{Multiplicity and perturbations} 

The  notion of multiple root that developed in the framework of Hermite's algebraic theory of forms supports an  analytic approach to the algebraic issues related to the secular equation. It indeed allows to analyze the variation of the number of distinct roots of an equation in function of the unknown $S$. Let us quote Darboux's 1874 memoir once again :

\begin{quote}
The number of positive squares in the form can only change if $S$ passes through a root of the equation [...], and in that case, the variation of the number of positive squares of the form cannot  be higher than the order of the multiplicity of the root under consideration. \cite{Darboux:1874} \footnote{Le nombre de carrés positifs de la forme ne peut changer que si $S$ passe par une racine de l'équation [...], et dans ce cas le nombre de carrés positifs de la forme ne peut varier d'une quantité supérieure à l'ordre de multiplicité de la racine considérée.}
\end{quote}

One finds some echoes of the above statement in Poincaré's investigation of the behavior of periodic trajectories in function of some  perturbations by a small parameter $\mu$ :

\begin{quote}
I must first observe that a periodic solution can disappear when $\mu$ passes from the value $-\epsilon$ to the value  $+\epsilon$ only if the equation has a multiple root for $\mu=0$ ; in other words, a periodic solution can only disappear after  mingling with another periodic solution [...] with which it will have exchanged its stability. Therefore, periodic solutions disappear by pairs similarly as the real roots of algebraic equations.\cite[p.83]{Poincar1892}\footnote{J'observe d'abord qu'une solution périodique ne peut disparaître quand $\mu$ passe de la valeur $-\epsilon$ à la valeur $+\epsilon$ que si pour $\mu=0$, l'équation admet une racine multiple ; en d'autres termes une solution périodique ne peut disparaître qu'après s'être confondue avec une autre solution périodique [...] avec laquelle elle aura échangé sa stabilité. Donc les solutions périodiques disparaissent par couples à la façon des racines réelles des équations algébriques.}  
\end{quote}

That the framework of Hermite's legacy was underlaying Poincaré's approach is made clear by the latter's use of the exact same notations  $\Delta_i$ as the former, or as his other followers such as Darboux.\cite[p.91-92]{Poincar1892} Moreover, Hermite's approach to Cauchy's reformulation of Lagrange's  procedures provides a structuration to the key section in which Poincaré introduces the equations of variations and the characteristic exponents.  Yet, in the \textit{Méthodes nouvelles}, the equation in $S$,  as well as its minors $\Delta_i$ may not only designate the determinant of a linear system with constant coefficient, but also  the  functional determinants extracted from the jacobian matrix associated to a linear system with periodic functions as coefficients. For instance, Poincaré used the implicit function theorem to prove that, if the time occurs explicitly in the equations,  then, given a periodic solution for $\mu=0$,  there is still a periodic solution  for small values of $\mu >0$ provided that the functional determinant corresponding to the equation in $S$ of the given periodic solution  does not vanish, i.e., if none of the characteristic exponent is zero.\cite[p.181]{Poincar1892}

 In transferring the notion of multiplicity from algebraic roots to periodic trajectories, it was  from  a preexisting algebraic mold that Poincaré was casting out the analysis of the perturbations of periodic solutions in function of the parameter  $\mu$. This algebraic cast is well illustrated by the key role played by the ``special discussions'' \cite[p.91, 159]{Poincar1892} that are devoted to the multiplicity of characteristic exponents, as well as by the way this issue is tackled, i.e., in discussing which of the $\Delta_i$ vanish simultaneously.\cite[p.91, 159, 173]{Poincar1892} More importantly, this algebraic cast plays a key role in most  statements relative to the existence of periodic solutions and to their behaviors after perturbations. Let us  exemplify this situation by quoting a few of these statements. For instance, Poincaré stated that,

\begin{quote}
 in the case when the differential equations do not include the time explicitly, if a periodic solution exists for $\mu=0$, one of the characteristic exponent of this solution has to be equal to zero ; moreover, if none of the other exponents is equal to zero, a periodic solutions will still exist for the small values of $\mu$. \cite[p.183]{Poincar1892} \footnote{ Ainsi, si les équations différentielles ne contiennent pas le temps explicitement, si elles admettent une solution périodique pour $\mu=0$, l'un des exposants caractéristiques de cette solution sera toujours nul ; si, de plus, aucun autre de ces exposants n'est nul, il y aura encore une solution périodique pour les petites valeurs de $\mu$.}
\end{quote}

In the case of the equations of dynamic, ``the characteristic exponents are two by two equals but of opposite signs.''\footnote{``Les exposants caractéristiques sont deux à deux égaux de signe contraire.''} In the case of the three-body problem,\footnote{In this case, two of the characteristic exponents always vanish because the original system is Hamiltonian. That the remaining two exponents do not vanish when $\mu >0$ implies that they can be expanded in convergent power series in $\sqrt{\mu}$.} ``the periodic solutions of the three-body problem have two, and only two, characteristic exponents.''\cite[p.218]{Poincar1892}\footnote{``Les solutions périodiques du problème des trois corps ont deux exposants caractéristiques nuls, mais elles n'en ont que deux.''}.

\subsection{The surface-of-section method}

The algebraic cast of Poincaré's strategy sheds light on several key innovations of the \textit{Méthodes nouvelles}, including the famed surface-of-section method, as we call it nowadays. This iterative method has often been celebrated as the first discrete recurrence to appear in dynamical systems (where time, no longer continuously varying, is symbolized by integers) and thereby as one of the origins of chaos theory. It is yet molded on the very same algebraic cast we have discussed in the previous section. 

It is well known that Poincaré forged the elements of a qualitative, geometric analysis making it possible, when differential equations are not solvable, to know the general look of the solutions and to state global results.\footnote{In modern parlance, this approach results in analyzing the 
phase portraits of the differential equations ; the phase space is the space of the bodies' position and momentum (velocity). It has thus $6n$ dimensions when $n$ is the number of bodies under consideration.} As a first step, he established a general classification of solutions in two dimensions in terms of singular points (centers, saddle points, nodes, and foci). His fundamental result was the following: among all the curves not ending in a singular point, some are periodic (they are limit cycles), and all the others wrap themselves asymptotically around limit cycles. Starting from behavior in the neighborhood of singular points, limit cycles and transverse arcs therefore provide a rather precise knowledge of trajectories. 

The surface-of-section method is based on  plane sections of a set of three-dimensional trajectories in the neighborhood of a periodical solution. A periodic trajectory is represented geometrically by a closed curve. One can thus consider a plane orthogonal to this curve. In this section plane, the periodic solution is represented by a fixed point $M$, while a non periodic trajectory intersects the plane in a point sequence $M_0$, $M_1$, $M_2$, .... Given a system of coordinates in the section plane, one can then define the transformation $T(z)$ that turns the coordinate of  $M_i$, i.e. a complex number $z$, into the one of  $M_{i+1}$, i.e. $T(z)$.

But linear approximations also play a crucial role in the  surface-of-section method, in a very similar way as the approach based on the ``equations of variations'' for approximations trajectories by periodic solutions. After proving that the transformation $T$ is holomorphic, Poincaré analyzed its development into power series, and eventually reduced this development to its linear term.  The sequence of points $M_i$ is then defined by the iterations $T^i$ of a linear operator. This approach leads to  a linear system of differential equations with constant coefficients. The ``equation in $S$'' of such a system, and its two  roots  $S_1=e^{\alpha_1 \omega}$ and $S_2=e^{\alpha_2 \omega}$,  provide the following analytic representation of $T$ :  
\[
x= A_1 e^{\alpha_1 t} \phi_1(t)+A_2 e^{\alpha_2 t} \phi_2(t) \ 
y= A_1 e^{\alpha_1 t} \psi_1(t)+A_2 e^{\alpha_2 t} \psi_2(t)
\]
(with $A_i$ constant ; $\phi$ and $\psi$ trigonometric sums).

Discussions on the algebraic nature of the roots of this equation followed, for the purpose of identifying different types of situations on the model of the criterion of stability for periodic solutions. Both discussions on stability were indeed molded on the same algebraic cast. Poincaré for instance stated that ``if the roots are real, positives and distinct, such that one is greater than 1 and the other lesser than 1,  then there exists two invariant curves in the section plane.''
%\footnote{Such a conclusion was extremely important for later research works. Hadamard (1901) extended Poincaré's  results to non-analytic cases. Lattès (1906) showed that the characteristic exponents are essential elements for the problem, and studied the   invariant curves  for more generic cases than those examined by Poincaré. Birkhoff (1920) gave a more precise definition of the surface of section and  of  the notion of invariant curves. See Roque (2007).} 
Such discussions on the algebraic nature (esp. the multiplicity) of the equation in $S$ underlain  the main results based on the surface-of-section method, especially the ones related to the stability of the solar system, which eventually lead to the introduction of homoclinic trajectories. 

In the plane of section, the issue of stability of flows of trajectories is related to the one of the existence of some invariant curves that would define some boundaries in which all the points $M_i$ would be trapped. It is well known that Poincaré tackled this issue in discussing the qualitative geometrical properties of curves in the section plane.\cite[p.199]{Poincar1886h} 
 Let us consider the case of a flow of asymptotic trajectories which slowly either approach or move away from an unstable given periodic solution, thereby generating families of curves which fill out surfaces and which asymptotically approach the curve representing the generating periodic solution. These surfaces correspond to curves in the transverse section for the investigation of which Poincaré developed his theory of invariant integrals.\footnote{Invariant integrals are differential forms whose integrals over suitable manifolds preserve their value when the manifolds are transported by the flow. This notion is introduced in the third volume of the \textit{Méthodes nouvelles} as the integrals of the equations of variations. It is discussed in respect with the multiplicity of the characteristic exponents.\cite[p.48]{Poincar1899b}.} Poincaré showed that  if the  corresponding curves meet in a closed curve, the flow remains confined in a certain region of space, which proves that the system is stable. 

 It  was precisely at this point that Poincaré committed his famous mistake in the memoir he addressed for the price of king Oscar II. Indeed, asymptotic trajectories do not correspond necessarily to closed curves in the section place. On the contrary,  neighborhoods of an unstable period oscillation can give rise to very complex trajectories.\footnote{Levi-Civita and Birkhoff showed later that such complexity also appears in the neighborhood of a stable periodic solutions, which forced a reassessment in the definition of stability. See \cite{Roque2011}.} After this error had been pointed out by Lars Edvard Phragmen, it was in discussing the multiplicity of characteristic exponents in connection with the convergence of the series expansions of the characteristic exponents in power of $\sqrt{\mu}$, that Poincaré eventually showed the existence of an infinity of doubly-periodic  (or homoclinic) solutions. Indeed, in the case of an autonomous Hamiltonian system, all the characteristic exponents are zero when $\mu=0$, and their series development in integer powers of $\sqrt{\mu}$ are divergent (i.e., these are asymptotic series in modern parlance),\cite[p.128]{Barrow-Greene1996} which implies the existence of trajectories with unpredictable long-time behavior. A doubly asymptotic trajectory can begin by being very close to the periodic solution when $t$ is large and negative; but then it moves away and deviates greatly from the periodic solution before getting close again to this solution when $t$ is large and positive. Moreover, the existence of a doubly asymptotic trajectory actually means that an infinite number of such trajectories exist.

\subsection{Open questions}

The various scales of collective dimensions we have discussed in the present paper have allowed us to analyze some individual innovations of Poincaré, such as his strategy of  small perturbations of periodic trajectories in regard with Lagrange's approach to small oscillations, or his transfer to differential equations of Hermite's developments  of what used to be a shared algebraic culture. Yet,  Poincaré's iterative processes  also point to some open questions in regard with the collective dimensions in which the latter's works took place. 

First, some iterative procedures had already been developed by  Lagrange in his works on the secular equation. These iterations aimed at  providing astronomers with an effective method for  integrating  linear differential systems, in contrast with the algebraic procedures based on the secular equation which required
 the resolution of an algebraic equation of degree $n$.  Lagrange first  investigated the issue of the stability of the solar system by decomposing a system of 5 planets into two sub-systems of 2 and 3 planets respectively, each associated to secular equations of degree $2$ and $3$ which he solved by radicals.  But he also developed another method for the needs of astronomers. In modern parlance, this method resorts to the iteration of a symmetric matrix for expressing its eigenvectors ; it is often designated nowadays as ``Le Verrier's method'' in linear algebra. Indeed, Lagrange's method has circulated in observatories and had especially been used in Le Verrier's investigations of the secular inequalities.

But the iterative processes of the surface-of-section method may also be compared to the ones of  the Newton method for approximating  algebraic roots by graphical iteration. Recall that this method starts with an approximation $a$ of a given root of the equation $P(x)=0$ under consideration. It consists in considering a small variation $a+\xi$ and in neglecting all non linear terms in the development of $P(a+\xi)$. One then gets an equation of the first degree in $\xi$, which resolution provides the value of $a+\xi$ with which  the procedure can be iterated.  

In the 1870s  the Newton method had actually been connected by  Edmond Laguerre to Hermite's approach to Sturm's theorem. We have seen that, already at the time of Sturm, the aim was to develop a general method for dealing with both algebraic and differential equations. Fourier, especially, had discussed the Fourier method in connection with Descartes' rule of signs, which played a key role in the statement of Sturm theorem. For this reason, several mathematicians, such as Cayley, rather designated the Newton method as the one of ``Newton-Fourier''. But Laguerre's aim was more specifically to generalize Hermite's approach to differential equations, i.e., a goal very similar to the one Poincaré would achieve in introducing the notion of multiple  periodic trajectories. Actually, this specific aspect of Laguerre's work was precisely the one Poincaré celebrated in his eulogy of the former in 1897. In the 1880s, prior to Poincaré's \textit{Méthodes nouvelles}, several works published in France took up with Laguerre's approach to differential equations. Most of them referred to Hermite's approach to Sturm theorem in connection with Descartes' rule of signs, the Newton method, and Fourier's upper bound for the roots of an algebraic equation. Yet, these works have not been analyzed up to now and the question of the collective framework in which they were developed remains open. Investigating this issue further would certainly shed light on some aspects of Poincaré's iterations processes, as well as on their early reception at the beginning of the 20th century by mathematicians such as Gabriel Koenigs, Hadamard, Samuel Lattès, Pierre Fatou, Gaston Julia,  Birkhoff, or Joseph Fels Ritt.

The open questions set above are all related to the fact that process of acculturations, that play a key role in the evolutions of any culture, cannot be separated from their social contexts. Interactional mathematical cultures such as the one of the secular equation raise issues related to the intertwining between the various infrastructures (networks of text and social spaces), and superstructures (institutions, journals, nations, etc.) in which any given individual's work take place.  As we shall see in the next section, processes of acculturation, i.e., the embedding of some external aspects into the internal coherence of a culture, always cause phenomena of chains reactions that  cause unexpected evolutions at each scale of a culture.\footnote{In this sense, processes of acculturations have been designated as ``total social phenomena'' in \cite{Bastide1956}.}

%Voir [Alexander 1994] pour une pré- histoire du sujet, de la méthode de Newton à ce dont il est question ici. Cf audin

%Raises the question of the connection between Poincaré's title and some other "Méthodes nouvelles" : "Les Méthodes nouvelles pour la résolution approchée des équation numériques" published in  1818 by Legendre.

% + anglais. Weyr ??? Préciser que ici = micro histoire

\section{Poincaré's specific algebraic practices as resulting from processes of acculturations to Jordan's algebraic culture}

We have seen that several statements of the \textit{Méthodes nouvelles} bear witness of the model role played by Hermite's specific approach to the secular equation for the strategy developed by Poincaré. Yet, this pespective is not sufficient for restoring the full individual specificity of Poincaré's own algebraic practices. As a matter of fact, a strong component of the algebraic mold from  which Poincaré casted out his new methods does not appear explicitly in the statement of theorems but much more implicitly in some procedures of proofs, such as the reductions of ``Tableaux'' to their canonical forms. As we shall see in this section, the main specificities of Poincaré's algebraic practices can be analyzed as resulting from a process of acculturation to Jordan's approach within an  cultural system dominated by Hermite's legacy.

 It is nevertheless not the place in the present paper to develop Jordan's approach in details.  We shall limit ourselves to identifying the specific cultural traits Poincaré pecked from Jordan's works, while ignoring its global coherence. For this reason, Hermite's legacy is more relevant than Jordan's for analyzing the strategy Poincaré developed in celestial mechanics. Actually,  the algebraic practices Poincaré had developed since the early 1880s in connection with Jordan's \textit{Traité des substitutions et des équations algébriques},\cite{Jordan1870} were greatly simplified, and even quite hidden,  in the \textit{Méthodes nouvelles}.

For instance, the introduction of the characteristic exponents was presented as an ``application of the theory of substitutions,''\cite[p.162]{Poincar1892} i.e., in a very vague allusion to Jordan's works with no further explanation.  In contrast, several of the great memoirs published in the 1880s in connection with Fuchsian functions opened with   some ``algebraic preliminaries'' devoted to detailed presentation of the  approach Poincaré had developed in  a crossbreeding of  Jordan's and Hermite's algebraic practices.

 Poincaré may  not have expected the readers of his treaties of celestial mechanics to be accustomed to substitutions group theory. In any case, this situation is a typical illustration of Roger Bastide's cut-off principle,\cite{Bastide1955} 
according to which individuals are  able to cut their own social space into several, coherent, sealed components in which they may act in very different ways.  As a matter of fact,  processes of acculturations do not create automatically some hybrid, or crossbred, individuals. Bastide's principle allows to analyze the discontinuities of social spaces and times that participate to the dynamic nature of algebraic cultures. As has already been observed in section 3, the fact that any individual is facing a plurality of contexts of cultural interactions results into the  inner heterogeneity of  any individual identity. Exactly as Sylvester or Cayley spared their continental publications from any symbolical algebraic issues, Poincaré's acculturation to  Jordan's group theory remained hidden in the context of celestial mechanics.

Let us now consider more closely Poincaré's  implicit use of Jordan's approach to substitutions. Recall that, in contrast with Lagrange, Hermite, and most authors concerned with  the secular equation, Poincaré did not deal with the  symmetric linear systems generated by the principles of mechanics but with  general linear systems generated by the equations of variations. In case of multiple roots, such systems cannot usually be turned into a diagonal form. Yet, the problem of multiple roots was completely solved in  the \textit{Méthodes nouvelles}  by  reducing  linear systems into a simplified Jordan canonical form.\cite[p.172-174]{Poincar1892} In doing so, Poincaré implicitly appealed to  a theorem Jordan had first stated in the case of linear groups in finite fields in 1868-1870. This theorem had been   generalized  in 1871 to linear systems of differential equations with constant coefficients on the field of complex numbers,\cite{Jordan1871} with a view on its application to the symmetric systems related to the secular equation. The canonical form theorem eventually allowed Jordan to prove, independently of Weierstrass, that the issue of the multiplicity of roots  is irrelevant for  Lagrange's criterion of stability.\cite{Jordan1872} 

\subsection{Jordan's algebraic culture}

The canonical form theorem was not an isolated result in Jordan's works. In contrast, this theorem played a key role in Jordan's own  reorganization of the cultural system in which  his early works took place from 1860 to 1868. This local mathematical culture was especially embodied by Jordan in connection with  the teaching (and textbooks) of Joseph Bertrand, Joseph-Alfred Serred, and Auguste Bravais, as well as by the study of papers of Cauchy, Victor Puiseux, Louis Poinsot, and Évariste Galois.\cite{Brechenmacher:2011} Even though Jordan himself  first designated this system  as ``the theory of order'' in his thesis of 1860, this designation did not point to what would be nowadays be considered as a ``theory'' but rather to an interactionist mathematical culture. 

 Jordan especially referred to Poinsot's characterization of the theory of order as maintaining a relation to algebra analogous to the relations between Gauss's higher arithmetic and usual arithmetic, or that between analysis situs and geometry.\cite{Jordan1860}  From 1808 to 1844, Poinsot had highlighted several times the transversal role played by the notion of ``order'' in the analogies encountered in various cyclic situations, such as the investigations of cyclotomic equations, congruences, symmetries, polyhedrons, and mechanical motions,\cite{Boucard:2011} to which Jordan added some analytic concerns for the groups of monodromy of  differential equations. 

Let us now characterize the global organization of this local  culture, as integrated by Jordan. As in the case of the secular equation, this  cultural system  is  characterized by the use of some specific  representations, ideals,  values, and forms of interactions between mathematical domains.

\subsubsection*{The analytic representation of substitutions}

Jordan's approach presents a particular declination of  the global uses of polynomial forms in the 19th century. Given a substitution $S$ operating on $p^n$ letters ($p$ prime), providing an  analytic representation to $S$ consists in indexing these letters  $a_k$ by a sequence of integers $k$ mod.$p$, for the purpose of finding a polynomial function  $f$ such that $S(a_k)=a_{f(k)}$.\footnote{In modern parlance, the function is defined in the finite field $F_{p^n}$.}

\subsubsection*{ Values : relations between general classes of objects}

While most of his contemporaries  focused on the particular objects associated with linear forms in one variable $(k \ ak+b)$, or linear fractional substitutions $(k \ \frac{ak+b}{ck+d})$,  Jordan aimed at dealing with  relations between general classes of objects.  Investigations of such relations  were  especially valued by Jordan. They were at the core of the latter's understanding of the ``theory of order,'' which he constrasted with classical concerns for quantities, magnitudes, or proportions.\cite{Jordan1881} This specificity of Jordan can be exemplified by contrasting the latter's approach to the analytic representation of substitutions with the one of Hermite.\footnote{These two approaches  were related to two very different readings of Galois's writings. See\cite{Brechenmacher:2011}.}  

On the one hand, Hermite provided in 1863 a complete characterization of the analytic representations of substitutions on  $p=3$, $p=5$ and $p=7$ letters. This issue was strongly connected to the particular groups of the modular equations of degree 3, 5, and 7 that Hermite, Kronecker, and Francesco Brioschi had been investigating in connection with Galois's works.\cite{Goldstein:2011} For instance, Hermite showed that any substitution on $5$ letters can be represented by combinations of the following polynomial forms :
\[
k \ ; k^2 \ ; k^3+ak
\]
Already in his thesis in 1860, Jordan, on the other hand, dealt with the problem of the analytic representation of substitutions in $n$ variables in introducing a chain of reductions  from the most general classes of groups to the most special ones (transitive groups, primitive groups, linear groups, symplectic groups, etc.).

\subsubsection*{  An ideal of generality : the successive reductions of the analytic representations of substitutions in $n$ variables}

The core of Jordan's approach  was the ``essential'' nature he attributed to a ``method of reduction'' for investigating the  relations between general classes of objects. In Jordan's first thesis, the main theorem introduces the general linear group in a finite field, $Gl_n(F_p)$, as generated from a two-step reduction of the  problem of finding the analytic representation of  general solvable groups. Given a set of indices $x=(x, x', x'',  ... , x^{(n)})$ ($x \in F_{p^n}$),  general linear groups are introduced as the ones in which substitutions take the following analytic form :
\[
\begin{vmatrix}
x & ax+bx'+cx''+...+d \\
x' & a'x+b'x'+c'x''+... +d'\\
x'' & a''x+b''x'+c''x''+...+d'' \\
.. & .....................
\end{vmatrix}
mod (p)
\]
 Later on, it was for  investigating further the reductions of general linear groups into chains of normal subgroups that Jordan stated the theorem on the invariance of the length and of the composition factors of the compositions series of a group, i.e., what is nowadays designated as the Jordan-Hölder theorem. It was also in this context that Jordan stated in 1868 the ``simplest'' form into which the analytical representation of a linear substitution can be reduced, whatever the multiplicity of its characteristic roots, i.e., what would later be designated as the Jordan canonical form theorem.\cite{Brechenmacher:2006a}

 In his investigations of general linear groups, Jordan aimed at reducing   any linear substitution on $F_{p^n}$ into an analytic ``form as simple as possible.''  In Galois's famous ``Mémoire sur les conditions de résolubilité  des équations par radicaux,''\cite{Galois1831} the main theorem was proven by appealing to the fact that in the case of one variable, the analytic form of the substitution $(k \ ak+b)$ can be easily decomposed into two cycles $(k \ gk)$ and $(k \ k+1)$. Yet, such a decomposition cannot be directly generalized to the case of  $n$ variables.\footnote{In modern parlance,  a matrix of  $n$ lines and $n$ columns can only be decomposed to a sequence of operations of the type  $(k \ gk)$ if this matrix can be diagonalized.} Let first consider the special case of linear  substitutions on $p^2$ letters (i.e. in 2 variables) that Jordan investigated in details in 1868 (thereby following Galois's second memoir\cite{Galois183?}). The determination of the simplest analytical forms was based on the polynomial decomposition of an equation of degree 2 (i.e., the characteristic equation of a matrix, in modern parlance). If this equation has two distinct roots, $\alpha$ and $\beta$, the letters can be reindexed in two blocks in such a way that the substitution is simply acting as a multiplication  on each block :\footnote{In modern parlance, one decomposes a vector space of dimension 2 into two subspaces each of dimension 1.}

\[
\begin{vmatrix}
z & \alpha z \\
u & \beta u 
\end{vmatrix}
\]
Yet, if the characteristic equation has a double root, the substitution cannot be reduced  to operations of multiplication as above, unless it is a trivial homothety. In the general case, the canonical form involves a combination of multiplications and additions :
\[
\begin{vmatrix}
z & \alpha z \\
u & \beta z+ \gamma u 
\end{vmatrix}
\]
The generalization of this canonical form  to $n$ variables would eventually lie beneath the global architecture of Jordan's 1870 \textit{Traité} :\cite{Brechenmacher:2011}

\begin{quote}
 This simple form
\[
\begin{vmatrix}
y_0, z_0, u_0, ..., y'_0, ... & K_0y_0, K_0(z_0+y_0), ... , K_0y'_0 \\
y_1, z_1, u_1, ..., y'_1, ... & K_1y_1, K_1(z_1+y_1), ... , K_1y'_1 \\
.... & ... \\
v_0, ... & K'_0v_0, ... \\
... & ... \\
\end{vmatrix}
\]
to which one can reduce the substitution à  $A$ by an adequate choice of indices,  will be designated as its canonical form.\cite[p.127]{Jordan1870}
\end{quote}
 From this point on, the reduction of the analytic representations of linear substitutions groups eventually replaced the notion of order in what Jordan considered as the ``very essence'' of his approach.

\subsubsection*{Specific interconnections between various branches of the mathematical sciences}

Analytic representations  also supported  analogies between the various issues Jordan tackled from 1860 to 1867: substitutions groups, algebraic equations, higher congruences, kinematics (motions of solid bodies), symmetries of polyhedrons, crystallography, the analysis situs of deformations of surfaces (including Riemann surfaces), and the groups of monodromy of linear differential equation (to which Jordan's second thesis was devoted). 

These interconnections were at first presented by Jordan as encompassed by the ``theory of order.'' Yet, this global organization changed between 1868 and 1870, as a consequence of  the specific approach Jordan had developed, and in close connection with some other contemporary works.\cite{Brechenmacher:2011b}  After 1868, it was actually Jordan's reduction of the analytic representation of linear substitutions to their  canonical form  that was supporting new interconnections. These were based on the transfer of the practices Jordan had  developed in the case of finite fields to situations involving the infinite field of complex numbers, such as  linear differential equations, and later  algebraic forms, a domain in which Jordan and Poincaré would eventually meet in the late 1870s.

\subsection{Poincaré's Tableaux and canonical forms}

Poincaré's acculturation to Jordan's  algebraic culture has  been  followed step by step in the context of the development of the theory of Fuchsian's functions from 1879 to 1884.\cite{Brechenmacher:2012a} This episode  illustrates once again the dynamic dimension of  the notion of culture : any individual appropriates his own culture progressively throughout his own life, without ever acquiring the totality of the  cultures of the various groups in which one belongs. Culture, therefore, is not a static heritage that would be transmitted as such from one generation to the other. It is a historical production built by the interactions between individuals and social groups.

Before analyzing further this situation, it is first customary to point out that Jordan's understanding of the ``theory of order'' was  for a long time almost completely disconnected from the local culture revolving around Hermite's legacy. To the point that, in 1873, Jordan's first intervention in the algebraic theory of forms caused  a very strong controversy with  Kronecker, whose approach was very close to the one of Hermite.\cite{Brechenmacher:2007a} When Jordan asked for Hermite's support, the latter eventually threatened to resign from the Academy if he was to be forced to read Jordan's works. Later on, from 1874 to 1880, Jordan struggled for getting acculturated to Hermite's approach, but he never completely adopted the latter's ideal of generality, which focused on the complete treatment of some particular cases,   values for effective computations, and specific ways to connect arithmetic, analysis, and algebra through specific objects such as modular equations.

Yet, Jordan's concerns for connecting his works on group theory to Hermite's algebraic theory of forms eventually established a point of contact with Poincaré's contemporary works, which in turn allowed the latter's acculturation to some traits of Jordan's approach.  

As has already been pointed out in section 3, processes of acculturation are key factors in the dynamic nature of any culture. These processes designate all the phenomena that result from a direct and continuous contact between groups of individuals of a different culture, and which especially cause some changes in the initial cultural models of each groups. On the one hand, these transformations usually result from the selection of some cultural elements from the alien culture. But on the other hand, the nature of this selection usually results from some deep tendency of the initial culture.\cite{Sapir1949}
The key role played by these processes call for the analysis of the ``reinterpretations''
by which ancient significations are attributed some new elements, or by which some new values change the cultural significations of some ancient forms.

This situation is illustrated by the way Poincaré eventually developed his own specific algebraic approach from the contact Jordan had created with Hermite's algebraic theory of forms. Some specific practices were appropriated, such as Jordan's canonical reduction, but with little concern for the global organization of Jordan's approach. Actually, these practices were embedded  within the global organization of Hermite's approach to algebraic forms.  But acculturation is nevertheless a total phenomenon that touches upon every level of a cultural system.  Even though Poincaré only pecked some specific traits of Jordan's algebraic practices, this  acculturation nevertheless resulted into a new organization of his initial culture. The notion of ``group'' was especially substituted to the one of ``form'' as the core element in Poincaré's approach. But  the internal logic of Hermite's theory of form remained predominant, with a permanence of its  representations, ideals,  values, and forms of interactions between various domains of  mathematics.

\subsubsection*{Forms of representations : Tableaux}

``Tableaux'' are a very visible consequence of the crossbreeding of Jordan's and Hermite's algebraic cultures. They indeed point to a form of representation  anchored in Hermite's legacy, and which Poincaré implemented with the procedures of reduction to canonical forms Jordan had initially developed with  his analytic representations. 

The notion of Tableau is close to, yet different from, what would be designated nowadays as a matrix.  This terminology  had  been used in France since the beginning of the 19th century for designating any ``form'' constituted by a complex of objects of the same type. It was still in this framework that Tableaux were used in the \textit{Méthodes nouvelles}, as for instance when Poincaré observed that a jacobian functional determinant ``can be considered as the tableau of the coefficients of a linear substitution.''\cite[p.175]{Poincar1892} In the framework of his works on the secular equation, Cauchy had already appealed to some operatory procedures on the Tableau formed by the determinant  $\Delta$  in order to extract some ``sub-Tableaux'' $\Delta_i$ from it, in connection with polynomial factorizations of the secular equation. The operatory character of this form of representation was later developed by Hermite in connection with the latter's investigations of various classes of equivalences of algebraic forms, whose coefficients were gathered into Tableaux that were manipulated by using  substitutions with either integer or real coefficients. In close connection with Hermite,  Sylvester made a specific use of this form of representation in his own approach to the secular equation, which eventually gave rise to the notion of ``minors'' extracted from a ``matrix'' as already alluded to before. In contrast with  contemporary mathematicians working in Hermite's legacy, such as Darboux, Jordan never made use of Tableaux before the late 1870s when the latter started publishing in the framework of Hermite's theory of forms. 

In the early 1880s, Jordan and Poincaré published  a series of memoirs on the algebraic theory of forms that were closely interconnected one with another. These memoirs not only document the crossbreeding of Tableaux with Jordan's canonical reduction  but also the persistent cultural differences between Jordan and Poincaré. While the former remained faithful to analytic representations of substitutions in addition to his new uses of  $n$ variables ``Tableaux,''  the latter explored in minute detail the various forms taken by Tableaux of a given small number of variables.

\subsubsection*{Ideals : generality through paradigmatic particular objects}

In coherence with Hermite's ideals for effective computations on specific objects, Poincaré rarely considered the general   canonical forms of  $n$ variables substitutions. In his works on Fuchsian (and hyperFuchsian) functions, Poincaré was dealing with particular analytic forms of real or complex substitutions, of 2 or 3 variables respectively, in the legacy of Hermite's works on modular equations :\footnote{In modern parlance, Poincaré's fuchsian (resp. hyperfuchsian) groups are discrete subgroups of $PSL_2(\mathbb{R})$ (resp.$PSL_3(\mathbb{R})$) while Poincaré''s kleinian groups are discrete subrgroups of  $PSL_2(\mathbb{C})$. See \cite{Gray:2000} } 

\[
(x, y \ ;  \ \frac{ax+by+c}{a''x+b''y+c''}, \frac{a'x+b'y+c'z}{a''x+b''y+c''}) 
 \]
 The classification of these substitution groups was based on their reduction to their  Jordan canonical forms, which Poincaré alternatively wrote analytically, \cite[p.349]{Poincar1884a}
\[
\begin{matrix}
(A) \ (x, y, z ; \ \alpha x, \beta y, \gamma z), \\
(B)  \ (x, y, z ; \ \alpha x, \beta y+z, \beta z), \\
(C) \ (x, y, z ;  \ \alpha x, \beta y, \beta z), \\
(D) \ (x, y, z, ; \  \alpha x+y, \alpha y+z, \alpha z), \\
(E) \ (x, y, z ;  \ \alpha x, \alpha y+z, \alpha z),
\end{matrix}
 \]
or by appealing to Tableaux :
\[
\begin{vmatrix}
\alpha & 0 & 0 \\
0 & \beta  & 0\\
0 & 0 & \gamma \\
 \end{vmatrix} 
\begin{vmatrix}
\beta & 0 & 0 \\
0 & \beta  & 0\\
0 & 1 & \alpha \\
 \end{vmatrix} 
\begin{vmatrix}
\beta & 0 & 0 \\
0 & \alpha  & 0\\
0 & 0 & \alpha \\
 \end{vmatrix} 
\begin{vmatrix}
\alpha & 0 & 0 \\
0 & \alpha  & 0\\
0 & 0 & \alpha \\
 \end{vmatrix} 
\begin{vmatrix}
\alpha & 0 & 0 \\
0 & \alpha  & 0\\
0 & 1 & \alpha \\
 \end{vmatrix} 
\]
The above Tableau notation was sometimes used by Poincaré for working with more general objects than linear fractional substitutions of 2 or 3 variables. Yet,  ``in order to avoid sacrificing clarity for the sake of generality,'' \cite[p. 28]{Poincar1881d} generality was usually expressed through the paradigmatic setting of the canonical forms of Tableaux with the smallest number of variables that allowed an exhaustive presentation of all possible cases.\footnote{On Poincaré's expression of generality through paradigm, see \cite{Robadey2004}, \cite{Robadey2006}.}

This framework is  exemplified by Poincaré's loose reference to the ``applications of the theory of linear substitutions'' in his \textit{Méthodes nouvelles}. Quite typically, Poincaré did not display the $n$ variable case but developed a paradigmatic example in the case of a linear system of 3 equations. As the latter concluded, ``We have supposed, for the sake of clarity, that we were dealing with a linear substitutions with only 3 variables ; but the same reasoning would apply whatever the number of variables.''\cite[p. 174]{Poincar1892}\footnote{Nous avons supposé, pour fixer les idées, que nous avions affaire à une substitution linéaire portant sur trois variables seulement ; mais le même raisonnement s'applique, quel que soit le nombre de variables.}

As said before, Poincaré's loose reference to the theory of substitutions mainly aimed at introducing a simplified presentation of the Jordan canonical form.  Poincaré first considered the following linear system,
\[
(1)
\begin{matrix}
\gamma_1=a_1\beta_1 +  a_2 \beta_2 + a_3 \beta_3 \\
\gamma_2=b_1\beta_1 +  b_2 \beta_2 + b_3 \beta_3 \\
\gamma_3=c_1\beta_1 +  c_2 \beta_2 + c_3 \beta_3 \\
 \end{matrix} 
\]
to which he associated a  linear substitution, ``linking the variables $\beta$ to the variables $\gamma$'',   of determinant,
\[
(2)
\begin{vmatrix}
a_1 & a_2 & a_3 \\
b_1 & b_2  & b_3\\
c_1 & c_2 & c_3 \\
 \end{vmatrix} 
\]
and of equation in $S$  : 
\[
\begin{vmatrix}
a_1-S & a_2 & a_3 \\
b_1 & b_2-S  & b_3\\
c_1 & c_2 & c_3-S \\
 \end{vmatrix} 
\]
According to the ``theory of linear substitutions,''  as Poincaré argued, the equation in $S$ and its minors are invariant for any linear substitutions applied simultaneously to the $\lambda$ and the $\beta$ of the system (1),\footnote{In modern parlance, these transformations $P$ of $Gl_n(\mathbb(C)$  define the  classes of matrices $P^{-1}AP$ similar to a given matrix $A$.} and for which the substitution (2) would turn into
\[
(3)
\begin{vmatrix}
a'_1 & a'_2 & a'_3 \\
b'_1 & b'_2  & b'_3\\
c'_1 & c'_2 & c'_3 \\
 \end{vmatrix} 
\]
Moreover, he wrote,
\begin{quote}
One can choose the $\lambda$ in such a way that the substitution is reduced to the simplest possible form, its \textit{canonical form}. This form consists in the following : 

If  all the roots of the equation in $S$ are simple, one can  nullify simultaneously $a'_2$, $a'_3$, $b'_1$, $b'_3$, $c'_1$, $c'_2$.
 
If  the equation in $S$ has a double root, one can  equal to zero  $a'_2$, $a'_3$, $b'_1$, $c'_1$ simultaneously, with  $b'_2$ = $c'_3$.

 If  the equation in $S$ has a triple root, one can  equal to zero $a'_2$, $a'_3$, $b'_3$  simultaneously, with  $a'_1$=$b'_2$ = $c'_3$.

In any case, one can alway suppose that the $\lambda$ have been chosen in such a way that
\[
a'_2=a'_3=b'_3=0
\]
\cite[p. 172]{Poincar1892}
\end{quote}

\subsubsection*{Interconnections  between domains through particular objects}

 Poincaré himself always distinguished carefully the notions of ``matrix'' and ``Tableau.''  In coherence  with Hermite's legacy,  matrices were conceived as the algebraic ``form'' underlaying the determinant from which one can extract a sequence of minors.\cite[p.90, 181, 187, 189]{Poincar1892}  ``Tableaux,'' on the other hand, were never defined precisely in a specific theoretical framework. The very function of Tableaux was actually to crossbreed various meanings in geometry, arithmetic, algebra, and analysis. Tableaux therefore supported  various analogies that participated of the links Poincaré established between distinct theories through particular objects, such as algebraic forms and fuchsian functions. This modality of interconnections shows the legacy of Hermite's ideals on the unity of mathematics.\cite[p.399]{Goldstein:2007} It is especially coherent with the way Hermite had constantly appealed to  key objects, such as forms and elliptic functions, to develop interactions between analysis, algebra, and arithmetic,\cite{GoldsteinSchappa:2007b} even though Poincaré's approach also included a strong geometric perspective. 

 As a consequence of Poincaré's acculturation to Jordan's algebraic practices, the notion of group gradually took over the ones of forms or functions as the key object for developing interconnections.\cite{Brechenmacher:2012a} Yet, in contrast with Jordan's approach, substitutions groups were always interconnected by Poincaré to functions, algebraic forms, and geometric objects. The reduction of a Tableau of 3 variables to its canonical form was typically understood  as a process of classification of the substitutions of a group (either finite, infinite, or even continuous), in regard with the analytic functions that such substitutions left invariants,  and, simultaneously, as a geometric process for finding the principal axes of a   surface, in regard with the  more arithmetical issue of  the identification of the classes of equivalence of an algebraic form. Algebraic forms especially still played a key role in  the \textit{Méthodes nouvelles}, especially in regard with the characterization of the multiplicity of characteristic exponents in regard with  particular mechanical situations.\cite[p.193]{Poincar1892}

\subsubsection*{Values : reduction and simplicity}

Because of their multivalent meanings, Tableaux could potentially be reduced to several kinds of ``canonical forms.'' These took the rather loose meaning of the ``simplest forms'' in regard with the problem under consideration. It was nevertheless precisely because of this loose meaning that the notion of ``canonical form,'' or ``reduced form,'' was   repetitively presented as an essential notion by Poincaré in the early 1880s. We have seen above that Jordan, also, had  attributed the ``essence'' of his approach to his ``method of reduction.'' Yet,  the terms  ``essence'' and ``reduction'' both took very different meanings in Jordan's and Poincaré's approaches. While the former appealed to an  abstract algebraic  approach to classes of groups, the latter followed the  loose signification of Hermite's  ``reduced form'' as a non-reified  norm of simplicity depending of the nature of the particular object under investigation :
 \begin{quote}
In order to represent each type, or  each sub-type [of classes of equivalence of cubic forms], we will choose, one of the forms of this type or sub-type that we shall designate as the canonical form $H$. The choice of the form $H$ is  nearly arbitrarily ; yet, in most cases, we shall be driven to prefer the simplest form of the type considered.\cite[p.203]{Poincar1881a}\footnote{On choisira dans chaque type ou dans chaque sous-type, pour le représenter, une des formes de ce type ou de ce sous-type que l'on appellera la forme canonique $H$. Le choix de la forme $H$ est peu près arbitraire ; toutefois on sera conduit, dans la plupart des cas, à choisir de préférence la forme la plus simple du type considéré.}
 \end{quote}

In his early works, Poincaré explicitly presented his approach as aiming at generalizing to cubic forms the ``very useful idea'' of ``reduced form'' developed by Hermite in his ``most elegant'' works on quadratic forms.\cite[p.28]{Poincar1881d} This loose notion of reduced form  is especially pervading the various canonical forms attributed to differential equations in the opening chapters of the \textit{Méthodes nouvelles}.\cite[p.9]{Poincar1892}

The acculturation of Jordan's approach within  Hermite's algebraic theory of forms eventually gave rise to some algebraic practices that were specific to Poincaré. These  were instrumental to the latter's capacity to intervene in a broad spectrum of mathematical issues in the 1880s.  The reduction of Tableaux to their Jordan canonical form especially provided Poincaré with an effective method for dealing with multiple roots in the equations in $S$ of linear systems.  The great memoirs of the time, especially the ones on fuchsian functions, usually open with some ``algebraic preliminaries'' devoted to the ``systems of definitions that will be useful hereafter,''\cite[p.203]{Poincar1881a} such as the analytic representation of linear substitutions, the notation of ``Tableaux,'' the equation in $S$,  canonical reductions of substitutions and forms, and their geometric and arithmetic interpretations.\footnote{See Poincaré's memoirs on algebraic forms,\cite[p.34]{Poincar1882a} fuchsian functions \cite[p.108]{Poincar1882e} continuous groups and partial differential equations,\cite{Poincar1883a} complex numbers,\cite{Poincar1884d} algebraic integration,\cite[p.300-313]{Poincar1884e},  integrations by series,\cite[p.316]{Poincar1886g} the arithmetic of fuchsian functions,\cite[p.463-505]{Poincar1887}, homologies in Analysis Situs,\cite[p.342-345]{Poincar1900} and continuous groups\cite[p.216-252]{Poincar1901}.}

 In contrast, the \textit{Méthodes nouvelles} open with preliminaries devoted to the existence of solutions of differential equations and to complex analysis. Yet, the reduction of Tableaux  is nevertheless implicitly underlaying the proofs of several crucial statements of this treaties.  Moreover, the transfer of the operatory procedures for manipulating Tableau  to issues of celestial mechanics  involved some innovations, such as their generalization to the infinite linear systems that Hill had previously considered.

\section*{Conclusions}
\addcontentsline{toc}{section}{Conclusions}

Looking upward at  Poincaré's approach to the three body problem provides a different picture than the retrospective celebrations of this approach as a starting point of chaos theory. Indeed, we have seen that the strategy developed by Poincaré in celestial mechanics can  be analyzed as molded on - or casted out - some specific algebraic practices for manipulating systems of linear equations. 

 The strategy of approximations by periodic trajectories, which is at the very core of the \textit{Méthodes nouvelles}, aims at introducing a very classical setting, i.e.,  linear systems with constant coefficients, in which Poincaré pushed a little toward a known difficulty, the one of multiple roots, for eventually finding something very new, thereby establishing new connections between algebraic forms, linear operations, analytical functions, probabilities, geometric and topologic interpretations, etc. This situation is actually quite typical of the way Poincaré was  approaching some new fields of researches in the 1880s.\footnote{This comment is related to the way  Catherine Goldstein has tackled the issue of what ``type  of great mathematician'' Poincaré was, with a case study of the latter's early works in number theory (centennial celebration of Poincaré held at IMPA, Rio de Janeiro, in November 2012). See also the way Poincaré created some new connections between matrices, lie algebras and associative algebras in 1884, as described in \cite{Brechenmacher:2012a}.} Algebra played a key role in the establishments of such links, but in a very different way than the abstract unifying algebraic structures modern mathematicians are used to appeal to.

As the structure of a cast-iron building  may be less noticeable than its creative façade, the linear algebraic cast of Poincaré's strategy was broken out of the mold in generating some new, non linear, methods of celestial mechanics.  But as the various components that are mixed in some casting process can still be detected in the resulting alloy, this algebraic cast points to some collective dimensions of Poincaré's methods, which sheds new light on  the novelty of Poincaré's \textit{Méthodes nouvelles}.

In the long run, and at a large scale, we have seen the key role played by a  shared algebraic culture based on the great mechanical treatises of the turn of the 19th century. References to the ``equation to the secular inequalities in Planetary theory'' were indeed used a the European level to identify some specific algebraic practices for manipulating linear systems. These practices were not limited to some polynomial procedures but went along with some collective representations,   ideals of generality, as well as with specific values  related to multiple roots. The specificity of the secular equation was instrumental to the circulation of these practices from one theory to another. In turn, new contexts made these practices evolve constantly, thereby showing the dynamic nature of the local, and individual, appropriations of a shared algebraic culture.

At a smaller scale, and during a more limited time-period, we have followed one of the lines of developments that gradually caused a decomposition of the shared culture of the secular equation. The legacy of Hermite's approach to Sturm's theorem  especially plays  a crucial role in Poincaré's approach to the variations and the perturbations of periodic solutions. But it was actually by mixing actively Hermite's theory of form with Jordan's approach to linear substitutions that Poincaré had developed some specific algebraic practices that were instrumental to his capacity to deal with various issues in a broad spectrum of the mathematical sciences, from arithmetic to celestial mechanics.

The way Poincaré introduced the equations of variations of  periodic trajectories to force linear systems with constant coefficients into his analysis of the three-body-problem echoes some results obtained by the so-called ``constructivist'' approach to the sociology and the history of sciences. In a word, those who take a ``constructivist position'' argue that scientists' decisions (such as for example to challenge a claim  or to open up a black boxed instrument) can be explained by reference to various active elements of their situation.
According to this approach, actors' choices are not only constrained by their aims,  and by the resources they select to advance them, but they are also guided by a complex of skills and technical competences that ``represent a  set of vested social interests \textit{within} the scientific community.''\cite[p.164]{Shapin1982}
They thus make the decisions they do because they seek to employ their particular specialist skills through developing new areas of work. In Andrew Pickering's classic investigations of high energy physics, this situation had been described as an attitude of ``opportunisms in context.''\cite{Pickering1984}
 In addition to such active components, the error in Poincaré's initial memoir also highlights the role played by some passive elements in constraining the production of scientific knowledge in ways that are beyond the control of those involved, a notion that was developed by Peter Galison in order to explain how experiments may have unexpected results.\cite[p.234-241]{Galison1987}
Moreover, Galison has emphasized the  importance of time in this process, in describing  constraints that are located at three different levels of temporal duration :  theoretical programs, such as the quest to unify all physical forces, are persisting as long-term constraints, whereas specific models of interpreting phenomena come and go more quickly, as well as the process of tinkering with instrumental set-ups. The situation we have seen in the present paper is very similar. We have especially contrasted the long-term dimensions of both the three-body problem and the practices attached to the secular equation, with the more local frameworks of Hermite's interpretation of Sturm's theorem or Jordan's analytic representations of substitutions. As with Galison's active, and time-embodied, interventions that shape phenomenal experience, we have described in this paper some elements that participated of Poincaré's own experience of mathematics and celestial mechanics.

Yet, in contrast with previous approaches, the present paper has analyzed some key aspects of the  individual specificity of Poincaré's algebraic practices as resulting from a casting process, understood as a  process of acculturation of interactional cultures within a spatialized, local culture. We have seen that the constructive and dynamic nature of such casting processes highlights  the  roles played by both individual's creativity and by some collective dimensions. The metaphor of the casting process aims at adapting to the history of sciences some key modern notions of social sciences, such as the ones of ``acculturation,''  ``syncretism,'' or ``crossbreeding.''  These notions   describe some original cultural configurations, which are not limited to some appropriations or assembly of heterogeneous elements through diffusions or circulations, but have to be regarded as genuine creations of new configurations.  Among the several types of such phenomena that have been described in social sciences, the notion of ``algebraic cast'' refers to the ``fusion'' model. 

On the hand, and from a retrospective point of view, the fusion model  contrasts with  the ``crossbreeding'' model in the sense that fusion makes it very difficult to distinguish  the initial elements involved, because these elements have transformed themselves into a new unified, and coherent system.\cite{Bastide1971} This approach allows to understand how Poincaré's \textit{Méthodes nouvelles} present simultaneously a  strong traditional dimension and a  complete novelty.  

On the other hand, and from a prospective point of view, the fusion model points to the emergence of a new sustainable cultural system. As a matter of fact, the algebraic practices that have been cast out  Jordan's approach in the melting pot of Hermite's algebraic culture did not only play a key role in most of Poincaré's works over the course of the latter's career. They were also quickly taken up by some other mathematicians and thereby gave rise to a new  local algebraic culture. It was, in a way, through the eyes of Poincaré that a new generation of mathematicians looked back at some more ancient works related to the secular equations, especially the ones of Jordan and Hermite. For mathematicians such as Léon Autonne, Hadamard, Edmond Maillet, Jean-Armand de Séguier, or Albert Châtelet,  the algebraic preliminaries of Poincaré's memoirs played a role similar to the one of a textbook,\footnote{Autonne especially followed Poincaré in opening systematically his own memoirs with algebraic preliminaries.\cite{Autonne:1885a} The great many treatises the former published at the turn of the 20th century actually institutionalized the uses of preliminaries  devoted to the Jordan canonical forms of ``Tableaux,'' with systematic interactions between interpretations in algebra, geometric, analysis, and arithmetic.
\cite[p.599-628]{Brechenmacher:2006a}  It was actually in this context that the Jordan  theorem was for the first time  stated with  a general matrix representation.\cite{Autonne:1905}

 Recall that  Poincaré had supervised Autonne's doctoral thesis on the algebraic integrals of linear differential equations at the École polytechnique from January to September 1881. This episode took place shortly after Poincaré got acculturated to Jordan's approach, and the former thus strongly advised his student to immerge himself into Jordan's 1870 \textit{Traité}. The epistolary communication between Poincaré and Autonne especially focuses on issues related to multiple roots in the equation in $S$,\cite{Autonne-Poincar} in relation with Jordan's works on linear differential equations.\cite[p.200]{Jordan1878}}
 through which they especially got acculturated to the algebraic practices of ``reductions'' of Tableaux to their multivalent ``canonical forms.'' These practices, in turn, circulated in a coherent network of texts from 1880 to 1914, which cultivated close contacts with  the works of the actors who were working either in Hermite's legacy, such as Hermann Minkowski, or the one of Jordan, such as Leonard Dickson.\cite{Brechenmacher2013b} These works played an important role in the institutionalization  of matrix theory at an international level in the 1920s.\cite{Brechenmacher:2010a}  From a prospective point of view, Poincaré's fusion of Hermite and Jordan's algebraic culture therefore gave rise to a sustainable algebraic culture, at least from 1880 to 1920. Moreover, Châtelet's 1950 critical edition of Poincaré's works in arithmetic and algebra exemplify that, from a retrospective point a view, it is very difficult to distinguish the various cultural elements involved in  the fusion : Poincaré's methods and results were indeed systematically translated by Châtelet in the new framework of matrix theory and linear algebra.

 In regard with  the local algebraic culture that developed from Poincaré's algebraic practices, the latter's works in celestial mechanics   appear as both isolated and not isolated. On the one hand, the algebraic cast of Poincaré's \textit{Méthodes nouvelles} was taken up by several mathematicians. On the other hand, and maybe as a consequence, the new methods Poincaré had cast out this algebraic mold in the specific framework of celestial mechanics were of little interest for most of the works of the network of the ``calcul des Tableaux.''

\bigskip 
Finally, the present paper has also shown that the relationships between celestial mechanics and the other branches of mathematical sciences in the 19th century was much more complex than a back-and-forth motion between application and abstraction. Not only did some specific procedures for dealing with linear systems emerge from some mechanical works. But the secular equation moreover generated a broadly shared algebraic culture in the 19th century by supporting the circulation of these procedures between various domains, thereby enriching them with new significations, and eventually returning to celestial mechanics with Poincaré's new methods. One may actually be tempted to describe the longue durée dimension of this situation in analogy with Poincaré's approach to the long run trajectories of celestial bodies : some trajectories may recede from the initial condition but nevertheless come back to their neighborhood over a long time.

\nocite{PoincarOeuvres}
\nocite{LaplaceOeuvres}
\nocite{LagrangeOeuvres}
\nocite{Franceschelli}
\nocite{Goldstein:1996}
\nocite{Gispert:1991}
\nocite{HarmanShapiro1992}

\addcontentsline{toc}{section}{Bibliography}
\bibliography{bibliographie}

\begin{thebibliography}{}
\expandafter\ifx\csname fonteauteurs\endcsname\relax
\def\fonteauteurs{\scshape}\fi
\expandafter\ifx\csname url\endcsname\relax
  \def\url#1{{\tt #1}}%
    \message{You should include the url package}\fi

\bibitem[Anderson, 1994]{Anderson1994}
\bgroup\fonteauteurs\bgroup Anderson\egroup\egroup{}, K. (1994).
\newblock Poincar{\'e}'s discovery of homoclinic points.
\newblock {\em Archive for History of Exact Sciences},
  48(2)\string:\penalty500\relax 133--147.

\bibitem[Archibald, 2011]{Archibald:2011}
\bgroup\fonteauteurs\bgroup Archibald\egroup\egroup{}, T. (2011).
\newblock Differential equations and algebraic transcendents: French efforts at
  the creation of a galois theory of differential equations (1880-1910).
\newblock {\em Revue d'histoire des math{\'e}matiques},
  17(2)\string:\penalty500\relax 371--399.

\bibitem[Aubin, 2009]{Aubin2008}
\bgroup\fonteauteurs\bgroup Aubin\egroup\egroup{}, D. (2009).
\newblock Observatory mathematics in the nineteenth century.
\newblock \emph{In} \bgroup\fonteauteurs\bgroup Robson\egroup\egroup{}, E. et
  \bgroup\fonteauteurs\bgroup Stedal\egroup\egroup{}, J., \'editeurs :  {\em
  Oxford Handbook for the History of Mathematics}, pages 273--298. Oxford
  University Press, Oxford.

\bibitem[Aubin \emph{et~al.}, 2010]{AubinBigg2010}
\bgroup\fonteauteurs\bgroup Aubin\egroup\egroup{}, D.,
  \bgroup\fonteauteurs\bgroup Bigg\egroup\egroup{}, C. et
  \bgroup\fonteauteurs\bgroup Sibum\egroup\egroup{}, O.~H., \'editeurs (2010).
\newblock {\em The Heavens on Earth: Observatories and Astronomy in
  Nineteenth-Century Science and Culture}.
\newblock Duke University Press, Durham.

\bibitem[Aubin et Dahan~Dalmedico, 2002]{AubinDahan}
\bgroup\fonteauteurs\bgroup Aubin\egroup\egroup{}, D. et
  \bgroup\fonteauteurs\bgroup Dahan~Dalmedico\egroup\egroup{}, A. (2002).
\newblock Writing the history of dynamical systems and chaos : Longue dur{\'e}e
  and revolution, disciplines and cultures.
\newblock {\em Historia Mathematica}, 29\string:\penalty500\relax 273--339.

\bibitem[Autonne, 1881]{Autonne-Poincar}
\bgroup\fonteauteurs\bgroup Autonne\egroup\egroup{}, L. (1881).
\newblock Correspondance d'{A}utonne (l{\'e}on) {\`a} {P}oincar{\'e} (henri).
\newblock Collection priv{\'e}e 75017. Reproduite {\`a} l'adresse
  http://www.univ-nancy2.fr/poincare/chp/.

\bibitem[Autonne, 1885]{Autonne:1885a}
\bgroup\fonteauteurs\bgroup Autonne\egroup\egroup{}, L. (1885).
\newblock Recherches sur les int{\'e}grales alg{\'e}briques des {\'e}quations
  diff{\'e}rentielles lin{\'e}aires {\`a} coefficients rationnels (second
  m{\'e}moire).
\newblock {\em Journal de l'{\'E}cole polytechnique},
  54\string:\penalty500\relax 1--30.

\bibitem[Autonne, 1905]{Autonne:1905}
\bgroup\fonteauteurs\bgroup Autonne\egroup\egroup{}, L. (1905).
\newblock {\em Sur les formes mixtes}.
\newblock A. Rey, and Gauthier-Villars, Lyon, and Paris.

\bibitem[Balandier, 1955]{Balandier1955}
\bgroup\fonteauteurs\bgroup Balandier\egroup\egroup{}, G. (1955).
\newblock La notion de "situation coloniale".
\newblock \emph{In} {\em La notion de "situation coloniale"}, pages 3--38. PUF,
  Paris.

\bibitem[Barnett, 1940]{Barnett1940}
\bgroup\fonteauteurs\bgroup Barnett\egroup\egroup{}, H. (1940).
\newblock Culture processes.
\newblock {\em American Anthropologist}, 42.

\bibitem[Barrow-Greene, 1994]{Barrow-Green1994}
\bgroup\fonteauteurs\bgroup Barrow-Greene\egroup\egroup{}, J. (1994).
\newblock Oscar {II}'s prize competition and the error in {P}oincar{\'e}'s
  memoir on the three body problem.
\newblock {\em Archive for History of Exact Sciences},
  48\string:\penalty500\relax 107--131.

\bibitem[Barrow-Greene, 1996]{Barrow-Greene1996}
\bgroup\fonteauteurs\bgroup Barrow-Greene\egroup\egroup{}, J. (1996).
\newblock {\em Poincare and the Three Body Problem}.
\newblock A.M.S. History of Mathematics.

\bibitem[Barth, 1969]{Barth1969}
\bgroup\fonteauteurs\bgroup Barth\egroup\egroup{}, F. (1969).
\newblock {\em Ethnic groups and boundaries. The social organization of culture
  difference}.
\newblock Universitetsforlaget, Oslo.

\bibitem[Bastide, 1954]{Bastide1955}
\bgroup\fonteauteurs\bgroup Bastide\egroup\egroup{}, R. (1954).
\newblock Le principe de coupure et le comportement afro-br{\'e}silien.
\newblock {\em Anais do XXXL Congresso Internacional de Americanistas, Sao
  Paulo}, 1\string:\penalty500\relax 493--503.

\bibitem[Bastide, 1956]{Bastide1956}
\bgroup\fonteauteurs\bgroup Bastide\egroup\egroup{}, R. (1956).
\newblock La causalit{\'e} externe et la causalit{\'e} interne dans
  l'explication sociologique.
\newblock {\em Cahiers internationaux de sociologie},
  21\string:\penalty500\relax 77--99.

\bibitem[Bastide, 1970b]{Bastide1970b}
\bgroup\fonteauteurs\bgroup Bastide\egroup\egroup{}, R. (1970b).
\newblock {\em Le {P}rochain et le {L}ointain}.
\newblock Cujas, Paris.

\bibitem[Bastide, 1971]{Bastide1971}
\bgroup\fonteauteurs\bgroup Bastide\egroup\egroup{}, R. (1971).
\newblock {\em Anthropologie appliqu{\'e}e}.
\newblock Payot, Paris.

\bibitem[Bastide, 1970a]{Bastide1970}
\bgroup\fonteauteurs\bgroup Bastide\egroup\egroup{}, R. (1996 (1970)a).
\newblock Continuit{\'e}s et discontinuit{\'e}s des soci{\'e}t{\'e}s et des
  cultures afro-am{\'e}ricaines.
\newblock {\em Bastidiana}, 13-14\string:\penalty500\relax 77--78.

\bibitem[Benedict, 1934]{Benedict1934}
\bgroup\fonteauteurs\bgroup Benedict\egroup\egroup{}, R. (1934).
\newblock {\em Patterns of Culture}.
\newblock Houghton Mifflin, New York.

\bibitem[B{\'e}n{\'e}ton, 1975]{Beneton1975}
\bgroup\fonteauteurs\bgroup B{\'e}n{\'e}ton\egroup\egroup{}, P. (1975).
\newblock {\em Histoire de mots : culture et civilisation}.
\newblock Presses de la FNSP, Paris.

\bibitem[Borchardt, 1846]{Borchardt1846}
\bgroup\fonteauteurs\bgroup Borchardt\egroup\egroup{}, K.~W. (1846).
\newblock Neue eigenschaft der gleichung, mit deren h{\"u}lfe man die
  saecularen storungen der planeten bestimmt.
\newblock {\em Journal f{\"u}r die reine und angewandte Mathematik},
  12\string:\penalty500\relax 38--45.

\bibitem[Boucard, 2011]{Boucard:2011}
\bgroup\fonteauteurs\bgroup Boucard\egroup\egroup{}, J. (2011).
\newblock Louis {P}oinsot et la th{\'e}orie de l'ordre : un cha{\^\i}non
  manquant entre {G}auss et {G}alois ?
\newblock {\em Revue d'histoire des math{\'e}matiques},
  17(1)\string:\penalty500\relax 41--138.

\bibitem[Bourdieu, 1980]{Bourdieu1980b}
\bgroup\fonteauteurs\bgroup Bourdieu\egroup\egroup{}, P. (1980).
\newblock {\em Le {S}ens pratique}.
\newblock Minuit.

\bibitem[Brechenmacher, 2006a]{Brechenmacher:2006a}
\bgroup\fonteauteurs\bgroup Brechenmacher\egroup\egroup{}, F. (2006a).
\newblock {\em Histoire du th{\'e}or{\`e}me de Jordan de la d{\'e}composition
  matricielle (1870-1930). Formes de repr{\'e}sentations et m{\'e}thodes de
  d{\'e}compositions.}
\newblock Th\`ese de doctorat, Ecole des hautes {\'e}tudes en sciences
  sociales.

\bibitem[Brechenmacher, 2006b]{Brechenmacher:2006d}
\bgroup\fonteauteurs\bgroup Brechenmacher\egroup\egroup{}, F. (2006b).
\newblock Les matrices : formes de repr{\'e}sentations et pratiques
  op{\'e}ratoires (1850-1930).
\newblock
  http://www.math.ens.fr/culturemath/histoire%20des%20maths/htm/Brechenmacher/matrices_index.htm.

\bibitem[Brechenmacher, 2007a]{Brechenmacher:2007a}
\bgroup\fonteauteurs\bgroup Brechenmacher\egroup\egroup{}, F. (2007a).
\newblock La controverse de 1874 entre {C}amille {J}ordan et {L}eopold
  {K}ronecker.
\newblock {\em Revue d'Histoire des Math{\'e}matiques},
  13\string:\penalty500\relax 187--257.

\bibitem[Brechenmacher, 2007b]{Brechenmacher:2007b}
\bgroup\fonteauteurs\bgroup Brechenmacher\egroup\egroup{}, F. (2007b).
\newblock L'identit{\'e} alg{\'e}brique d'une pratique port{\'e}e par la
  discussion sur l'{\'e}quation {\`a} l'aide de laquelle on d{\'e}termine les
  in{\'e}galit{\'e}s s{\'e}culaires des plan{\`e}tes (1766-1874).
\newblock {\em Sciences et techniques en perspective},
  1\string:\penalty500\relax 5--85.

\bibitem[Brechenmacher, 2010]{Brechenmacher:2010a}
\bgroup\fonteauteurs\bgroup Brechenmacher\egroup\egroup{}, F. (2010).
\newblock Une histoire de l'universalit{\'e} des matrices math{\'e}matiques.
\newblock {\em Revue de Synth{\`e}se}, 4\string:\penalty500\relax 569--603.

\bibitem[Brechenmacher, 2011]{Brechenmacher:2011}
\bgroup\fonteauteurs\bgroup Brechenmacher\egroup\egroup{}, F. (2011).
\newblock Self-portraits with {{\'E}}variste {G}alois (and the shadow of
  {C}amille {J}ordan).
\newblock {\em Revue d'histoire des math{\'e}matiques}, 17(fasc.
  2)\string:\penalty500\relax 271--369.

\bibitem[Brechenmacher, 2012a]{Brechenmacher:2012c}
\bgroup\fonteauteurs\bgroup Brechenmacher\egroup\egroup{}, F. (2012a).
\newblock Linear groups in galois fields. a case study of tacit circulation of
  explicit knowledge.
\newblock {\em Oberwolfach Reports}, 4-2012\string:\penalty500\relax 48--54.

\bibitem[Brechenmacher, 2012b]{Brechenmacher:2011b}
\bgroup\fonteauteurs\bgroup Brechenmacher\egroup\egroup{}, F. (2012b).
\newblock On {J}ordan's measurements.
\newblock {\em To appear},
  http://hal.archives-ouvertes.fr/aut/Frederic+Brechenmacher.

\bibitem[Brechenmacher, 2013a]{Brechenmacher:2012a}
\bgroup\fonteauteurs\bgroup Brechenmacher\egroup\egroup{}, F. (2013a).
\newblock Autour de pratiques alg{\'e}briques de {P}oincar{\'e} : h{\'e}ritages
  de la r{\'e}duction de {J}ordan.
\newblock {\em To appear},
  http://hal.archives-ouvertes.fr/aut/Frederic+Brechenmacher.

\bibitem[Brechenmacher, 2013b]{Brechenmacher2013b}
\bgroup\fonteauteurs\bgroup Brechenmacher\egroup\egroup{}, F. (2013b).
\newblock A history of galois fields.
\newblock {\em Chronos}, To appear.

\bibitem[Brechenmacher, 2008]{Brechenmacher:2008}
\bgroup\fonteauteurs\bgroup Brechenmacher\egroup\egroup{}, F. (To appear
  (2008)).
\newblock Algebraic generality vs arithmetic generality in the controversy
  between {C}. {J}ordan and {L}. {K}ronecker (1874).
\newblock \emph{In} \bgroup\fonteauteurs\bgroup Chemla\egroup\egroup{}, K.,
  \bgroup\fonteauteurs\bgroup Cambefort\egroup\egroup{}, Y.,
  \bgroup\fonteauteurs\bgroup Chorlay\egroup\egroup{}, R. et
  \bgroup\fonteauteurs\bgroup Rabouin\egroup\egroup{}, D., \'editeurs :  {\em
  Algebraic generality vs arithmetic generality in the controversy between {C}.
  {J}ordan and {L}. {K}ronecker (1874)}. Oxford University Press.

\bibitem[Cauchy, 1829]{Cauchy:1829}
\bgroup\fonteauteurs\bgroup Cauchy\egroup\egroup{}, A.-L. (1829).
\newblock Sur l'{\'e}quation {\`a} l'aide de laquelle on d{\'e}termine les
  in{\'e}galit{\'e}s s{\'e}culaires du mouvement des plan{\`e}tes.
\newblock {\em Exercices de math{\'e}matiques 4 in {\OE}uvres compl{\`e}tes
  d'Augustin Cauchy, Paris : Gauthier-Villars et fils, 1882-1974,},
  9(2)\string:\penalty500\relax 174)195.

\bibitem[Cauchy, 1831]{Cauchy1831}
\bgroup\fonteauteurs\bgroup Cauchy\egroup\egroup{}, A.-L. (1831).
\newblock M{\'e}moire sur les rapports qui existent entre le calcul des
  r{\'e}sidus et le calcul des limites et sur les avantages qu'offrent ces deux
  nouveaux calculs dans la r{\'e}solution des {\'e}quations alg{\'e}briques ou
  transcendantes.
\newblock {\em Lithography published par R. Taton in [Cauchy, \OE uvres, 2,
  XV]}, 182-361.

\bibitem[Cauchy, 1839a]{Cauchy1839b}
\bgroup\fonteauteurs\bgroup Cauchy\egroup\egroup{}, A.-L. (1839a).
\newblock M{\'e}moire sur l'int{\'e}gration des {\'e}quations lin{\'e}aires.
\newblock {\em Comptes rendus hebdomadaires des s{\'e}ances de l'Acad{\'e}mie
  des sciences}, 8\string:\penalty500\relax 369--426.

\bibitem[Cauchy, 1839b]{Cauchy1839a}
\bgroup\fonteauteurs\bgroup Cauchy\egroup\egroup{}, A.-L. (1839b).
\newblock M{\'e}thode g{\'e}n{\'e}rale propre {\`a} fournir les conditions
  relatives aux limites des corps dans les probl{\`e}mes de physique
  math{\'e}matique.
\newblock {\em Comptes rendus hebdomadaires des s{\'e}ances de l'Acad{\'e}mie
  des sciences}, 8\string:\penalty500\relax 374, 432, 459.

\bibitem[Cauchy, 1850]{Cauchy1850}
\bgroup\fonteauteurs\bgroup Cauchy\egroup\egroup{}, A.-L. (1850).
\newblock M{\'e}moire sur les perturbations produites dans les mouvements
  vibratoires d'un syst{\`e}me de mol{\'e}cules par l'influence d'un autre
  syst{\`e}me.
\newblock {\em Comptes rendus hebdomadaires des s{\'e}ances de l'Acad{\'e}mie
  des sciences}, 30\string:\penalty500\relax 17--24.

\bibitem[Cauchy, 1826]{Cauchy1826}
\bgroup\fonteauteurs\bgroup Cauchy\egroup\egroup{}, A.-L. (1882-1974 (1826)).
\newblock Application du calcul des r{\'e}sidus {\`a} l'int{\'e}gration des
  {\'e}quations diff{\'e}rentielles lin{\'e}aires {\`a} coefficients constants.
\newblock \emph{In} {\em {\OE}uvres compl{\`e}tes d'Augustin Cauchy}, volume
  (2)6, pages 252--255. Gauthier-Villars et fils, Paris.

\bibitem[Cayley, 1858]{Cayley:1858}
\bgroup\fonteauteurs\bgroup Cayley\egroup\egroup{}, A. (1858).
\newblock A memoir on the theory of matrices.
\newblock {\em Philosophical Transactions of the Royal Society of London},
  148\string:\penalty500\relax 17--37.

\bibitem[Chabert, 1992]{Chabert1992}
\bgroup\fonteauteurs\bgroup Chabert\egroup\egroup{}, J.-L. (1992).
\newblock Hadamard et les g{\'e}od{\'e}siques des surfaces {\`a} courbure
  n{\'e}gative.
\newblock \emph{In} \bgroup\fonteauteurs\bgroup Chabert\egroup\egroup{}, J.-L.,
  \bgroup\fonteauteurs\bgroup Chemla\egroup\egroup{}, K. et
  \bgroup\fonteauteurs\bgroup Dahan~Dalmedico\egroup\egroup{}, A., \'editeurs :
   {\em Chaos et d{\'e}terminisme,, collection Points Sciences, mai 1992, pp.
  306-330.} Seuil.

\bibitem[Chabert et Dahan~Dalmedico, 1992]{ChabertDahan1992}
\bgroup\fonteauteurs\bgroup Chabert\egroup\egroup{}, J.-L. et
  \bgroup\fonteauteurs\bgroup Dahan~Dalmedico\egroup\egroup{}, A. (1992).
\newblock Les id{\'e}es nouvelles de {P}oincar{\'e}.
\newblock \emph{In} {\em In Dahan-Dalmedico et al.}, pages 274--305.

\bibitem[Chenciner, 2007]{Chenciner2007}
\bgroup\fonteauteurs\bgroup Chenciner\egroup\egroup{}, A. (2007).
\newblock De la m{\'e}canique c{\'e}leste {\`a} la th{\'e}orie des syst{\`e}mes
  dynamiques, aller et retour: {P}oincar{\'e} et la g{\'e}om{\'e}trisation de
  l'espace des phases.
\newblock \emph{In} {\em Franceschelli et al.}, pages 13--36.

\bibitem[Christoffel, 1864]{Christoffel1864}
\bgroup\fonteauteurs\bgroup Christoffel\egroup\egroup{}, E.~B. (1864).
\newblock Ueber die kleinen schwingungen eines periodisch eingerichteten
  systems materieller punkte.
\newblock {\em Journal f{\"u}r die reine und angewandte Mathematik},
  63\string:\penalty500\relax 273--288.

\bibitem[Dahan~Dalmedico, 1984]{Dahan1980b}
\bgroup\fonteauteurs\bgroup Dahan~Dalmedico\egroup\egroup{}, A. (1984).
\newblock La math{\'e}matisation des th{\'e}ories de l'{\'e}lasticit{\'e} par
  {A}.-{L}. {C}auchy et les d{\'e}bats dans la physique math{\'e}matique
  fran{\c c}aise (1800-1840).
\newblock {\em Sciences et techniques en perspective},
  9\string:\penalty500\relax 1--100.

\bibitem[Dahan~Dalmedico, 1996]{Dahan1996}
\bgroup\fonteauteurs\bgroup Dahan~Dalmedico\egroup\egroup{}, A. (1996).
\newblock Le difficile h{\'e}ritage de henri poincar{\'e} en syst{\`e}mes
  dynamiques.
\newblock \emph{In} \bgroup\fonteauteurs\bgroup Greffe\egroup\egroup{}, J.-L.,
  \bgroup\fonteauteurs\bgroup Heinzmann\egroup\egroup{}, G. et
  \bgroup\fonteauteurs\bgroup Lorenz\egroup\egroup{}, K., \'editeurs :  {\em
  Dahan-Dalmedico, A. 1996. Le difficile h{\'e}ritage de Henri Poincar{\'e} en
  syst{\`e}mes dynamiques. In Greffe et al. (1996), 13--33.}, pages 13--33.
  Blanchard, Akademie Verlag, Paris.

\bibitem[Darboux, 1874]{Darboux:1874}
\bgroup\fonteauteurs\bgroup Darboux\egroup\egroup{}, G. (1874).
\newblock M{\'e}moire sur la th{\'e}orie alg{\'e}brique des formes
  quadratiques.
\newblock {\em Journal de math{\'e}matiques pures et appliqu{\'e}es},
  19(2)\string:\penalty500\relax 347--396.

\bibitem[Darigol, 2002]{Darigol2002}
\bgroup\fonteauteurs\bgroup Darigol\egroup\egroup{}, O. (2002).
\newblock Stability and instability in nineteenth century fluid mechanics.
\newblock {\em Revue d'histoire des math{\'e}matiques},
  8\string:\penalty500\relax 5--65.

\bibitem[de~la MSH, 1996]{Goldstein:1996}
de~la \bgroup\fonteauteurs\bgroup MSH\egroup\egroup{}, {\'E}., \'editeur
  (1996).
\newblock {\em L'Europe math{\'e}matique. Histoires, mythes, identit{\'e}s}.
\newblock GOLDSTEIN, Catherine and GRAY, Jeremy and RITTER, Jim, Paris.

\bibitem[Durand-Richard, 1996]{Durand-Richard:1996}
\bgroup\fonteauteurs\bgroup Durand-Richard\egroup\egroup{}, M.-J. (1996).
\newblock L'{\'e}cole alg{\'e}brique anglaise : les conditions conceptuelles et
  institutionnelles d'un calcul symbolique comme fondement de la connaissance.
\newblock \emph{In} {\em [Goldstein, Gray, Ritter 1996]}, pages 445--477.

\bibitem[Elias, 1939]{Elias1939}
\bgroup\fonteauteurs\bgroup Elias\egroup\egroup{}, N. (1939).
\newblock {\em {\"U}ber den Proze{\ss} der Zivilisation.}
\newblock Verlag Haus zum Falken, Basel.

\bibitem[Franceschelli \emph{et~al.}, 2007]{Franceschelli}
\bgroup\fonteauteurs\bgroup Franceschelli\egroup\egroup{}, S.,
  \bgroup\fonteauteurs\bgroup Paty\egroup\egroup{}, M. et
  \bgroup\fonteauteurs\bgroup Roque\egroup\egroup{}, T., \'editeurs (2007).
\newblock {\em Chaos et Syst{\`e}mes Dynamiques, {\'e}l{\'e}ments pour une
  {\'e}pist{\'e}mologie}.
\newblock Hermann, Paris.

\bibitem[Galison, 1987]{Galison1987}
\bgroup\fonteauteurs\bgroup Galison\egroup\egroup{}, P. (1987).
\newblock {\em How Experiments End}.
\newblock Chicago University Press, Chicago.

\bibitem[Galison, 2003]{Galison2003}
\bgroup\fonteauteurs\bgroup Galison\egroup\egroup{}, P. (2003).
\newblock {\em Einstein's clocks, {P}oincar{\'e}'s maps: empires of time}.
\newblock W. W. Norton, New York.

\bibitem[Galois, 1831a]{Galois183?}
\bgroup\fonteauteurs\bgroup Galois\egroup\egroup{}, {\'E}. (1831(?)a).
\newblock Fragment d'un second {M}{\'e}moire. {D}es {\'e}quations primitives
  qui sont solubles par radicaux.
\newblock \emph{In} {\em [Galois 1846]}, pages 434--444.

\bibitem[Galois, 1831b]{Galois1831}
\bgroup\fonteauteurs\bgroup Galois\egroup\egroup{}, {\'E}. (1831b).
\newblock M{\'e}moire sur les conditions de r{\'e}solubilit{\'e} des
  {\'e}quations par radicaux.
\newblock \emph{In} {\em [Galois 1846]}, pages 417--433.

\bibitem[Gilain, 1977]{Gilain1977}
\bgroup\fonteauteurs\bgroup Gilain\egroup\egroup{}, C. (1977).
\newblock {\em La th{\'e}orie g{\'e}om{\'e}trique des {\'e}quations
  diff{\'e}rentielles de {P}oincar{\'e} et l'histoire de l'analyse}.
\newblock Th\`ese de doctorat, Universit{\'e} de Paris7, Paris.

\bibitem[Gilain, 1991]{Gilain:1991}
\bgroup\fonteauteurs\bgroup Gilain\egroup\egroup{}, C. (1991).
\newblock La th{\'e}orie qualitative de poincar{\'e} et le probl{\`e}me de
  l'int{\'e}gration des {\'e}quations diff{\'e}rentielles.
\newblock \emph{In} {\em [Gispert 1991]}, pages 215--242.

\bibitem[Gispert, 1991]{Gispert:1991}
\bgroup\fonteauteurs\bgroup Gispert\egroup\egroup{}, H. (1991).
\newblock {\em La France math{\'e}matique. La Soci{\'e}t{\'e} Math{\'e}matique
  de France (1870-1914)}.
\newblock Cahiers d'histoire et de philosophie des sciences, Belin, Paris.

\bibitem[Goldstein, 1999]{Goldstein:1999}
\bgroup\fonteauteurs\bgroup Goldstein\egroup\egroup{}, C. (1999).
\newblock Sur la question des m{\'e}thodes quantitatives en histoire des
  math{\'e}matiques : le cas de la th{\'e}orie des nombres en france (1870-
  1914).
\newblock {\em Acta historiae rerum necnon technicarum},
  3\string:\penalty500\relax 187--214.

\bibitem[Goldstein, 2007]{Goldstein:2007}
\bgroup\fonteauteurs\bgroup Goldstein\egroup\egroup{}, C. (2007).
\newblock The {H}ermitian {F}orm of {R}eading the {D}isquisitiones.
\newblock \emph{In} {\em The {H}ermitian {F}orm of {R}eading the
  {D}isquisitiones}, pages 377--410.

\bibitem[Goldstein, 2011]{Goldstein:2011}
\bgroup\fonteauteurs\bgroup Goldstein\egroup\egroup{}, C. (2011).
\newblock Charles {H}ermite's {S}troll through the {G}alois fields.
\newblock {\em Revue d'histoire des math{\'e}matiques},
  17\string:\penalty500\relax 135--152.

\bibitem[Goldstein et Schappacher, 2007]{GoldsteinSchappa:2007b}
\bgroup\fonteauteurs\bgroup Goldstein\egroup\egroup{}, C. et
  \bgroup\fonteauteurs\bgroup Schappacher\egroup\egroup{}, N. (2007).
\newblock A book in search of a discipline (1801-1860).
\newblock \emph{In} {\em [Goldstein, Schappacher, Schwermer, 2007]}.

\bibitem[Gray, 1992]{Gray1992}
\bgroup\fonteauteurs\bgroup Gray\egroup\egroup{}, J. (1992).
\newblock Poincar{\'e}, topological dynamics, and the stability of the solar
  system.
\newblock \emph{In} {\em Harman and Shapiro}, pages 502--524.

\bibitem[Gray, 2000]{Gray:2000}
\bgroup\fonteauteurs\bgroup Gray\egroup\egroup{}, J. (2000).
\newblock {\em Linear differential equations and group theory from Riemann to
  Poincar{\'e}}, volume 2nd ed.
\newblock Birkh{\"a}user, Boston.

\bibitem[Hadamard, 1897]{Hadamard1897}
\bgroup\fonteauteurs\bgroup Hadamard\egroup\egroup{}, J. (1897).
\newblock Sur certaines propri{\'e}t{\'e}s des trajectoires en dynamique.
\newblock {\em Journal de math{\'e}matiques pures et appliqu{\'e}es}, (5)
  3\string:\penalty500\relax 331--387.

\bibitem[Hadamard, 1901]{Hadamard1901}
\bgroup\fonteauteurs\bgroup Hadamard\egroup\egroup{}, J. (1901).
\newblock Sur l'it{\'e}ration et les solutions asymptotiques des {\'e}quations
  diff{\'e}rentielles.
\newblock {\em Bulletin de la Soci{\'e}t{\'e} math{\'e}matique de France},
  29\string:\penalty500\relax 224--228.

\bibitem[Hadamard, 1913]{Hadamard1913}
\bgroup\fonteauteurs\bgroup Hadamard\egroup\egroup{}, J. (1913).
\newblock L'\oe uvre d'henri poincar{\'e}. le math{\'e}maticien.
\newblock {\em Revue de m{\'e}taphysique et de morale}, pages 617--658.

\bibitem[Harman et Shapiro, 1992]{HarmanShapiro1992}
\bgroup\fonteauteurs\bgroup Harman\egroup\egroup{}, P. et
  \bgroup\fonteauteurs\bgroup Shapiro\egroup\egroup{}, A., \'editeurs (1992).
\newblock {\em The investigation of difficult things: Essays on Newton and the
  history of exact sciences in honour of D. T. Whiteside}.
\newblock Cambridge University Press, Cambridge.

\bibitem[Hawkins, 1972]{Hawkins1972}
\bgroup\fonteauteurs\bgroup Hawkins\egroup\egroup{}, T. (1972).
\newblock Hypercomplex numbers, lie groups, and the creation of group
  representation theory.
\newblock {\em Archive for History of Exact Sciences,},
  8\string:\penalty500\relax 243--87.

\bibitem[Hawkins, 1975]{Hawkins1975}
\bgroup\fonteauteurs\bgroup Hawkins\egroup\egroup{}, T. (1975).
\newblock Cauchy and the spectral theory of matrices.
\newblock {\em Historia Mathematica}, 2\string:\penalty500\relax 1--20.

\bibitem[Hawkins, 1977]{Hawkins1977}
\bgroup\fonteauteurs\bgroup Hawkins\egroup\egroup{}, T. (1977).
\newblock Weierstrass and the theory of matrices.
\newblock {\em Archive for History of Exact Sciences,},
  17\string:\penalty500\relax 119--163.

\bibitem[Hermite, 1853]{Hermite1853}
\bgroup\fonteauteurs\bgroup Hermite\egroup\egroup{}, C. (1853).
\newblock Sur la d{\'e}composition d'un nombre en quatre carr{\'e}s.
\newblock {\em Comptes rendus hebdomadaires des s{\'e}ances de l'Acad{\'e}mie
  des sciences}, 37\string:\penalty500\relax 133--134.

\bibitem[Hermite, 1854]{Hermite1854}
\bgroup\fonteauteurs\bgroup Hermite\egroup\egroup{}, C. (1854).
\newblock Sur la th{\'e}orie des formes quadratiques.
\newblock {\em Journal f{\"u}r die reine und angewandte Mathematik}, 47.

\bibitem[Hermite, 1855]{Hermite1855}
\bgroup\fonteauteurs\bgroup Hermite\egroup\egroup{}, C. (1855).
\newblock Remarque sur un th{\'e}or{\`e}me de m. cauchy.
\newblock {\em Comptes rendus hebdomadaires des s{\'e}ances de l'Acad{\'e}mie
  des sciences}, 41\string:\penalty500\relax 181--183.

\bibitem[Hermite, 1857]{Hermite1857}
\bgroup\fonteauteurs\bgroup Hermite\egroup\egroup{}, C. (1857).
\newblock Sur l'invariabilit{\'e} du nombre des carr{\'e}s positifs et des
  carr{\'e}s n{\'e}gatifs dans la transformation des polyn{\^o}mes
  homog{\`e}nes du second degr{\'e}.
\newblock {\em Journal f{\"u}r die reine und angewandte Mathematik},
  53\string:\penalty500\relax 271--274.

\bibitem[Hill, 1877]{Hill1877}
\bgroup\fonteauteurs\bgroup Hill\egroup\egroup{}, G.~W. (1877).
\newblock {\em On the part of the motion of the lunar perigee which is a
  function of the mean motions of the Sun and Moon}.
\newblock Wilson, Cambridge.

\bibitem[Hill, 1878]{Hill1878}
\bgroup\fonteauteurs\bgroup Hill\egroup\egroup{}, G.~W. (1878).
\newblock Researches in the lunar theory.
\newblock {\em American Journal of Mathematics}, 1(5)\string:\penalty500\relax
  129--145.

\bibitem[Hill, 1886]{Hill1886}
\bgroup\fonteauteurs\bgroup Hill\egroup\egroup{}, G.~W. (1886).
\newblock On the part of the motion of the lunar perigee which is a function of
  the mean motions of the sun and moon.
\newblock {\em Acta Mathematica}, 8\string:\penalty500\relax 1--36.

\bibitem[Jacobi, 1834]{Jacobi1834}
\bgroup\fonteauteurs\bgroup Jacobi\egroup\egroup{}, C. G.~J. (1834).
\newblock De binis quibuslibet functionibus homogeneis secundi ordinis per
  substitutiones lineares.
\newblock {\em Journal f{\"u}r die reine und angewandte Mathematik},
  12\string:\penalty500\relax 191--268.

\bibitem[Jacobi, 1857]{Jacobi1857}
\bgroup\fonteauteurs\bgroup Jacobi\egroup\egroup{}, C. G.~J. (1857).
\newblock Ueber eine elementare transformation eines in bezug auf jedes von
  zwei variabeln-systemen linearen und homogenen ausdr{\"u}cks.
\newblock {\em Journal f{\"u}r die reine und angewandte Mathematik},
  53\string:\penalty500\relax 583--590.

\bibitem[Jordan, 1860]{Jordan1860}
\bgroup\fonteauteurs\bgroup Jordan\egroup\egroup{}, C. (1860).
\newblock {\em Sur le nombre des valeurs des fonctions. Th{\`e}ses
  pr{\'e}sent{\'e}es {\`a} la Facult{\'e} des sciences de {P}aris par {C}amille
  {J}ordan, 1re th{\`e}se}.
\newblock Mallet-Bachelier, Paris.

\bibitem[Jordan, 1870]{Jordan1870}
\bgroup\fonteauteurs\bgroup Jordan\egroup\egroup{}, C. (1870).
\newblock {\em Trait{\'e} des substitutions et des {\'e}quations
  alg{\'e}briques}.
\newblock Gauthier-Villars, Paris.

\bibitem[Jordan, 1871]{Jordan1871}
\bgroup\fonteauteurs\bgroup Jordan\egroup\egroup{}, C. (1871).
\newblock Sur la r{\'e}solution des {\'e}quations diff{\'e}rentielles
  lin{\'e}aires.
\newblock {\em Comptes rendus hebdomadaires des s{\'e}ances de l'Acad{\'e}mie
  des sciences}, 73\string:\penalty500\relax 787--791.

\bibitem[Jordan, 1872]{Jordan1872}
\bgroup\fonteauteurs\bgroup Jordan\egroup\egroup{}, C. (1872).
\newblock Sur les oscillations infiniment petites des syst{\`e}mes
  mat{\'e}riels.
\newblock {\em Comptes rendus hebdomadaires des s{\'e}ances de l'Acad{\'e}mie
  des sciences}, 74\string:\penalty500\relax 1395--1399.

\bibitem[Jordan, 1878]{Jordan1878}
\bgroup\fonteauteurs\bgroup Jordan\egroup\egroup{}, C. (1878).
\newblock M{\'e}moire sur les {\'e}quations diff{\'e}rentielles lin{\'e}aires
  {\`a} int{\'e}grale alg{\'e}brique.
\newblock {\em Journal f{\"u}r die reine und angewandte Mathematik},
  84\string:\penalty500\relax 89--215.

\bibitem[Jordan, 1881]{Jordan1881}
\bgroup\fonteauteurs\bgroup Jordan\egroup\egroup{}, C. (1881).
\newblock {\em Notice sur les travaux de {M}. {C}amille {J}ordan {\`a} l'appui
  de sa candidature {\`a} l'{A}cad{\'e}mie des sciences}.
\newblock Gauthier-Villars, Paris.

\bibitem[Kellert, 1993]{Kellert1993}
\bgroup\fonteauteurs\bgroup Kellert\egroup\egroup{}, S.~H. (1993).
\newblock {\em In the wake of chaos: unpredictable order in dynamical systems}.
\newblock Chicago University Press, Chicago.

\bibitem[Kronecker, 1874]{Kronecker1874}
\bgroup\fonteauteurs\bgroup Kronecker\egroup\egroup{}, L. (1874).
\newblock Ueber schaaren von quadratischen und bilinearen formen.
\newblock {\em Monatsberichte der K\"oniglich Preussischen Akademie der
  Wissenschaften zu Berlin}, pages 59--76, 149--156, 206--232.

\bibitem[Lagrange, 1766]{Lagrange1766}
\bgroup\fonteauteurs\bgroup Lagrange\egroup\egroup{}, J.-L. (1766).
\newblock Solution de diff{\'e}rents probl{\`e}mes de calcul int{\'e}gral.
\newblock \emph{In} {\em \OE uvres}, volume~1, page 471.

\bibitem[Lagrange, 1775]{Lagrange1775}
\bgroup\fonteauteurs\bgroup Lagrange\egroup\egroup{}, J.-L. (1775).
\newblock Nouvelle solution du probl{\`e}me du mouvement de rotation d'un corps
  de figure quelconque qui n'est anim{\'e} par aucune force
  acc{\'e}l{\'e}ratrice.
\newblock \emph{In} {\em in [Lagrange, \OE uvres, 3, p. 577-616].}

\bibitem[Lagrange, 1778]{Lagrange1778}
\bgroup\fonteauteurs\bgroup Lagrange\egroup\egroup{}, J.-L. (1778).
\newblock Recherches sur les {\'e}quations s{\'e}culaires des mouvements des
  n{\oe}uds, et des inclinaisons des orbites des plan{\`e}tes.
\newblock \emph{In} {\em [Lagrange, \OE uvres, 6, p. 635-709].]}, page 177.

\bibitem[Lagrange, 1774]{Lagrange1774}
\bgroup\fonteauteurs\bgroup Lagrange\egroup\egroup{}, J.-L. (1778 (1774)).
\newblock {\em Recherches sur les {\'e}quations s{\'e}culaires des mouvements
  des noeuds et des inclinaisons des plan{\`e}tes}.
\newblock M{\'e}moires de l'Acad{\'e}mie des Sciences de Paris.

\bibitem[Lagrange, 1783]{Lagrange1783}
\bgroup\fonteauteurs\bgroup Lagrange\egroup\egroup{}, J.-L. (1783).
\newblock Th{\'e}orie des variations s{\'e}culaires des {\'e}l{\'e}ments des
  plan{\`e}tes; {P}remi{\`e}re partie.
\newblock \emph{In} {\em Th{\'e}orie des variations s{\'e}culaires des
  {\'e}l{\'e}ments des plan{\`e}tes; {P}remi{\`e}re partie}.

\bibitem[Lagrange, 1784]{Lagrange1784}
\bgroup\fonteauteurs\bgroup Lagrange\egroup\egroup{}, J.-L. (1784).
\newblock Th{\'e}orie des variations s{\'e}culaires des {\'e}l{\'e}ments des
  plan{\`e}tes ; {S}econde partie.
\newblock \emph{In} {\em Th{\'e}orie des variations s{\'e}culaires des
  {\'e}l{\'e}ments des plan{\`e}tes ; {S}econde Partie}.

\bibitem[Lagrange, 1788]{Lagrange1788}
\bgroup\fonteauteurs\bgroup Lagrange\egroup\egroup{}, J.-L. (1788).
\newblock {\em M{\'e}canique analytique}.
\newblock Veuve Desaint, Paris.

\bibitem[Lagrange, 1892]{LagrangeOeuvres}
\bgroup\fonteauteurs\bgroup Lagrange\egroup\egroup{}, J.-L. (1867-1892).
\newblock {\em \OE uvres de Lagrange}.
\newblock Gauthier-Villars, Paris.

\bibitem[Laplace, 1775]{Laplace1775}
\bgroup\fonteauteurs\bgroup Laplace\egroup\egroup{}, P.-S. (1775).
\newblock M{\'e}moire sur les solutions particuli{\`e}res des {\'e}quations
  diff{\'e}rentielles et sur les in{\'e}galit{\'e}s s{\'e}culaires des
  plan{\`e}tes.
\newblock \emph{In} {\em [Laplace, \OE uvres, 8, p. 325-366]}. M{\'e}moires de
  l'Acad{\'e}mie des Sciences de Paris.

\bibitem[Laplace, 1776]{Laplace1776}
\bgroup\fonteauteurs\bgroup Laplace\egroup\egroup{}, P.-S. (1776).
\newblock Recherches sur le calcul int{\'e}gral et sur le syst{\`e}me du monde.
\newblock \emph{In} {\em [\OE uvres, 8, p.369-501]}.

\bibitem[Laplace, 1787]{Laplace1787}
\bgroup\fonteauteurs\bgroup Laplace\egroup\egroup{}, P.-S. (1787).
\newblock M{\'e}moire sur les in{\'e}galit{\'e}s s{\'e}culaires des
  plan{\`e}tes et des satellites.
\newblock \emph{In} {\em [Laplace, \OE uvres, 8, p. 49-92]}.

\bibitem[Laplace, 1789]{Laplace1789}
\bgroup\fonteauteurs\bgroup Laplace\egroup\egroup{}, P.-S. (1789).
\newblock M{\'e}moire sur les variations s{\'e}culaires des orbites des
  plan{\`e}tes.
\newblock \emph{In} {\em [Laplace, \OE uvres, 11, p. 295-306]}.

\bibitem[Laplace, 1799]{Laplace1799}
\bgroup\fonteauteurs\bgroup Laplace\egroup\egroup{}, P.-S. (1799).
\newblock {\em Trait{\'e} de m{\'e}canique c{\'e}leste}, volume~1.
\newblock Paris.

\bibitem[Laplace, 1912]{LaplaceOeuvres}
\bgroup\fonteauteurs\bgroup Laplace\egroup\egroup{}, P.-S. (1878-1912).
\newblock {\em \OE uvres de compl{\`e}tes de Laplace}.
\newblock Gauthier-Villars, Paris.

\bibitem[Laskar, 1992]{Laskar1992}
\bgroup\fonteauteurs\bgroup Laskar\egroup\egroup{}, J. (1992).
\newblock La stabilit{\'e} du syst{\`e}me solaire.
\newblock \emph{In} {\em Chaos et d{\'e}terminisme}, pages 170--212. A. Dahan
  Dalmedico, J-L. Chabert, K. Chemla (dir.).

\bibitem[Le~{V}errier, 1856]{LeVerrier1856}
\bgroup\fonteauteurs\bgroup Le~{V}errier\egroup\egroup{}, U. J.~J. (1856).
\newblock {\em Annales de l'{O}bservatoire de {P}aris}, volume~2.
\newblock Mallet {B}achelet, Paris.

\bibitem[Lejeune-Dirichlet, 1842]{Dirichlet1842}
\bgroup\fonteauteurs\bgroup Lejeune-Dirichlet\egroup\egroup{}, J. P.~G. (1842).
\newblock Recherches sur les formes quadratiques {\`a} coefficients et {\`a}
  ind{\'e}termin{\'e}es complexes.
\newblock {\em Journal f{\"u}r die reine und angewandte Mathematik},
  24\string:\penalty500\relax 533--618.

\bibitem[Lejeune-Dirichlet, 1846]{Dirichlet1846}
\bgroup\fonteauteurs\bgroup Lejeune-Dirichlet\egroup\egroup{}, J. P.~G. (1846).
\newblock Ueber die stabilit{\"a}t des gleichgewichts.
\newblock {\em Journal f{\"u}r die reine und angewandte Mathematik},
  32\string:\penalty500\relax 3--8.

\bibitem[L{\'e}vy-Strauss, 1950]{Strauss1950}
\bgroup\fonteauteurs\bgroup L{\'e}vy-Strauss\egroup\egroup{}, C. (1950).
\newblock Introduction {\`a} l'\oe uvre de {M}arcel {M}auss.
\newblock \emph{In} {\em Introduction {\`a} l'\oe uvre de {M}arcel {M}auss}.
  PUF, Paris.

\bibitem[L{\'e}vy-Strauss, 1958]{Strauss1958}
\bgroup\fonteauteurs\bgroup L{\'e}vy-Strauss\egroup\egroup{}, C. (1958).
\newblock {\em Anthropologie structurale}.
\newblock Plon, Paris.

\bibitem[L{\"u}tzen, 1984]{Lutzen1984}
\bgroup\fonteauteurs\bgroup L{\"u}tzen\egroup\egroup{}, J. (1984).
\newblock Joseph {L}iouville's work on the figures of equilibrium of a rotating
  mass of fluid.
\newblock {\em Archive for History of Exact Sciences},
  30(2)\string:\penalty500\relax 113--166.

\bibitem[Malinowski, 1944]{Malinowski1944}
\bgroup\fonteauteurs\bgroup Malinowski\egroup\egroup{}, B.~K. (1944).
\newblock {\em A Scientific Theory of Culture and Others Essays}.
\newblock The University of North Carolina Press, Chapel Hill.

\bibitem[Mawhin, 1996]{Mawhin1996}
\bgroup\fonteauteurs\bgroup Mawhin\egroup\egroup{}, J. (1996).
\newblock The early reception in france of the work of {P}oincar{\'e} and
  {L}yapunov in the qualitative theory of differential equations.
\newblock {\em Philosophia Scienti\ae}, 1(4)\string:\penalty500\relax 119--133.

\bibitem[Nabonnand, 2000]{Nabonnand2000}
\bgroup\fonteauteurs\bgroup Nabonnand\egroup\egroup{}, P. (2000).
\newblock Les recherches sur l'oeuvre de poincar{\'e}.
\newblock {\em La gazette des math{\'e}maticiens}, 85\string:\penalty500\relax
  34--54.

\bibitem[Nabonnand, 2005]{Nabonnand2005}
\bgroup\fonteauteurs\bgroup Nabonnand\egroup\egroup{}, P. (2005).
\newblock Bibliographie des travaux sur l'oeuvre de {P}oincar{\'e} (2001-2005).
\newblock {\em Philosophia Scienti\ae}, 9(1),\string:\penalty500\relax
  195--206.

\bibitem[Pickering, 1984]{Pickering1984}
\bgroup\fonteauteurs\bgroup Pickering\egroup\egroup{}, A. (1984).
\newblock {\em Constructing Quarks : A Sociological History of Particle
  Physics}.
\newblock Edinburg Univesity Press, Edinburg.

\bibitem[Poincar{\'e}, 1881a]{Poincar1881f}
\bgroup\fonteauteurs\bgroup Poincar{\'e}\egroup\egroup{}, H. (1881a).
\newblock M{\'e}moire sur les courbes d{\'e}finies par les {\'e}quations
  diff{\'e}rentielles.
\newblock {\em Journal de math{\'e}matiques pures et appliqu{\'e}es}, III,
  7\string:\penalty500\relax 375--422.

\bibitem[Poincar{\'e}, 1881b]{Poincar1881d}
\bgroup\fonteauteurs\bgroup Poincar{\'e}\egroup\egroup{}, H. (1881b).
\newblock Sur la repr{\'e}sentation des nombres par les formes.
\newblock {\em Comptes rendus hebdomadaires des s{\'e}ances de l'Acad{\'e}mie
  des sciences}, 92\string:\penalty500\relax 333--335.

\bibitem[Poincar{\'e}, 1881c]{Poincar1881a}
\bgroup\fonteauteurs\bgroup Poincar{\'e}\egroup\egroup{}, H. (1881c).
\newblock Sur les formes cubiques ternaires et quaternaires.
\newblock {\em Journal de l'{\'E}cole polytechnique},
  50\string:\penalty500\relax 199--253.

\bibitem[Poincar{\'e}, 1882a]{Poincar1882f}
\bgroup\fonteauteurs\bgroup Poincar{\'e}\egroup\egroup{}, H. (1882a).
\newblock M{\'e}moire sur les courbes d{\'e}finies par une {\'e}quation
  diff{\'e}rentielle (2e partie).
\newblock {\em Journal de Math{\'e}matiques}, 8\string:\penalty500\relax
  251--296.

\bibitem[Poincar{\'e}, 1882b]{Poincar1882a}
\bgroup\fonteauteurs\bgroup Poincar{\'e}\egroup\egroup{}, H. (1882b).
\newblock Sur les formes cubiques ternaires et quaternaires. deuxi{\`e}me
  partie.
\newblock {\em Journal de l'{\'E}cole polytechnique},
  51\string:\penalty500\relax 45--91.

\bibitem[Poincar{\'e}, 1882c]{Poincar1882e}
\bgroup\fonteauteurs\bgroup Poincar{\'e}\egroup\egroup{}, H. (1882c).
\newblock Th{\'e}orie des groupes fuchsiens,.
\newblock {\em Acta Mathematica}, 1\string:\penalty500\relax 1--62.

\bibitem[Poincar{\'e}, 1883]{Poincar1883a}
\bgroup\fonteauteurs\bgroup Poincar{\'e}\egroup\egroup{}, H. (1883).
\newblock Sur les groupes des {\'e}quations lin{\'e}aires.
\newblock {\em Comptes rendus hebdomadaires des s{\'e}ances de l'Acad{\'e}mie
  des sciences}, 96\string:\penalty500\relax 691--694.

\bibitem[Poincar{\'e}, 1884a]{Poincar1884e}
\bgroup\fonteauteurs\bgroup Poincar{\'e}\egroup\egroup{}, H. (1884a).
\newblock Sur les groupes des {\'e}quations lin{\'e}aires.
\newblock {\em Acta Mathematica}, 4\string:\penalty500\relax 202--311.

\bibitem[Poincar{\'e}, 1884b]{Poincar1884d}
\bgroup\fonteauteurs\bgroup Poincar{\'e}\egroup\egroup{}, H. (1884b).
\newblock Sur les nombres complexes.
\newblock {\em Comptes rendus hebdomadaires des s{\'e}ances de l'Acad{\'e}mie
  des sciences}, 99\string:\penalty500\relax 740--742.

\bibitem[Poincar{\'e}, 1884c]{Poincar1884a}
\bgroup\fonteauteurs\bgroup Poincar{\'e}\egroup\egroup{}, H. (1884c).
\newblock Sur les substitutions lin{\'e}aires.
\newblock {\em Comptes rendus hebdomadaires des s{\'e}ances de l'Acad{\'e}mie
  des sciences}, 98\string:\penalty500\relax 349--352.

\bibitem[Poincar{\'e}, 1885]{Poincar1885b}
\bgroup\fonteauteurs\bgroup Poincar{\'e}\egroup\egroup{}, H. (1885).
\newblock M{\'e}moire sur les courbes d{\'e}finies par une {\'e}quation
  diff{\'e}rentielle (3e partie).
\newblock {\em Journal de Math{\'e}matiques}, (4e s{\'e}rie)
  1\string:\penalty500\relax 167--244.

\bibitem[Poincar{\'e}, 1886a]{Poincar1886h}
\bgroup\fonteauteurs\bgroup Poincar{\'e}\egroup\egroup{}, H. (1886a).
\newblock M{\'e}moire sur les courbes d{\'e}finies par une {\'e}quation
  diff{\'e}rentielle (4e partie).
\newblock {\em Journal de Math{\'e}matiques}, (4e s{\'e}rie)
  2\string:\penalty500\relax 151--217.

\bibitem[Poincar{\'e}, 1886b]{Poincar1886g}
\bgroup\fonteauteurs\bgroup Poincar{\'e}\egroup\egroup{}, H. (1886b).
\newblock Sur les int{\'e}grales irr{\'e}guli{\`e}res des {\'e}quations
  lin{\'e}aires.
\newblock {\em Acta Math{\'e}matica}, 8\string:\penalty500\relax 295--344.

\bibitem[Poincar{\'e}, 1887]{Poincar1887}
\bgroup\fonteauteurs\bgroup Poincar{\'e}\egroup\egroup{}, H. (1887).
\newblock Les fonctions fuchsiennes et l'arithm{\'e}tique.
\newblock {\em Journal de math{\'e}matiques pures et appliqu{\'e}es}, (4)
  3\string:\penalty500\relax 405--464.

\bibitem[Poincar{\'e}, 1889]{Poincar1889}
\bgroup\fonteauteurs\bgroup Poincar{\'e}\egroup\egroup{}, H. (1889).
\newblock Sur le probl{\`e}me des trois corps et les {\'e}quations de la
  dynamique avec des notes par l'auteur.
\newblock M{\'e}moire couronn{\'e} du prix de {S}.{M}. le {R}oi {O}scar II.

\bibitem[Poincar{\'e}, 1890]{Poincar1890}
\bgroup\fonteauteurs\bgroup Poincar{\'e}\egroup\egroup{}, H. (1890).
\newblock Sur le probl{\`e}me des trois corps et les {\'e}quations de la
  dynamique.
\newblock {\em Acta Mathematica}, 13\string:\penalty500\relax 1--270.

\bibitem[Poincar{\'e}, 1891]{Poincar1891}
\bgroup\fonteauteurs\bgroup Poincar{\'e}\egroup\egroup{}, H. (1891).
\newblock Le probl{\`e}me des trois corps.
\newblock {\em Revue g{\'e}n{\'e}rale des sciences pures et appliqu{\'e}es},
  2\string:\penalty500\relax 1--5.

\bibitem[Poincar{\'e}, 1892]{Poincar1892}
\bgroup\fonteauteurs\bgroup Poincar{\'e}\egroup\egroup{}, H. (1892).
\newblock {\em Les m{\'e}thodes nouvelles de la m{\'e}canique c{\'e}leste},
  volume tome 1.
\newblock Gauthier-Villars, Paris.

\bibitem[Poincar{\'e}, 1893]{Poincar1893}
\bgroup\fonteauteurs\bgroup Poincar{\'e}\egroup\egroup{}, H. (1893).
\newblock {\em Les m{\'e}thodes nouvelles de la m{\'e}canique c{\'e}leste},
  volume~2.
\newblock Gauthier-Villars, Paris.

\bibitem[Poincar{\'e}, 1899]{Poincar1899b}
\bgroup\fonteauteurs\bgroup Poincar{\'e}\egroup\egroup{}, H. (1899).
\newblock {\em Les m{\'e}thodes nouvelles de la m{\'e}canique c{\'e}leste},
  volume~3.
\newblock Gauthier-Villars, Paris.

\bibitem[Poincar{\'e}, 1900]{Poincar1900}
\bgroup\fonteauteurs\bgroup Poincar{\'e}\egroup\egroup{}, H. (1900).
\newblock Second compl{\'e}ment {\`a} l'analysis situs.
\newblock {\em Proceedings of the London Mathematical Society},
  32\string:\penalty500\relax 277--308.

\bibitem[Poincar{\'e}, 1901]{Poincar1901}
\bgroup\fonteauteurs\bgroup Poincar{\'e}\egroup\egroup{}, H. (1901).
\newblock Quelques remarques sur les groupes continus.
\newblock {\em Rendiconti del circolo matematico del Palermo},
  15\string:\penalty500\relax 321--368.

\bibitem[Poincar{\'e}, 1950]{PoincarOeuvres}
\bgroup\fonteauteurs\bgroup Poincar{\'e}\egroup\egroup{}, H. (1950).
\newblock {\em \OE uvres de Henri Poincar{\'e} publi{\'e}e sous les auspices de
  l'Acad{\'e}mie des sciences. Tome V, publi{\'e} avec la collaboration de M.
  Albert Ch{\^a}telet}.
\newblock Gauthier-Villars, Paris.

\bibitem[Poincar{\'e}, 2013]{Poincar2013}
\bgroup\fonteauteurs\bgroup Poincar{\'e}\egroup\egroup{}, H. (2013).
\newblock La correspondance de poincar{\'e} avec les astronomes et les
  g{\'e}od{\'e}siens.
\newblock \emph{In} \bgroup\fonteauteurs\bgroup Walter\egroup\egroup{}, S.,
  \bgroup\fonteauteurs\bgroup Kr{\"o}mer\egroup\egroup{}, R.,
  \bgroup\fonteauteurs\bgroup Nabonnand\egroup\egroup{}, P. et
  \bgroup\fonteauteurs\bgroup Schiavon\egroup\egroup{}, M., \'editeurs :  {\em
  La correspondance de Poincar{\'e} avec les astronomes et les
  g{\'e}od{\'e}siens}, volume~3. Birkh\"auser, B{\^a}le.

\bibitem[Poincar{\'e} et Mittag-Leffler, 1999]{Nabonnand1999}
\bgroup\fonteauteurs\bgroup Poincar{\'e}\egroup\egroup{}, H. et
  \bgroup\fonteauteurs\bgroup Mittag-Leffler\egroup\egroup{}, G. (1999).
\newblock {\em La correspondance entre Henri Poincar{\'e} and G{\"o}sta Mittag-
  Leffler, avec en annexe les lettres {\'e}chang{\'e}es par Poincar{\'e} avec
  Fredholm, Gyld{\`e}n et Phragm{\'e}n}.
\newblock Birkha{\"u}ser, Basel.

\bibitem[Poisson, 1809]{Poisson1809}
\bgroup\fonteauteurs\bgroup Poisson\egroup\egroup{}, S.~D. (1809).
\newblock M{\'e}moire sur les in{\'e}galit{\'e}s s{\'e}culaires des moyens
  mouvements des plan{\`e}tes, lu {\`a} l'{A}cad{\'e}mie le 20 juin 1808.
\newblock {\em Journal de l'{\'E}cole polytechnique}, Cahier
  XV(VIII)\string:\penalty500\relax 1--56.

\bibitem[Revel, 1996]{Revel1996}
\bgroup\fonteauteurs\bgroup Revel\egroup\egroup{}, J. (1996).
\newblock {\em Jeux d'echelles. La micro-analyse {\`a} l'exp{\'e}rience}.
\newblock Gallimard, Le Seuil, EHESS.

\bibitem[Robadey, 2004]{Robadey2004}
\bgroup\fonteauteurs\bgroup Robadey\egroup\egroup{}, A. (2004).
\newblock Exploration d'un mode d'{\'e}criture de la g{\'e}n{\'e}ralit{\'e}:
  l'article de poincar{\'e} sur les lignes g{\'e}od{\'e}siques des surfaces
  convexes (1905).
\newblock {\em Revue d'histoire des math{\'e}matiques},
  10(2)\string:\penalty500\relax 257--318.

\bibitem[Robadey, 2006]{Robadey2006}
\bgroup\fonteauteurs\bgroup Robadey\egroup\egroup{}, A. (2006).
\newblock {\em Diff{\'e}rentes modalit{\'e}s de travail sur le g{\'e}n{\'e}ral
  dans les recherches de {P}oincar{\'e} sur les syst{\`e}mes dynamiques}.
\newblock Th\`ese de doctorat, Universit{\'e} {P}aris7, Paris.

\bibitem[Roque, 2007]{Roque2007}
\bgroup\fonteauteurs\bgroup Roque\egroup\egroup{}, T. (2007).
\newblock Les enjeux du qualitatif dans la d{\'e}finition d'un syst{\`e}me
  dynamique.
\newblock \emph{In} \bgroup\fonteauteurs\bgroup Franceschelli\egroup\egroup{},
  S., \bgroup\fonteauteurs\bgroup Pary\egroup\egroup{}, M. et
  \bgroup\fonteauteurs\bgroup Roque\egroup\egroup{}, T., \'editeurs :  {\em
  Chaos et Syst{\`e}mes Dynamiques. {\'e}l{\'e}ments pour une
  {\'e}pist{\'e}mologie}, pages 37--66. Hermann, Paris.

\bibitem[Roque, 2011]{Roque2011}
\bgroup\fonteauteurs\bgroup Roque\egroup\egroup{}, T. (2011).
\newblock Stability of trajectories from {P}oincar{\'e} to {B}irkhoff :
  approaching a qualitative definition.
\newblock {\em Archive for History of Exact Sciences},
  65\string:\penalty500\relax 295--342.

\bibitem[Roque, 2008]{Roque2008}
\bgroup\fonteauteurs\bgroup Roque\egroup\egroup{}, T. (To appear (2008)).
\newblock The meaning of genericity in the classification of dynamical systems.
\newblock \emph{In} \bgroup\fonteauteurs\bgroup Chemla\egroup\egroup{}, K.,
  \bgroup\fonteauteurs\bgroup Cambefort\egroup\egroup{}, Y.,
  \bgroup\fonteauteurs\bgroup Chorlay\egroup\egroup{}, R. et
  \bgroup\fonteauteurs\bgroup Rabouin\egroup\egroup{}, D., \'editeurs :  {\em
  Perspectives on generality}. Oxford University Press.

\bibitem[Sapir, 1949]{Sapir1949}
\bgroup\fonteauteurs\bgroup Sapir\egroup\egroup{}, E. (1949).
\newblock Selected writings in language, culture and personality.
\newblock \emph{In} \bgroup\fonteauteurs\bgroup Mandelbaum\egroup\egroup{}, D.,
  \'editeur :  {\em Selected writings in language, culture and personality}.
  University of California Press, Berkeley.

\bibitem[Serret, 1849]{Serret1849}
\bgroup\fonteauteurs\bgroup Serret\egroup\egroup{}, J.-A. (1849).
\newblock {\em Cours d'alg{\`e}bre sup{\'e}rieure}.
\newblock Bachelier, Paris.

\bibitem[Shapin, 1982]{Shapin1982}
\bgroup\fonteauteurs\bgroup Shapin\egroup\egroup{}, S. (1982).
\newblock History of science and its sociological reconstructions.
\newblock {\em History of Science}, 20\string:\penalty500\relax 157--211.

\bibitem[Sinaceur, 1991]{Sinaceur1991}
\bgroup\fonteauteurs\bgroup Sinaceur\egroup\egroup{}, H. (1991).
\newblock {\em Corps et Mod{\`e}les, Essai sur l'histoire de l'alg{\`e}bre
  r{\'e}elle}.
\newblock Vrin, Paris.

\bibitem[Sinaceur, 1992]{Sinaceur1992}
\bgroup\fonteauteurs\bgroup Sinaceur\egroup\egroup{}, H. (1992).
\newblock Cauchy, sturm et les racines des {\'e}quations.
\newblock {\em Revue d'histoire des sciences}, 45(1)\string:\penalty500\relax
  51--68.

\bibitem[Smith, 1861]{Smith1861}
\bgroup\fonteauteurs\bgroup Smith\egroup\egroup{}, H. J.~S. (1861).
\newblock On systems of linear indeterminate equations and congruences.
\newblock {\em Philosophical Transactions of the Royal Society of London}, 151.

\bibitem[Sturm, 1829a]{Sturm1829a}
\bgroup\fonteauteurs\bgroup Sturm\egroup\egroup{}, C. (1829a).
\newblock Analyse d'un m{\'e}moire sur la r{\'e}solution des {\'e}quations
  num{\'e}riques.
\newblock {\em Bulletin de F{\'e}russac}, 11(271)\string:\penalty500\relax
  419--422.

\bibitem[Sturm, 1829b]{Sturm1829}
\bgroup\fonteauteurs\bgroup Sturm\egroup\egroup{}, C. (1829b).
\newblock Extrait d'un m{\'e}moire sur l'int{\'e}gration d'un syst{\`e}me
  d'{\'e}quations diff{\'e}rentielles lin{\'e}aires, pr{\'e}sent{\'e} {\`a}
  l'acad{\'e}mie des sciences le 27 juillet 1829 par m. sturm.
\newblock {\em Bulletin de F{\'e}russac}, 12\string:\penalty500\relax 313--322.

\bibitem[Sturm, 1835]{Sturm1835}
\bgroup\fonteauteurs\bgroup Sturm\egroup\egroup{}, C. (1835).
\newblock M{\'e}moire sur la r{\'e}solution des {\'e}quations num{\'e}riques.
\newblock {\em M{\'e}moires pr{\'e}sent{\'e}s par divers {S}avants
  {\'e}trangers {\`a} l'{A}cad{\'e}mie royale des sciences}, VI(273-318).

\bibitem[Sturm, 1836]{Sturm1836}
\bgroup\fonteauteurs\bgroup Sturm\egroup\egroup{}, C. (1836).
\newblock M{\'e}moire sur les {\'e}quations diff{\'e}rentielles lin{\'e}aires
  du second ordre.
\newblock {\em Journal de math{\'e}matiques pures et appliqu{\'e}es},
  1\string:\penalty500\relax 106--186.

\bibitem[Sturm et Liouville, 1836]{Sturm1836b}
\bgroup\fonteauteurs\bgroup Sturm\egroup\egroup{}, C. et
  \bgroup\fonteauteurs\bgroup Liouville\egroup\egroup{}, J. (1836).
\newblock D{\'e}monstration d'un th{\'e}or{\`e}me de cauchy.
\newblock {\em Journal de math{\'e}matiques pures et appliqu{\'e}es}, 278-289 ;
  290-308.

\bibitem[Sylvester, 1850]{Sylvester1850}
\bgroup\fonteauteurs\bgroup Sylvester\egroup\egroup{}, J.~J. (1850).
\newblock On the intersections, contacts, and other correlations of two conics
  expressed by indeterminante coordinates.
\newblock {\em Cambridge and Dublin Mathematical Journal},
  5\string:\penalty500\relax 119--137.

\bibitem[Sylvester, 1851]{Sylvester1851}
\bgroup\fonteauteurs\bgroup Sylvester\egroup\egroup{}, J.~J. (1851).
\newblock Enumeration of the contacts of lines and surfaces of the second
  order.
\newblock {\em Philosophical Magazine}, pages 119--140.

\bibitem[Sylvester, 1852]{Sylvester1852}
\bgroup\fonteauteurs\bgroup Sylvester\egroup\egroup{}, J.~J. (1852).
\newblock Sur une propri{\'e}t{\'e} nouvelle de l'{\'e}quation qui sert {\`a}
  d{\'e}terminer les in{\'e}galit{\'e}s s{\'e}culaires des plan{\`e}tes.
\newblock {\em Nouvelles Annales de Math{\'e}matiques}, pages 438--440.

\bibitem[Thomson et Tait, 1883]{ThomsonTait}
\bgroup\fonteauteurs\bgroup Thomson\egroup\egroup{}, W. et
  \bgroup\fonteauteurs\bgroup Tait\egroup\egroup{}, P. (1879--1883).
\newblock {\em Treatise on Natural Philosophy (2nd ed)}.
\newblock Cambridge University Press, Cambridge.

\bibitem[Tisserand, 1896]{Tisserand1889}
\bgroup\fonteauteurs\bgroup Tisserand\egroup\egroup{}, F. (1889--1896).
\newblock {\em Trait{\'e} de {M}{\'e}canique {C}{\'e}leste}.
\newblock Gauthier-Villars et fils, Paris.

\bibitem[Turner, 2012]{Turner2012}
\bgroup\fonteauteurs\bgroup Turner\egroup\egroup{}, L.~E. (2012).
\newblock The {M}ittag-{L}effler theorem: The origin, evolution and reception
  of a mathematical result, 1876--1884.
\newblock {\em Historia Mathematica}.

\bibitem[Weierstrass, 1858]{Weierstrass1858}
\bgroup\fonteauteurs\bgroup Weierstrass\egroup\egroup{}, K. (1858).
\newblock Ueber ein die homogenen functionen zweiten grades betreffendes
  theorem.
\newblock {\em Monatsberichte der K\"oniglich Preussischen Akademie der
  Wissenschaften zu Berlin}, pages 207--220.

\bibitem[Weierstrass, 1868]{Weierstrass1868}
\bgroup\fonteauteurs\bgroup Weierstrass\egroup\egroup{}, K. (1868).
\newblock Zur theorie der quadratischen und bilinearen formen.
\newblock {\em Monatsberichte der K{\"o}niglich Preussischen Akademie der
  Wissenschaften zu Berlin}, pages 310--338.

\bibitem[Wise, 2005]{Wise2005}
\bgroup\fonteauteurs\bgroup Wise\egroup\egroup{}, N. (2005).
\newblock William thomson and peter guthrie tait, treatise on natural
  philosophy, 1st edn (1867).
\newblock \emph{In} {\em In Grattan-Guinness and Cooke (2005)}, pages 521--533.

\bibitem[Yvon-Villarceau, 1870]{Villarceau1870}
\bgroup\fonteauteurs\bgroup Yvon-Villarceau\egroup\egroup{}, A. (1870).
\newblock Note sur les conditions des petites oscillations d'un corps solide de
  figure quelconque et la th{\'e}orie des {\'e}quations diff{\'e}rentielles
  lin{\'e}aires.
\newblock {\em Comptes rendus hebdomadaires des s{\'e}ances de l'Acad{\'e}mie
  des sciences}, 71\string:\penalty500\relax 762--766.

\end{thebibliography}

\bibliographystyle{apalike-fr}

\end{document}